\documentclass[twoside,11pt]{article}

\usepackage{jmlr2e_arxiv}

\usepackage[utf8]{inputenc} 
\usepackage[T1]{fontenc}    
\usepackage{hyperref}       
\usepackage{url}            
\usepackage{booktabs}       
\usepackage{amsfonts}       
\usepackage{nicefrac}       
\usepackage{microtype}      

\usepackage{amsfonts,amssymb,amsmath,amsbsy}
\usepackage{tabularx}
\usepackage{graphicx}
\usepackage{float}
\usepackage{verbatim}
\usepackage{dsfont}
\usepackage{placeins}
\usepackage{color}
\usepackage{tikz}
\usepackage{multimedia}
\usepackage{url}
\usepackage{bbm}
\usepackage{latexsym,mathrsfs}
\usetikzlibrary{arrows,shapes}
\usepackage{algorithm}
\usepackage{algpseudocode}
\usepackage{algorithmicx}
\usepackage{textcomp}
\usepackage{gensymb} 
\usepackage{placeins}

\def\R{\mathbb{R}}
\newcommand{\x}{\mathbf{x}}
\newcommand{\X}{\mathbf{X}}
\newcommand{\y}{\mathbf{y}}
\newcommand{\f}{\mathbf{f}}
\newcommand{\Y}{\mathbf{Y}}

\newcommand{\s}{\mathbf{s}}

\newcommand{\Rset}{\mathbb{R}}
\newcommand{\Xset}{\mathbb{X}}

\newcommand{\pfront}{\mathcal{P}^*}
\newcommand{\pset}{\mathbb{X}^*}
\newcommand{\yr}{\mathcal{Y}}
\newcommand{\fr}{\mathcal{F}}
\newcommand{\Fr}{\boldsymbol{\mathcal{F}}}
\newcommand{\Yr}{\boldsymbol{\mathcal{Y}}}
\newcommand{\Gr}{\boldsymbol{\mathcal{G}}}

\newcommand{\esp}{\mathbb{E}}
\newcommand{\prob}{\mathbb{P}}

\newcommand{\+}{_{n+1}}
\newcommand{\uto}{\mathbf{u}}
\newcommand{\dis}{\mathbf{d}}
\newcommand{\central}{\text{central}}
\newcommand{\nobj}{p}
\newcommand{\vecu}{\mathbf{u}}
\newcommand{\EIinv}{EI\textsubscript{nad}}
\newcommand{\PND}{p\textsubscript{ND}}

\def\cov{\operatorname{cov}}

\def\RR{\textsf{R}\/}

\def\R{\mathbb{R}}

\begin{document}

\title{The Kalai-Smorodinsky solution for many-objective Bayesian optimization}

%
%
%
\author{%
  \name M.~Binois
  \email mbinois@mcs.anl.gov\\ 
  \addr Mathematics and Computer Science Division, Argonne National Laboratory, 
  Lemont, USA
  \AND
  \name V.~Picheny
  \email victor@prowler.io\\
  \addr PROWLER.io, Cambridge, UK and\\
  \addr MIAT, Universit\'e de Toulouse, INRA,
  Castanet-Tolosan, France 
    \AND
  \name P.~Taillandier
  \email patrick.taillandier@inra.fr\\
  \addr MIAT, Universit\'e de Toulouse, INRA,
  Castanet-Tolosan, France 
  \AND
  \name A.~Habbal 
  \email habbal@unice.fr\\
  \addr Universit\'e C\^ote d'Azur, Inria, CNRS, LJAD, UMR 7351, 
  Parc Valrose, 06108 Nice, France}

\maketitle

\begin{abstract}  

An ongoing aim of research in multiobjective Bayesian optimization is to
extend its applicability to a large number of objectives.  While coping with a
limited budget of evaluations, recovering the set of optimal compromise
solutions generally requires numerous observations and is less interpretable
since this set tends to grow larger with the number of objectives. We thus
propose to focus on a specific solution originating from game theory, the
Kalai-Smorodinsky solution, which possesses attractive properties. In
particular, it ensures equal marginal gains over all objectives. We further
make it insensitive to a monotonic transformation of the objectives by
considering the objectives in the copula space.  A novel tailored algorithm is
proposed to search for the solution, in the form of a Bayesian optimization
algorithm: sequential sampling decisions are made based on acquisition
functions that derive from an instrumental Gaussian process prior. Our
approach is tested on four problems with respectively four, six, eight, and
nine objectives. The method is available in the \RR{} package \texttt{GPGame}
available on CRAN at \url{https://cran.r-project.org/package=GPGame}.
\end{abstract}

\begin{keywords}
  Gaussian process, Game theory, Stepwise uncertainty reduction
\end{keywords}

\section{Introduction}

Bayesian optimization (BO) is recognized as a powerful tool for global
optimization of expensive objective functions and now has a broad range of
applications in engineering and machine learning \citep[see, for
instance,][]{Shahriari2016}. A typical example is the calibration of a complex
numerical model:  to ensure that the model offers an accurate representation
of the system it emulates, some input parameters need to be chosen so that
some model outputs match real-life data \citep{walter1997identification}.
Another classical application is the optimization of performance of large
machine learning systems via the tuning of its hyperparameters
\citep{bergstra2011algorithms}.  In both cases, the high cost of evaluating
performance drastically limits the optimization budget,  and the high
sample-efficiency of BO compared with alternative black-box optimization
algorithms makes it highly competitive.

Many black-box problems, including those cited, involve several (or many)
performance metrics that are typically conflicting, so no common minimizer to
them exists. This is the field of multiobjective optimization, where one aims
at minimizing simultaneously a set of $\nobj$ objectives with respect to a set
of input variables over a bounded domain $\Xset \subset \mathbb{R}^d$:

\begin{equation}
 \min_{\x \in \Xset} \left\{ y^{(1)}(\x), \ldots, y^{(\nobj)}(\x) \right\}. 
\end{equation}
We assume that the $y^{(i)}: \Xset \rightarrow \mathbb{R}$ functions are
expensive to compute, nonlinear, and potentially observed in noise. Defining
that a point $\x^*$ dominates another point $\x$ if \textit{all} its
objectives are better, the usual goal of multiobjective optimization (MOO) is
to uncover the Pareto set, that is, the subset $\pset \subset \Xset$
containing all the Pareto nondominated solutions:
\begin{eqnarray}
 \forall \x^* \in \pset, \forall \x \in \mathbb{X}, \exists k \in \{1, \ldots, \nobj \} \text{ such that } \nonumber y^{(k)}(\x^*) \leq y^{(k)}(\x). \label{eq:defpset}
\end{eqnarray}
The image of the Pareto set in the objective space, $\pfront = \{
y^{(1)}(\pset), \ldots, y^{(\nobj)}(\pset) \}$, is called the Pareto front.
Since $\pset$ is in general not finite, most MOO algorithms aim at obtaining a
discrete representative subset of it.

Most of the well-established Pareto-based algorithms, such as evolutionary
\citep{deb2002fast,Chugh2017}, descent-based \citep{das1998normal}, or
Bayesian optimization \citep{wagner2010expected,Hernandez-Lobato2015}, perform
well on two or three objectives problems but poorly when $\nobj \ge 4$.
Indeed, the complexity of computing Pareto fronts dramatically increases when
the number of objectives increases to four and more
\citep{ishibuchi2008evolutionary}, such as the exponential increase in the
number of points necessary to approximate the --possibly singular-- Pareto
front hyper-surface and the difficulties in its graphical representation.
Moreover, one has to deal with a more MaO intrinsic problem, which is that
almost any admissible design becomes nondominated. In this paper, our aim is
precisely to address these ``many-objective optimization problems'' (MaO)
where ``many'' does refer to four or more objectives (that is, as soon as the
dominance relation looses its significance). Typical many-objective benchmark
problems have a number of objectives that range between 3 and 10
\citep{ishibuchi2016performance} or 3 and 15  \citep{cheng2017benchmark},
usually with far less complex (typically convex or linear) functions that the
ones we consider herein (test-cases with 4, 6, 8 and 9 nonlinear and/or non
convex objectives).

To circumvent these issues, \citet{Kukkonen2007} advocated the use of ranks
instead of nondomination and \citet{bader2011hype} contributions to the
hypervolume. However, such algorithms require many objective evaluations and
hence do not adapt well to expensive black boxes. In addition, they do not
solve the problem of exploiting the resulting very large Pareto set. Some
authors proposed methods to reduce the number of objectives to a manageable
one \citep{singh2011pareto}, to use the so-called decomposition-based
approaches \citep{asafuddoula2015decomposition}, or to rely on a set of fixed
and adaptive reference vectors \citep{Chugh2016}. Remarkably, while the
difficulty of representing the Pareto front is highlighted, the question of
selecting a particular solution on the Pareto front is mostly left to the
user.

An alternative way of solving this problem is by scalarisation, aggregating
all objective functions via weights, or by weighted metrics and goal
programming \citep{charnes1977goal}, or by the reference point methodology
\citep{wierzbicki1979use} which generalizes goal programming. This allows the
use of any technique dedicated to optimization of expensive black-boxes,
generally using surrogates \citep[see, e.g.,][]{Knowles2006,Zhang2010}.
Nevertheless, the aggregated function may become harder to model than its
components, and the relation between weights and the corresponding Pareto
optimal solution is generally not trivial. It may even be harder for
hyperparameter optimization tasks \citep{smithson2016neural} that have no
physical intuition backing weights. In addition, these issues worsen with more
objectives, advocating taking objectives separately.

Our present proposition amounts to searching for a single, but remarkable in
some sense, solution. To do so, we adopt a game-theoretic perspective, where
the selected solution arises as an equilibrium in a (non) cooperative game
played by $\nobj$ virtual agents who own the respective $\nobj$ objectives
\citep{JAD13-AIAA-Special}. In the following, we show that the so-called
Kalai-Smorodinsky (KS) solution \citep{KSE1975} is an appealing alternative.
We emphasize on the fact that the following development is built on the
assumption that all objectives have a priori equal importance. If any
preponderance was to be set as a rule in the theoretic game framework or as a
decision maker preference, then KS would be useless and ad hoc approaches
should be considered, see e.g., \cite{el2010optimisation} for a
fluid-structure application with preponderant aerodynamic criterion.

Intuitively, solutions at the ``center'' of the Pareto front are preferable
compared with those at extremities --which is precisely what the KS solution
consists of. Yet, the notion of center is arbitrary, since transforming the
objectives (nonlinearly, e.g., with a log scale) would modify the Pareto front
shape and affect the decision. Consider for instance the tuning of
hyperparameters of neural networks with accuracies, entropies and times as
objectives. Is a 0.05 improvement in entropy worth 0.05 second additional
prediction time, when most designs in the domain have entropies below 0.1, or
is some rescaling needed? Still, most MOO methods are sensitive to a rescaling
of the objective, which is not desirable \citep{Svenson2011}. Our second
proposition is thus to make the KS solution insensitive to monotone
transformations, by operating in the copula space
\citep{Nelsen2006,Binois2015b} with what we called the copula-KS (CKS)
solution. As we show in an experiment on hyperparameter tuning, this solution
provides a complementary alternative to the KS one.

Uncovering the KS and CKS solutions are nontrivial tasks for which, to our
knowledge, no algorithm is available in an expensive black-box framework. Our
third contribution is a novel Gaussian-process-based algorithm, building on
the stepwise uncertainty reduction (SUR) paradigm \citep{bect2012sequential}.
SUR, which is closely related to information-based approaches
\citep{Hennig2012,hernandez2016general}, has proven to be efficient for
solving single- and multiobjective optimization problems
\citep{villemonteix2009informational,picheny2013multi}, while enjoying strong
asymptotic properties \citep{bect2019supermartingale}.

The rest of the paper is organized as follows. Section \ref{sec:KS} describes
the KS solution and its extension in the copula space. Section \ref{sec:SUR}
presents the Bayesian optimization algorithm developed to find KS solutions.
Section \ref{sec:results} reports empirical performances of our
approach on three challenging problems with respectively four, six, eight, and nine
objectives. Section \ref{sec:conclusion} summarizes our conclusions and 
discusses areas for future work.

\section{The Kalai-Smorodinsky solution}\label{sec:KS}

\subsection{The standard KS solution} 
The Kalai-Smorodinsky solution was first proposed by Kalai and Smorodinsky in
1975 as an alternative to the Nash bargaining solution in cooperative
bargaining. The problem is as follows: starting from a {\em disagreement} or
{\em status quo} point $\dis$ in the objective space, the players aim at
maximizing their own benefit while moving from $\dis$ toward the Pareto front
(i.e., the efficiency set). The KS solution is of egalitarian inspiration
\citep{conley1991bargaining} and states that the selected efficient solution
should yield equal benefit ratio to all the players. Indeed, given the utopia
(or ideal, or shadow) point $\mathbf{u}\in\mathbb{R}^\nobj$ defined by
\[u^{(i)} = \min_{\x \in \mathbb{X}^*}\ y^{(i)}(\x),
\]
selecting any compromise solution $\s = [y^{(1)}(\x), \ldots, y^{(p)}(\x)]$
would yield, for objective $i$, a benefit ratio
\[r^{(i)}(\s)= \frac{d^{(i)} - s^{(i)}}{d^{(i)} - u^{(i)}}.\]

\begin{figure}[htbp]
\centering
 \includegraphics[trim=0mm 3mm 0mm 0mm, width=\textwidth, clip]{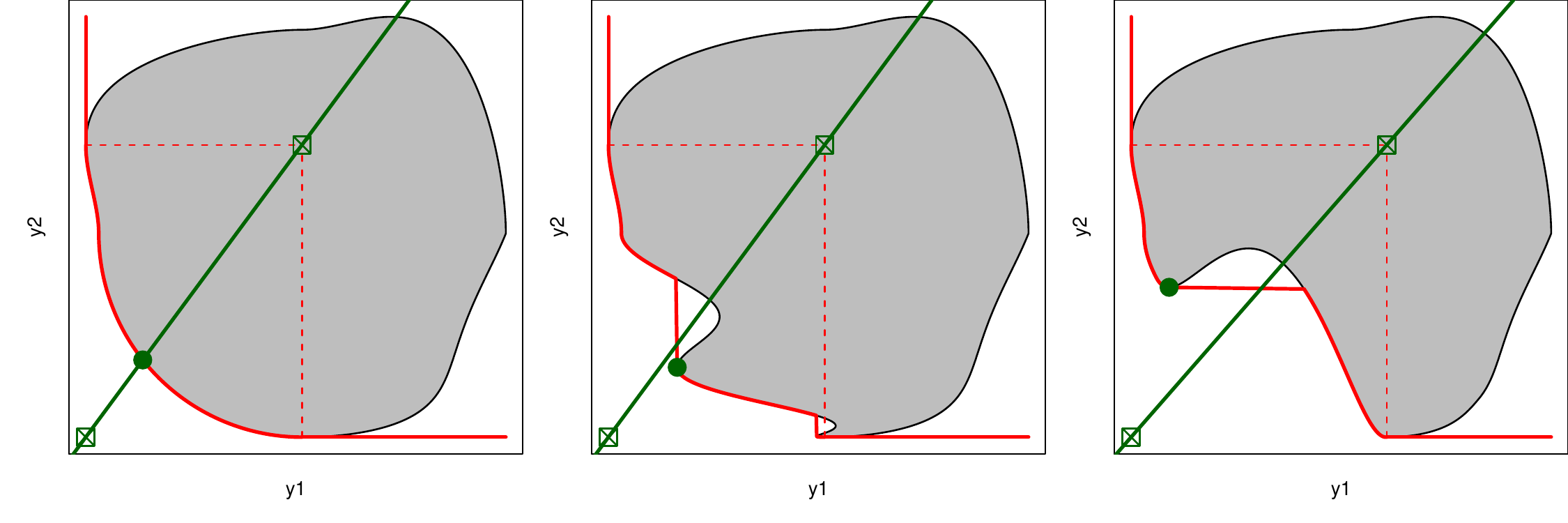}
 \caption{KS solution for a continuous (left) and discontinuous (center and
 right) $\pfront$. The KS is shown with a green disk, along with the utopia
 and nadir points (squares). The shaded area shows the feasible space, and the
 Pareto front is depicted in red.}
\label{fig:KSillustrated}
\end{figure}

Notice that the benefit from staying at $\dis$ is zero, while it is maximal
for the generically unfeasible point $\mathbf{u}$. The KS solution is the
Pareto optimal choice $\s^* = [y^{(1)}(\x^*), \ldots, y^{(p)}(\x^*)]$ for
which all the benefit ratios $r^{(i)}(\s)$ are equal. One can easily show that
$\s^*$ is the intersection point of the Pareto front and the line
$(\dis,\uto)$ (Figure \ref{fig:KSillustrated}, left).

We use here the extension of the KS solution to discontinuous fronts proposed
in \citet{hougaard2003nonconvex} under the name \textit{efficient maxmin
solution}. Indeed, for discontinuous fronts the intersection with the
$(\dis,\uto)$ line might not be feasible, so there is a necessary trade-off
between Pareto optimality and centrality. The efficient maxmin solution is
defined as the Pareto-optimal solution that maximizes the smallest benefit
ratio among players, that is:
\begin{equation}\label{eq:defKS}
  \s^{**} \in  \arg \max_{\y \in \pfront} \min_{1 \leq i \leq \nobj} r^{(i)}(\s).
\end{equation}
It is straightforward that when the intersection is feasible, then $\s^{*}$
and $\s^{**}$ coincide.

Figure \ref{fig:KSillustrated} shows $\s^{**}$ in two situations when the
feasible space is nonconvex. $\s^{**}$ is always on the Pareto front (hence
not on the $(\dis,\uto)$ line). In the central graph, it corresponds to the
point on the front closest to the $(\dis,\uto)$ line, which is intuitive. In
the right graph, the closest point on the front (on the right hand side of the
line) is actually a poor solution, as there exists another point with
identical performance in $y_2$ but much better in terms of $y_1$ (on the right
hand side of the line): the latter corresponds to $\s^{**}$.

In the following, we refer indifferently to $\s^{*}$ (if it exists) and
$\s^{**}$ as the KS solution. Note that this definition also extends to
discrete Pareto sets, which will prove useful in Section \ref{sec:SUR}.

For $\nobj \ge 3$, the KS solution, defined as the intersection point above,
fulfills some of the bargaining axioms: Pareto optimality, affine invariance,
and equity in benefit ratio. Moreover, for $\nobj=2$, KS is the {\em unique}
solution that fulfills all the bargaining solution axioms that are Pareto
optimality, symmetry, affine invariance, and restricted monotonicity
\citep{KSE1975}.

It is particularly attractive in a many-objective context since it scales
naturally to a large number of objectives and returns a single solution,
avoiding the difficulty of exploring and approximating large
$\nobj$-dimensional Pareto fronts--especially with a limited number of
observations.

The KS solution is known to depend strongly on the choice of the disagreement
point $\dis$. A standard choice is the nadir point $N$ given by $N_i =
\max_{\x \in \mathbb{X}^*}\ y^{(i)}(\x)$. Some authors introduced alternative
definitions, called extended KS, to alleviate the dependence on choice of
$\dis$ \citep{bozbay2012bargaining}, for instance by taking as disagreement
point the Nash equilibrium arising from a noncooperative game, but such a
choice would need a prebargaining split of the decision variable $\x$ among
the $\nobj$ players.

Closely related to the Kalai and Smorodinsky solution is the reference point
methodology, developed by \cite{wierzbicki1980use}. The reference point
methodology uses achievement functions, which refer to aspiration and
reservation references. Ideal and nadir points respectively could be taken as
such reference points, and achievement functions play the same role as the
benefit ratio of our present setting. Moreover, the notion of neutral
compromise and max-min (of achievement functions) were introduced by
Wierzbicki in the cited reference, yielding a framework objectively close to
the one of Kalai and Smorodinsky. One important difference between the
reference point and KS approaches is that the latter relies on a game
theoretic axiomatic construction, in view of palliating the Nash equilibrium
inefficiency.

\subsection{Incorporating preferences or constraints}
The egalitarian inspiration of the KS solution is based on the assumption that
all the objectives have equal importance. However, in many practical cases,
the end user may want to favor a subset of objectives which are of primary
importance, while still retaining the others for optimization. One way of
incorporating those \textit{preferences}
\citep{junker2004preference,thiele2009preference} is to discard solutions with
extreme values and actually solve the following problem:
\begin{eqnarray}
 \min_{\x \in \Xset} & \left\{ y^{(1)}(\x), \ldots, y^{(\nobj)}(\x) \right\}\nonumber \\
 \text{s.t.} & y^{(i)}(\x) \leq c_i, \qquad i \in J \subset \left[1, \ldots, \nobj \right],\label{eq:constraints}
\end{eqnarray}
with $c_i$'s predefined constants. Choosing a tight value (i.e., difficult to
attain) for $c_i$ may discard a large portion of the Pareto set and favor the
i-th objective over the others.

Incorporating such preferences in the KS solution can simply be done by using
$\mathbf{c}$ as the disagreement point if all objectives are constrained, or
by replacing the coordinates of the nadir with the $c_i$ values. In a game
theory context, this would mean that each player would state a limit of
acceptance for his objective before starting the cooperative bargaining.

From geometrical considerations, one may note that $\dis$ (and hence
$\mathbf{c}$) does not need to be a feasible point. As long as $\dis$ is
dominated by the utopia point $\uto$, the KS solution remains at the
intersection of the Pareto front and the $(\dis,\uto)$ line (optimistic $\dis$
may simply lead to negative ratios but Eq.\ \ref{eq:defKS} still applies). If
$\dis$ is not dominated by $\uto$, the $(\dis,\uto)$ direction is no more
relevant and the solution of Eq.\ \ref{eq:defKS} may not be Pareto-efficient.
A classical way to mitigate this issue is to allow some interactivity (so that
the players can reset their preferences to more realistic values). We refer to
\cite{hakanen2017using,tabatabaei2019interactive,miettinen2012nonlinear} for
further ideas on interactive preference-based algorithms.

\subsection{Robust KS using copulas}

A drawback of KS is that it is not invariant under a monotonic (nonaffine)
transformation of the objectives (this is not the case of the Pareto set,
since a monotonic transformation preserves ranks, hence domination relations).
In a less-cooperative framework, some players could be tempted to rely on such
transformations to influence the choice of a point on the Pareto front. It may
even be involuntary, for instance when comparing objectives of different
natures.

To circumvent this problem, we use copula theory, which has been linked to
Pareto optimality by \citet{Binois2015b}. In short, from a statistical point
of view, the Pareto front is related to the zero level-line $\partial F_0$ of
the multivariate cumulative density function $F_Y$ of the objective vector $Y
=
\left(y_1(X), \dots, y_p(X) \right)$ with $X$ any random vector with support equal to
$\Xset$. That is, $\partial F_0 = \lim_{\alpha \rightarrow 0^+} \left\{\y \in
\R^\nobj, F_Y(\y) = \prob(Y_1 \leq y_1, \dots, Y_p \leq y_p) = \alpha\right\}$. 
Notice that solving the MaO problem by random
sampling (with $X$ uniformly distributed), as used for hyperparameter
optimization by \citet{Bergstra2012}, amounts to sample from $F_Y$. 
Extreme-level lines of $F_Y$ actually indicate how likely it is to be close to the
Pareto front (or on it if the Pareto set has a nonzero probability mass).

A fundamental tool for studying multivariate random variables is the copula.
A copula $C_Y$ is a function linking a multivariate cumulative distribution
function (CDF) to its univariate counterparts, such that for $\y \in
\R^\nobj$, $F_Y(\y) = C_Y(F_1(y_1),
\dots, F_p(y_p))$ with $F_i = \prob(Y_i \leq y_i)$, $1 \leq i \leq \nobj$
\citep{Nelsen2006}. When $F_Y$ is continuous, $C_Y$ is unique, from Sklar's
theorem \citep[Theorem 2.3.3]{Nelsen2006}.  Consequently, as shown by
\citet{Binois2015b}, learning the marginal CDFs $F_i$'s and extreme levels of
the copula $C_Y$, that is $\partial C_0 = \lim_{\alpha \rightarrow 0^+}
\left\{\vecu \in [0,1]^\nobj, C_Y(\vecu) = F_Y(F_1^{-1}(u_1), \dots,
F_\nobj^{-1}(u_\nobj)) = \alpha \right\}$, amounts to learning $\partial F_0$,
hence the Pareto front.

A remarkable property of copulas is their invariance under monotone increasing
transformations of univariate CDFs \citep[Theorem 2.4.3]{Nelsen2006}.
Extending their proof to the $\nobj$-dimensional case, suppose that $g_1,
\dots, g_p$ are strictly increasing transformations on the ranges of $Y_1,
\dots, Y_p$, respectively; hence they are invertible. Denote $G_1, \dots, G_p$
the marginal distribution functions of $\mathbf{g}(Y)$. It then holds that
$G_i(y_i) = \prob(g_i(Y_i) \leq y_i) = \prob(Y_i \leq g_i^{-1}(y_i)) =
F_i(g_i^{-1}(y_i))$. Then $C_{\mathbf{g}(Y)}(G_1(y_1), \dots, G_\nobj(y_p)) =
\prob(g_1(Y_1) \leq y_1,
\dots, g_p(Y_p) < y_p) = \prob(Y_1 \leq g_1^{-1}(y_1), \dots, Y_p \leq
g_p^{-1}(y_p)) = C_Y(F_1(g_1^{-1}(y_1)), \dots, F_p(g_p^{-1}(y_p))) =
C_{Y}(G_1(y_1), \dots, G_p(y_p))$.  

Now, our proposition is to consider the KS solution in the copula space, that
is, taking $F_1, \ldots, F_p$ as objectives instead of  $y_1, \ldots, y_p$.
This ``copula-KS solution'' (henceforth CKS) is Pareto-efficient and invariant
to any monotonic transformation of the objectives. In addition, in the copula
space the utopia point is always $(0, \ldots, 0)$. While the nadir remains
unknown, the point $(1, \ldots, 1)$ may serve as an alternative disagreement
point, since it corresponds to the worst solution for each objective. This
removes the difficult task of learning the $(\dis,\uto)$ line, at the expense
of learning the marginal distribution and the copula function. With this
choice, one may remark that when $\Xset$ is finite, our proposed solution is
very intuitive, as it amounts to choosing the maxmin solution over ranks, as:
\begin{equation}\label{eq:defCKS}
    r_C^{(i)}(\s) = \frac{1}{\text{Card}(\Xset)}\sum_{j=1}^{\text{Card}(\Xset)} \delta \left(s^{(i)} \leq y^{(i)}(\x_j) \right),
\end{equation}
with $\delta()$ the Kronecker delta function.

Importantly, CKS now depends on the instrumental law of $X$. For instance,
this is the probability distribution that would be used to solve the problem
via random sampling, such as a uniform distribution for a bounded domain and a
normal one for unbounded domain. Without loss of generality, we can always
assume that $X$ is uniformly distributed over $\Xset$, and using a different
distribution with the same support is equivalent to introducing a
transformation over $X$. Hence, contrary to KS, CKS is sensitive to the input
definition (e.g., the MOO problems $\min(y^{(1)}(x), y^{(2)}(x))$ and
$\min(y^{(1)}(x^2), y^{(2)}(x^2))$ with $x \in [0, 1]$, with equal Pareto
fronts, would lead to different CKS solutions). Despite this, in a MaO
context, CKS remains much more robust to reformulations than KS, in the sense
that transformations in the input space will affect all of the objectives and
not a single one. This makes it much more difficult to bias the solution in
favor of, say, a particular objective: it can only be done by increasing the
mass of the distribution of $X$ around the optimum of this objective, which is
presumably unknown for expensive black-box objectives.

In the following, we only consider the CKS computed using the uniform
distribution over $\Xset$. By doing so, we assume that the problem formulation
is appropriate, in the sense that the outputs are roughly stationary with
respect to the inputs (as opposed to an ill-defined problem, that would
contain large plateaus for instance).

\subsection{Illustration}
\begin{figure}[htpb]
 \includegraphics[trim=0mm 2mm 0mm 0mm, width=.32\textwidth, clip]{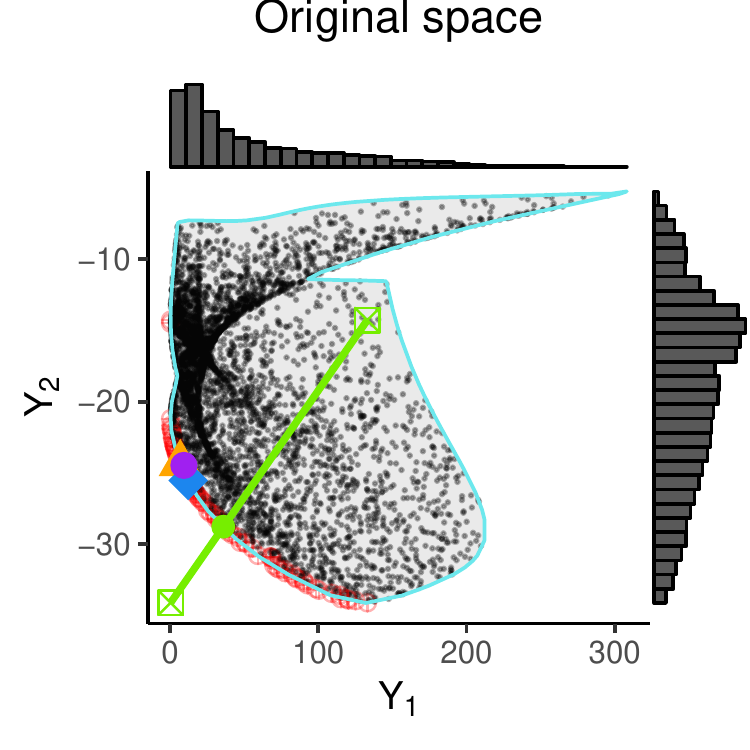}
 \includegraphics[trim=0mm 2mm 0mm 0mm, width=.32\textwidth, clip]{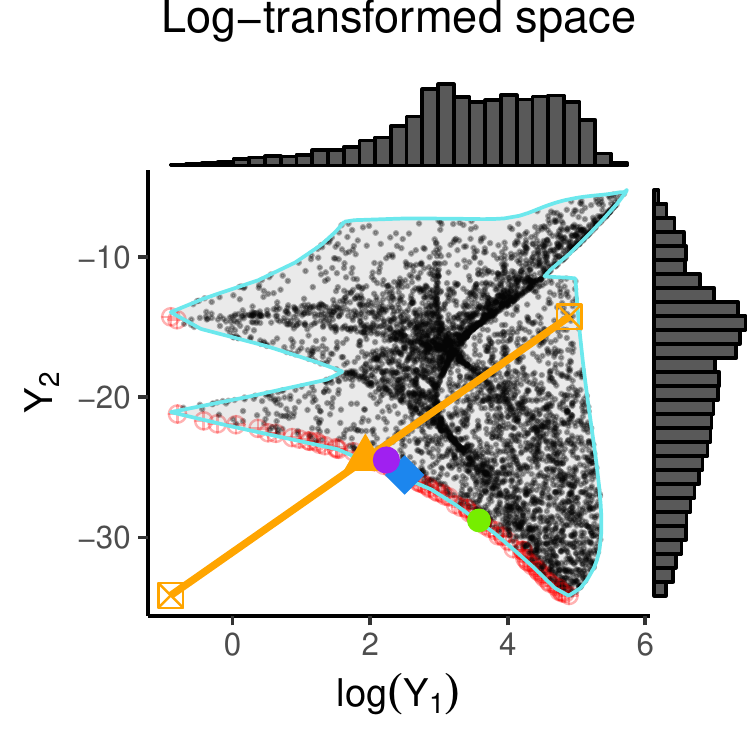}
 \includegraphics[trim=0mm 2mm 0mm 0mm, width=.32\textwidth, clip]{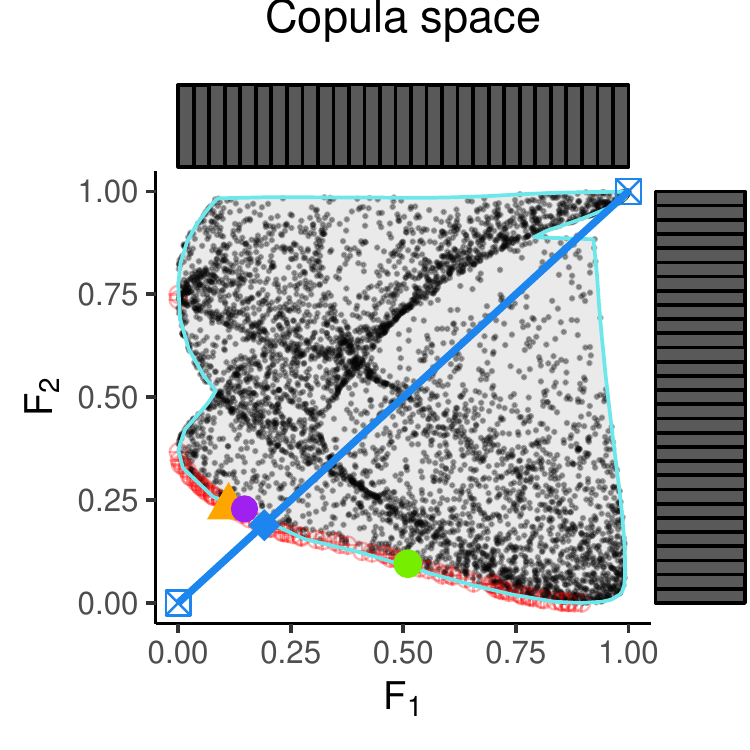}\\
 \includegraphics[trim=0mm 2mm 0mm 0mm, width=.32\textwidth, clip]{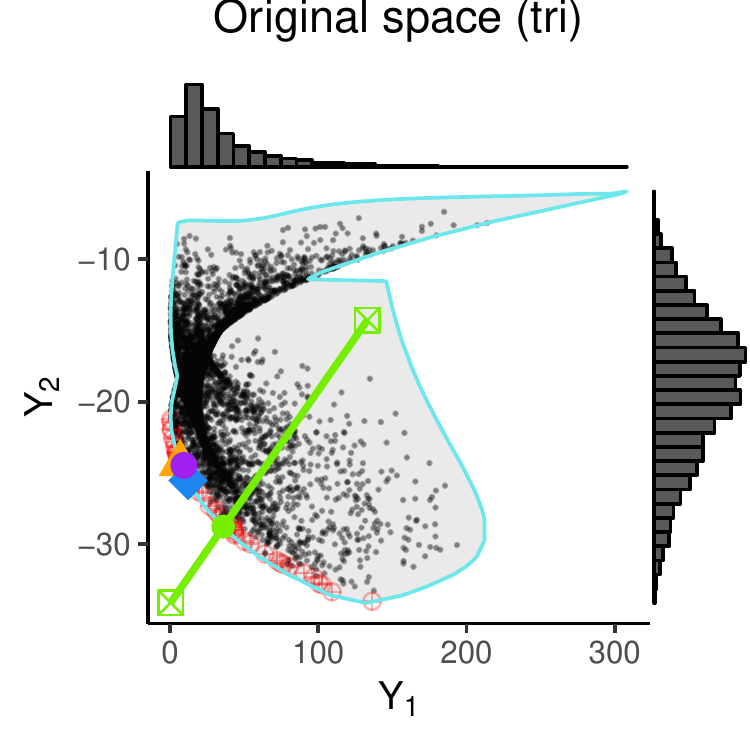}
 \includegraphics[trim=0mm 2mm 0mm 0mm, width=.32\textwidth, clip]{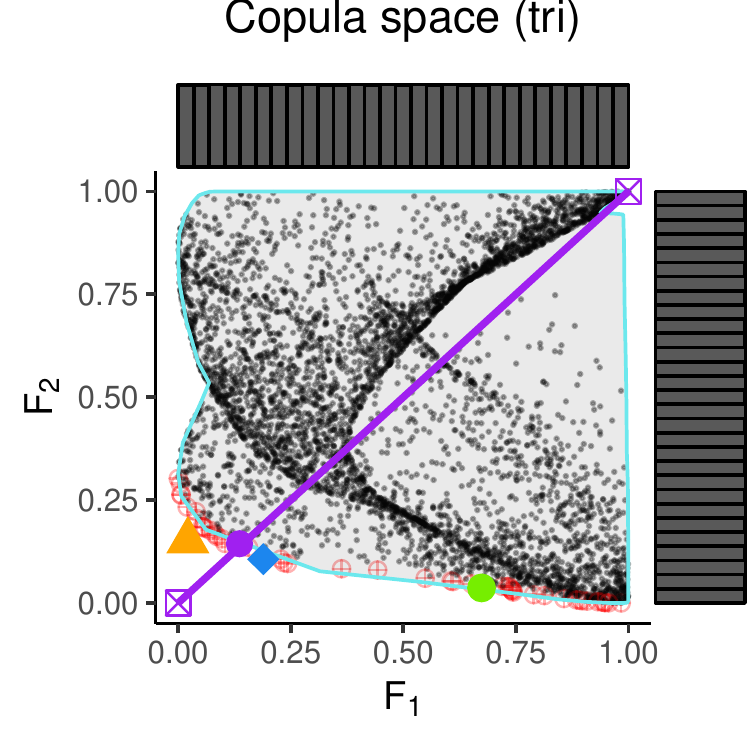}
 \includegraphics[trim=0mm 2mm 0mm 0mm, width=.32\textwidth, clip]{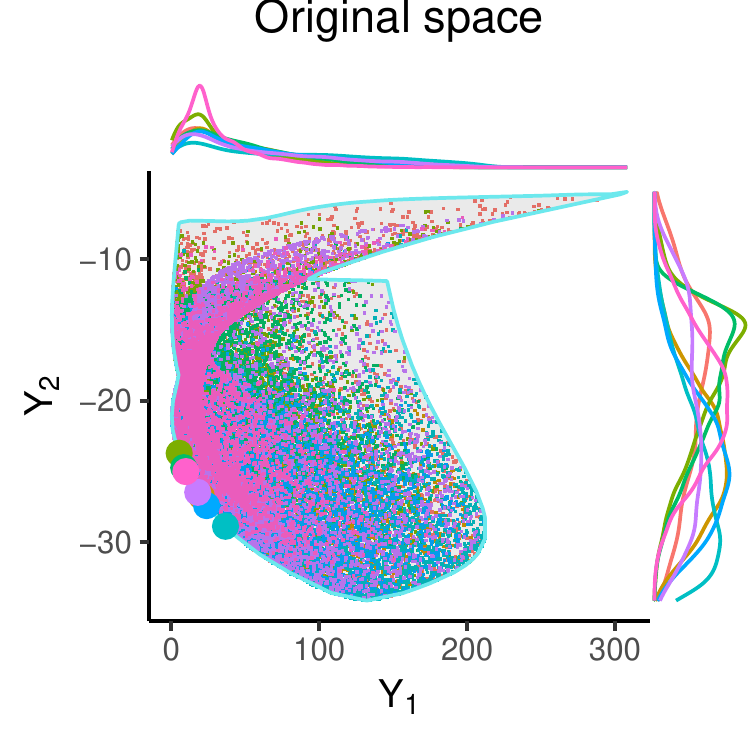}
  \caption{Top: KS (green disk), KS in log-scale (orange triangle) and CKS
  (blue diamond) solutions for a biobjective problem, based on $5,000$
  uniformly sampled designs in $\Xset=[0,1]^2$, shown in the objective space.
  The black dots show all the dominated values of the grid and the red crossed
  circles the Pareto-optimal ones. The shaded area shows the feasible space.
  Marginal objective densities are reported on the corresponding axes. The
  $(\dis,\uto)$ lines are shown with matching colors. Bottom left: same
  problem, but with $5,000$ samples in $\Xset$ from a triangular distribution.
  The corresponding CKS solution (bottom, center) is shown in purple. Bottom
  right: CKS solutions corresponding to 8 different input distributions.}
\label{fig:KSandCKS}
\end{figure}

Let us consider a classical two-variable, biobjective problem from the
multiobjective literature \citep[P1, see][]{parr2013improvement}. We first
compute the two objectives on $5,000$ uniformly sampled points in
$\Xset=[0,1]^2$, out of which the feasible space, Pareto front, and KS
solution are extracted (Figure \ref{fig:KSandCKS}, top left). The KS solution
has a visibly central position in the Pareto front, which makes it a
well-balanced compromise.

Applying a log transformation of the first objective does not change the
Pareto set but modifies here substantially the shape of the Pareto front (from
convex to concave) and the KS solution, leading to a different compromise
(Figure \ref{fig:KSandCKS}, top center), which compared to the previous one
favors $Y_1$ over $Y_2$.

Both original and rescaled problems share the same image in the copula space
(Figure \ref{fig:KSandCKS}, top right), which provides a third compromise.
Seen from the original space, the CKS solution seems here to favor the first
objective: this is due to the high density of points close to the minimum on
$Y_1$ (Figure \ref{fig:KSandCKS}, top left, top histogram). It is, however,
almost equivalent to the KS solution under a log transformation of $Y_1$. From
a game perspective, the two players agree on a solution with equal ranks: here
roughly the $1,200$th best ($F_1 \approx F_2 \approx 0.24$) out of $5,000$,
independently of the gains in terms of objective values.

Now, we consider the same problem, but the $5,000$ points are sampled in
$\Xset$ using a triangular distribution for each marginal (between 0 and 1,
asymmetric with a mode at $x=0.25$). Using a distribution does not change the
KS solution, but modifies the distributions of the objectives (Figure
\ref{fig:KSandCKS}, bottom left, both histograms). The image in the copula
space (Figure \ref{fig:KSandCKS}, bottom center) is thus different, which
leads to a slightly different CKS solution.

Finally, we compare CKS solutions obtained by using independent beta
distributions on the inputs, with parameters randomly chosen between 0.5 and 4
(Figure \ref{fig:KSandCKS}, bottom right). Despite being based on extremely
different distributions, the CKS solutions are relatively similar.

\section{A Bayesian optimization algorithm to find KS solutions}\label{sec:SUR}

Computing the KS or CKS solutions are challenging problems. It requires for KS
learning the ideal and disagreement points $\uto$ and $\dis$ \citep[which are
very difficult problems on their own; see, for
instance,][]{bechikh2010estimating} and for CKS the marginals and copula, as
well as (for both) the part of the Pareto front that intersects the
$(\dis,\uto)$ line. An additional difficulty arises when the objective
functions cannot be computed exactly but only through a noisy process.

In this section, we consider that one has access to observations of the form
\begin{equation}
 \f_i = \y(\x_i) + \boldsymbol{\varepsilon}_i,
\end{equation}
where $\boldsymbol{\varepsilon}_i$ is a zero-mean noise vector with
independent components of variances $(\tau_i^{(1)}, \ldots,
\tau_i^{(\nobj)})$. Our approach readily adapts to the deterministic case by
setting $\forall j: \tau_i^{(j)}=0$. We assume here that all objectives are
collected at the same time (using a single black-box), hence sharing the
experimental points $\x_i$. The case of several black-boxes, as presented by
\citet{hernandez2016general}, is not considered here and is deferred to future
work.

\subsection{Elements of Bayesian optimization}
For our algorithm, we consider a classical BO framework, where independent
Gaussian process (GP) priors are put on the objectives:
\begin{equation}\label{eq:Y}
  \forall i \in 1, \ldots, \nobj, \quad Y^{(i)}(.) \sim \mathcal{GP} \left( {\mu^{(i)}}(.),{\sigma^{(i)}}\left(.,.\right) \right),
\end{equation}
where the mean $\mu^{(i)}$ and covariance $\sigma^{(i)}$ have predetermined
parametric forms whose parameters are estimated by maximum likelihood
\citep{Rasmussen2006}. Conditioning on a set of observations $\left\{ \f_1,
\ldots, \f_n \right\}$, GPs provide flexible response fits associated with
uncertainty estimates. They enable operating sequential design decisions via
an \textit{acquisition function} $J(\x)$, which balances between exploration
and exploitation in seeking global optima. Hence, the design consists of a
first set of $n_0$ observations generated by using a space-filling design,
generally from a variant of Latin hypercube design
\citep[LHD,][]{mckay1979comparison} to obtain a first predictive distribution
of $Y^{(i)}(.)$, and a second set of sequential observations chosen as
\begin{equation}\label{eq:boloop}
 \x_{n+1} \in \arg \max_{\x \in \Xset} J(\x), \quad (n \geq n_0).
\end{equation}

GP equations are deferred to Appendix A. In the following, we use the
subscript $n$ to denote quantities conditional on the set of $n$ observations
(e.g., $Y_n^{(i)}$, $\mu_n^{(i)}$ or $\sigma_n^{(i)}$).

\subsection{Baseline algorithms}\label{sec:baseline}

We propose first a simple algorithm for KS, that alternatively searches for
the utopia point, the nadir, the KS solution, and so on until the budget is
met.

Looking for the utopia point is the simplest task, as it amounts to looking
for the individual minima of the objective functions. We can use $\nobj$ times
an acquisition function designed for single objective problems, e.g., expected
improvement, EI, \cite{jones1998efficient} in the noiseless setting, or one of
its ``noisy'' counterpart, see \cite{picheny2013benchmark} for a review.

The nadir point is constructed from worst objective function values of Pareto
optimal solutions. For $\nobj=2$, looking for the utopia directly provides the
nadir. However, when $\nobj\geq3$, it does not coincide with individual
objective extrema. Accurately estimating the nadir is a difficult task,
especially when the number of objective functions increases, for which only a
handful of algorithms are available
\citep{deb2006towards,bechikh2010estimating}. Noting that each coordinate is
then the solution of a single objective constrained problem (maximum under
constraint that $\x$ belongs to the Pareto front), we propose to use the
constrained version of EI \citep[Expected Feasible Improvement,
EFI,][]{schonlau1998global}. Our acquisition function is the expected
improvement of minus the objective (\EIinv, to point towards maximal values
instead of minimal ones) multiplied by the probability of being nondominated
by the observation set (\PND, which is available in closed form from the GP
marginal \citep{couckuyt2014fast}). See Appendix A for expressions of EI,
\EIinv~ and \PND.

Finally, given an $\uto$ and $\dis$, the KS solution can be sought by using a
GP-UCB \citep{srinivas2012information} analog for Eq.\ \ref{eq:defKS}, that
is, replacing the $r^{(i)}(\s)$ values by a lower confidence bound provided by
the GP model, which is simply the GP posterior mean minus a positive constant
$\beta$ (typically, $\beta=1.96$) times the GP posterior standard deviation:
\begin{equation*}
    r_{LCB}^{(i)}(y(\x))= \frac{d^{(i)} - \mu_n^{(i)}(\x) + \beta \sqrt{\sigma_n^{(i)}(\x, \x)}}{d^{(i)} - u^{(i)}}, \quad \beta \geq 0.
\end{equation*}

In total, for $\nobj$ objectives we have $2p+1$ different acquisition
functions. Some corresponding tasks may be more important than others, but it
is difficult to prioritize them in a principled way. Hence, as a baseline
algorithm, we propose to optimize and sample the corresponding maximizer of
each acquisition function sequentially (see Algorithm \ref{alg:baseline}).

\begin{algorithm}[!ht]
\caption{Pseudo-code for the baseline algorithm for KS}\label{alg:baseline}
\begin{algorithmic}[1]
\Require  $\nobj$ GP models trained on $\left\{ \x_1, \ldots, \x_n \right\}$, $\left\{ \f_1, \ldots, \f_n \right\}$
\While{$n \leq N$}
    \For{$i \gets 1$ to $\nobj$}
        \State Choose $\x_{n+1} = \arg \max_{\s \in \Xset} EI^{(i)}(\x)$
        \State Compute $\f_{n+1}$ using the expensive black-box
        \State Update the GP models by conditioning on $\{\x_{n+1}, \f_{n+1}\}$.
        \State $n \leftarrow n+1$
    \EndFor
    \For{$i \gets 1$ to $\nobj$}
        \State Choose $\x_{n+1} = \arg \max_{\x \in \Xset} \EIinv^{(i)}(\x) \times \PND(\x)$
        \State Compute $\f_{n+1}$ using the expensive black-box
        \State Update the GP models by conditioning on $\{\x_{n+1}, \f_{n+1}\}$.
        \State $n \leftarrow n+1$
    \EndFor
    \State Choose $\x_{n+1}$ solution of Eq.\ \ref{eq:defKS} with $r_{LCB}^{(i)}$.
    \State Compute $\f_{n+1}$ using the expensive black-box
    \State Update the GP models by conditioning on $\{\x_{n+1}, \f_{n+1}\}$.
    \State $n \leftarrow n+1$
\EndWhile
\end{algorithmic}
\end{algorithm}

Finding the CKS solution requires learning the copula and marginals, which
does not easily convert into acquisition functions. However, both quantities
correspond to global features of the objectives. As a baseline, we propose to
sample where the GP models are most uncertain, that is, using as acquisition
functions the posterior variances $\sigma_n^{(i)}$. Using the same structure
as Algorithm \ref{alg:baseline}), we propose $2p$ maximizations of the
$\sigma_n^{(i)}$'s, followed by an exploitation step. Instead of using the
LCB, which is not available in closed form for CKS, we propose here to sample
at the CKS solution of the GP posterior mean.

\subsection{Stepwise uncertainty reduction} \label{sec:sur}

The drawback of the presented baseline algorithms are that they do not balance
the different learning tasks in a principled way, which may result in very
subefficient sampling as the number of objectives increases. Instead, we
propose here to design a single acquisition $J(\x)$ tailored to our problem.
To do so, we follow a step-wise uncertainty reduction (SUR) approach, similar
to the one proposed by \citet{picheny2016bayesian} to solve Nash equilibria
problems.

Let us denote by $\Psi: \mathbb{Y} \rightarrow \Rset^\nobj$ the mapping that
associates a KS or CKS solution with any multivariate function $\y$ defined
over $\Xset$. Given a distribution $\Y_n()$, $\Psi(\Y_n)$ is a random vector
(of unknown distribution). The first step of designing a SUR strategy is to
define an uncertainty measure $\Gamma(\Y_n)$ that reflects our lack of
knowledge of our problem solution (\ref{eq:defKS}). A simple measure of
variability of a vector is the determinant of its covariance matrix
\citep{fedorov1972theory}:
\begin{equation}\label{eq:Gamma}
 \Gamma(\Y_n) = \det \left[\cov \left(\Psi(\Y_n) \right) \right].
\end{equation}
Intuitively, $\Gamma(\Y_n)$ tends to zero when all the components of
$\Psi(\Y_n)$ become known accurately.

The SUR strategy greedily chooses the next observation that reduces the most
this uncertainty:
\begin{equation*}
 \max_{\x \in \Xset} \Gamma(\Y_n) - \Gamma(\Y_{n, \x}),
\end{equation*}
where $\Y_{n, \x}$ is the GP conditioned on $\{\y(\x_1), \ldots, \y(\x_n),
\y(\x)\}$. However, such an ideal strategy is not tractable as it would
require evaluating all $\y(\x)$ while maximizing over $\Xset$.

For a tractable strategy, we consider the expected uncertainty reduction:
\begin{equation}
  \Gamma(\Y_n) - \esp_{\Y_n(\x)} \left[ \Gamma(\Y_{n, \x}) \right],
\end{equation}
where $\esp_{\Y_n(\x)}$ denotes the expectation taken over $\Y_n(\x)$.
Removing the constant term $\Gamma(\Y_n)$, our policy is defined with
\begin{equation}\label{eq:nextpoint}
  \x\+ \in \arg \min_{\x \in \Xset} J(\x) = \esp_{\Y_n(\x)} \left[ \Gamma(\Y_{n, \x}) \right].
\end{equation}

Note that $J(\x)$ defines a trade-off between exploration and exploitation, as
well as a trade-off between the different learning tasks ($\dis$ and $\uto$
points, copula, Pareto front).

\subsection{Computational aspects} \label{sec:comput}
Given a distribution $\Y_n()$, solving exactly Eq.\ \ref{eq:nextpoint} is
impossible for two reasons. First, the mapping $\Psi$ can be very complex and
there are no algorithms to approach it accurately. Second, the nonlinearity of
$\Gamma$ prevent us from getting a closed-form expression for $J(\x)$.

We propose to solve both problems by replacing the (continuous) design space
$\Xset$ by a discrete representation $\Xset^\star$ of size $N$, and evaluate
our acquisition function by Monte-Carlo. The two key-points for our approach
to work in practice are, first, that the cost of the Monte-Carlo evaluation
can be substantially alleviated by specific pre-computations and update
formulae: this is exposed in Section \ref{sec:computingsur}. Second, the
Monte-Carlo approach has a cubic cost which only allows the use of small
discrete sets, but carefully choosing those sets allows us to estimate our
acquisition function with precision: this is proposed in Section
\ref{sec:integpts}.

\subsubsection{Computing and optimizing the SUR
criterion}\label{sec:computingsur} Let $\Xset^\star = \{\x_1^\star, \ldots,
\x_N^\star\} \subset \Xset$. Since $\Xset^\star$ is discrete, one can easily
generate independent drawings of $\Y_n\left(\Xset^\star \right)$:
$\Yr_1^\star, \ldots, \Yr_M^\star$. For each sample, the corresponding KS
solution $\Psi(\Yr_i^\star)$ can be computed by exhaustive search. The
following empirical estimator of $\Gamma(\Y_n)$ is then available:
\begin{equation*}
 \hat \Gamma \left( \Yr_1^\star, \ldots, \Yr_M^\star \right) = \det \left[ \mathbf{Q}_{\yr}\right],
\end{equation*}
with $\mathbf{Q}_{\yr}$ the sample covariance of $\Psi(\Yr_1^\star), \ldots, \Psi(\Yr_M^\star)$.

Now, let $\x$ be a candidate observation point, and $\Fr_1, \ldots, \Fr_K$ be
independent drawings of $\Y_n(\x)$. Denoting $\Yr_1^\star | \Fr_i, \ldots,
\Yr_M^\star | \Fr_i$ independent drawings of $\Y_{n, \x}(\Xset^\star) =
\Y_n(\Xset^\star) | \Y_n(\x)=\Fr_i$, an estimator of $J(\x)$ is then obtained
by using the empirical mean:
\begin{equation*}
 \hat J(\x) = \frac{1}{K}\sum_{i=1}^K \hat \Gamma \left( \Yr_1^\star|\Fr_i, \ldots, \Yr_M^\star|\Fr_i \right).
\end{equation*}

\begin{figure}[htbp]
\centering
 \includegraphics[trim=0mm 0mm 0mm 13mm, width=\textwidth, clip]{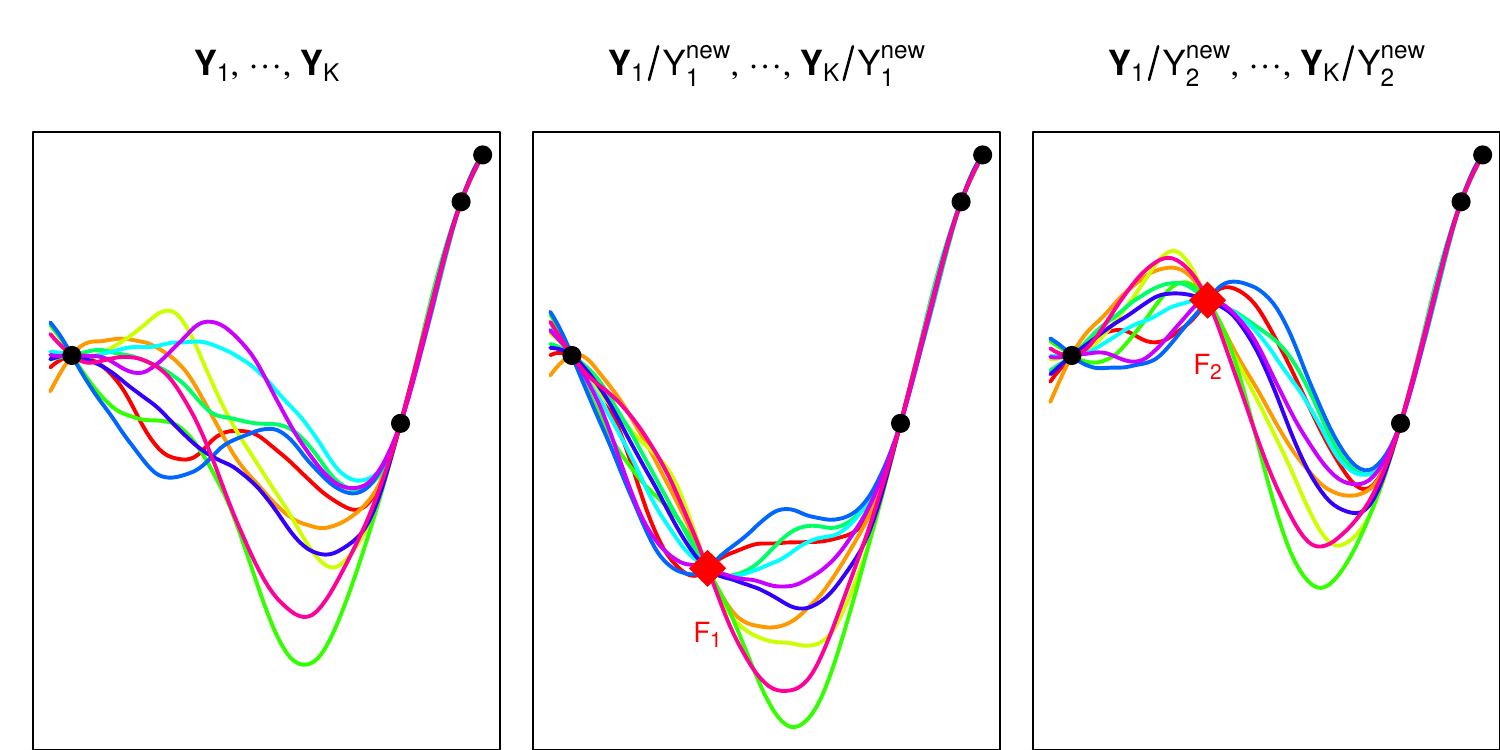}
 \caption{Left: initial drawings $\Yr_1^\star, \ldots, \Yr_M^\star$. Middle
 and right: same drawings, but conditioned further on an observation $\fr^1$
 (middle) and $\fr^2$ (right).}
\label{fig:GPS}
\end{figure}

Hence, computing $\hat J(\x)$ requires $M \times K$ samples
$\Yr_j^\star|\Fr_i$. This can simply be done by nesting two Monte-Carlo loops,
but would be very intensive computationally. To overcome this problem, we
employ here the \textit{fast update of conditional simulation ensemble}
algorithm proposed by \citet{chevalier2015fast}, as detailed below and
illustrated in Figure \ref{fig:GPS}. As shown by \citet{chevalier2015fast},
samples $\Yr^\star | \Fr$ of $\Y_{n, \x}(\Xset^\star)$ (conditioned on $n+1$
observations) can be obtained efficiently from a joint sample $\{\Yr^\star,
\Fr\}$ of $\Y_{n}(\{\Xset^\star, \x\})$ (conditioned only on $n$ observations)
by working on residuals using, with $1 \leq i \leq \nobj$:
\begin{eqnarray}\label{eq:updatesim}
 \yr^{\star(i)} | \fr^{(i)} &=& \yr^{\star(i)} +  \boldsymbol{\lambda}^{(i)}(\x) \left( \fr^{(i)} - \yr^{(i)}(\x) \right),
\end{eqnarray}
with $\Yr^\star = [\yr^{\star (1)}, \ldots, \yr^{\star (p)}] \in \Rset^{N \times \nobj}$, $\Yr(\x) = [\yr^{(1)}(\x), \ldots, \yr^{(p)}(\x)]  \in \Rset^{\nobj}$,\\ $\Fr = [\fr^{\star (1)}, \ldots, \fr^{\star (p)}]  \in \Rset^{\nobj}$ 
and 
\begin{eqnarray*}
\boldsymbol{\lambda}^{(i)}(\x) = \frac{1}{\boldsymbol{\sigma}_n^{(i)}(\x, \x)} \left[\boldsymbol{\sigma}_n^{(i)}(\x^\star_1, \x), \ldots, \boldsymbol{\sigma}_n^{(i)}(\x^\star_N, \x) \right].
\end{eqnarray*}
Hence, the two expensive operations are computing
$\boldsymbol{\lambda}^{(i)}(\x)$ and drawing samples of $\Y_{n}(\{\Xset^\star,
\x\})$, the later having an $\mathcal{O}((N+1)^3)$ complexity when using the
standard decomposition procedure based on Cholesky, \citep[see,
e.g.,][]{Diggle2007}.

A first decisive advantage of the update approach is that computing
$\boldsymbol{\lambda}^{(i)}(\x)$ needs to be done only once for each $\x$,
regardless of the values taken by $\yr^{\star(i)}$ and $\fr^{(i)}$. This
considerably eases the nested Monte-Carlo loop.

Moreover, if we restrict $\x$ to belong to $\Xset^\star$, drawing
$\{\Yr^\star, \Fr\}$ reduces to drawing $\Yr^\star$ and needs to be done only
once, independently of $\x$, hence prior to minimizing the acquisition
function. Hence, although optimizing $\hat J(\x)$ over the continuous space
$\Xset$ is definitely feasible, it requires completing the simulation over
$\Xset^\star$ with $\x$, incurring an additional computational cost. We thus
restrict $\x$ to $\Xset^\star$, and solve Eq.\ \ref{eq:nextpoint} by
exhaustive search.

\subsubsection{Choosing integration points \texorpdfstring{$\Xset^\star$}{X*}}\label{sec:integpts}

Due to the cubic complexity, generating $\Yr_1^\star, \ldots, \Yr_M^\star$
limits $\Xset^\star$ in practice to at most a couple thousand points. In small
dimension (say, $1 \leq d \leq 3$), $\Xset^\star$ may consist simply of a
Cartesian grid or a dense space-filling design \citep{niederreiter1988low}.
For larger $d$, such approaches do not scale and $\Xset^\star$ may not be
large enough to cover $\Xset$ sufficiently. However, accurate approximations
of $J(\x)$ can still be obtained, as long as $\Xset^\star$ covers the
influential parts of $\Xset$ with respect to $J$, that is, regions where (a)
the KS or CKS solution is likely to lie, (b) (for KS) minimal values of $\Y$
are likely (to estimate the utopia), (c) (for KS, without preferences) maximal
non-dominated values of $\Y$ are likely (to estimate the nadir).

To design a set $\Xset^\star$ of limited size that contains those three
components, we proceed as follows. First, a large space-filling set
$\Xset_{\text{large}}$ is generated, and the marginal posterior distribution
is computed for each element of $\Xset_{\text{large}}$. Since we do not
consider the joint distribution, this has a cost of only
$\mathcal{O}(N_\text{large}^2)$. Then, we use this marginal distribution to
pick a small subset of $\Xset_{\text{large}}$ that is likely to contain the
components a, b, and c. For robustness, the choice of the subset corresponds
to an exploration (large $\Y_n$ variance) / exploitation (following the $\Y_n$
mean) trade-off, as we show below.

\paragraph{Central part}
For the first BO iteration, we obtain a rough estimate of the solution by
computing the KS or CKS solution of the posterior GP mean. As the mean is
likely to be smoother than the actual function, a better solution, available
as soon as one iteration has been performed, is to use the set of simulated
KS/CKS solutions generated during the previous computation of $J$:
$\Psi(\Yr_1^\star), \ldots, \Psi(\Yr_M^\star)$. Then, using the marginal
distribution, we compute the probability $p_{\text{box}}$ for each element of
$\Xset_{\text{large}}$ to belong to a box defined by the extremal values of
the simulated KS (see Appendix A for formulas). A first set $\Xset_\central$
is obtained by sampling from $\Xset_{\text{large}}$ randomly with
probabilities proportional to $p_{\text{box}}$.

\paragraph{Utopia}
To include points with an exploration / exploitation trade-off, we compute for
each objective the expected improvement for all elements in
$\Xset_{\text{large}}$. The $\nobj$ EI maximizers constitute a second set
$\Xset_\text{utopia}$.

\paragraph{Nadir} Following the baseline strategy (Section
\ref{sec:baseline}), to capture the nadir we compute for each objective the
expected improvement of minus the objective multiplied by the probability of
non-domination. The $\nobj$ maximizers of these criteria constitute a third
set  $\Xset_\text{nadir}$.
 
Finally, the reduced set $\Xset^\star$ is taken as 
$\Xset_\central \cup \Xset_\text{utopia} \cup \Xset_\text{nadir}$ for KS, 
$\Xset_\central \cup \Xset_\text{utopia}$ for KS with a predefined disagreement point and
$\Xset_\central$ for CKS.

$\Xset^\star$ may be renewed every time a new observation is added to the GP
model, firstly to account for the changes in the GP distribution, but also to
improve robustness (a critical region can be missed at an iteration but
accounted for during the next).

Algorithm \ref{alg:surloop} provides the pseudo-code of the SUR loop for KS,
that replaces the standard acquisition function maximisation of Eq.\
\ref{eq:boloop} in the classical BO structure. The CKS algorithm is almost
identical, except that $\Xset^\star = \Xset_\central$, and the copula requires
a specific treatment, as we describe next.

\begin{algorithm}[!ht]
\caption{Pseudo-code for the SUR loop for KS}\label{alg:surloop}
\begin{algorithmic}[1]
\Require  $\nobj$ GP models trained on $\left\{ \x_1, \ldots, \x_n \right\}$, $\left\{ \f_1, \ldots, \f_n \right\}$
\Require $\Xset_{\text{large}}$, $N_\star$, $M$
\State Compute the marginal distribution of $\Y(\Xset_{\text{large}})$
\State Using the marginal compute: $p_{\text{box}}$, EI, \EIinv, \PND
\State Choose: 
\begin{itemize}
    \item $\Xset_\central \subset \Xset_{\text{large}}$ with $p_{\text{box}}$, size $N_\star - 2p$
    \item $\Xset_\text{utopia} \subset \Xset_{\text{large}}$ with $EI$
    \item $\Xset_\text{nadir} \subset \Xset_{\text{large}}$ with $\EIinv \times \PND$
\end{itemize}
\State Define $\Xset^\star = \Xset_\central \cup \Xset_\text{utopia} \cup \Xset_\text{nadir}$
\State Generate $M$ draws of $\Y(\Xset^\star)$: $(\Yr_1^\star, \ldots, \Yr_M^\star)$  
\For{$i \gets 1$ to $N_\star$}
    \State From $(\Yr_1^\star, \ldots, \Yr_M^\star)$, extract $M$ draws of $\Y(\x_\star^i)$: $\fr_1, \ldots, \fr_M$
    \For{$k \gets 1$ to $M$}
        \State Compute $(\Yr_1| \fr^k, \ldots, \Yr_M| \fr^k)$ using Eq.\ (\ref{eq:updatesim})
        \State Compute $\Gamma_{ik} = \hat \Gamma \left( \Yr_1|\fr^k, \ldots, \Yr_M|\fr^k \right)$
    \EndFor
    \State Compute $\hat J_i = \frac{1}{M}\sum_{k=1}^M \Gamma_{ik}$
\EndFor
\State Choose $\x_{n+1} = \x^{\star}_I$ with $I = \arg \min_{i \in [1, N_\star]} \hat J_{i}$
\State Compute $\f_{n+1}$ using the expensive black-box
\State Update the GP models by conditioning on $\{\x_{n+1}, \f_{n+1}\}$.
\end{algorithmic}
\label{alg:SUR}
\end{algorithm}

\subsubsection{Copula estimation}\label{sec:copula}

In the continuous case, unless additional information is available about the
marginal distributions and copula function, empirical estimators may be used;
see, for instance, \citet{Omelka2009}. Note that as we are interested in CKS,
only a good approximation of the diagonal of the copula is critical, which
makes it an easier problem than regular copula inference. Additionally, the
CKS may not be in the tails of the marginal distributions, which also eases
the inference.

Nonetheless, empirical estimates are based on a sample of $\Y(X)$, with $X$
i.i.d. according to the instrumental distribution (i.e., uniform). Hence, both
the observation set and $\Xset^\star$ are inappropriate, as they are too small
sets, and more importantly, not uniformly distributed. Instead, the copula and
marginals are estimated on a large auxiliary i.i.d. sample
$\{\Xset_{\text{aux}}, \Y_{\text{aux}}\}$. One may choose $\Xset_{\text{aux}}
= \Xset_{\text{large}}$. Since jointly sampling $\Y(\Xset_{\text{aux}})$ is
out of reach for a large sample (Section \ref{sec:integpts}), we follow
\citet{oakley2004estimating} and use the conditional simulations $\Yr_j$ on
$\Xset^\star$ as pseudo-observations to update the GP predictive mean, and we
take our sample $\Y_{\text{aux}}$ for this updated mean $\esp \left[
Y(\Xset_{\text{aux}}) | \Yr_j \right]$ \citep[this approach has also been
referred to as \textit{hallucinated
observations},][]{desautels2014parallelizing}.

With that setting, in practice there is no need to explicitly estimate the
copula and marginals. The CKS is simply given by the ranks with respect to the
auxiliary set, i.e.,
\begin{equation*}
  r_C^{(i)}(\s) = \frac{1}{N_{\text{aux}}} \sum_{j=1}^{N_{\text{aux}}} \delta \left(s^{(i)} \leq y_\text{aux}^{(i)}(\x_j) \right),
\end{equation*}
with $N_{\text{aux}} = \text{Card}(\Xset_{\text{aux}})$.

\section{Experiments}\label{sec:results}

This section details numerical experiments on four test problems: two toy
problems from the multiobjective and BO literature, a problem of
hyperparameter tuning, and the calibration of a complex simulator. All
experiments were conducted in \texttt{R} \citep{R2016}, by using the dedicated
package \texttt{GPGame}; see the work of \citet{picheny2017} for details.

\subsection{Synthetic problems}
As a proof of concept, we consider the DTLZ2 function \citep{deb2002fast},
with five variables and four objectives, that has a concave dome-shaped Pareto
front, and a six variables and six objectives problem obtained by rotating and
rescaling six times the classical mono-objective function hartman
\citep{dixon1978towards} (see Appendix D for the problem formulation).

The objectives of these experiments are 1) to compare the performance of SUR
with baseline strategies, and 2) to evaluate the effect of discretizations. On
both cases, we consider two finite domains $\Xset$ with respectively $10^5$
and $10^7$ elements uniformly distributed in $[0,1]^d$. This allows the
computation of the exact KS and CKS solutions.

For the smaller domain, we set $\Xset_{\text{large}} = \Xset$, while for the
larger one we take $\Xset_{\text{large}}$ as a random subset of size $10^5$ of
$\Xset$ (renewed every iteration), in order to emulate the loss incurred by
SUR for working with discrete representations of continuous spaces.

For SUR, we test two configurations corresponding respectively to a small or
large computational effort to estimate the acquisition function. The
parameters are reported in Table \ref{tab:configs}, along with the
computational times required for a single SUR iteration on ten CPUs @2.3Ghz.

As competitors, we use the baseline algorithms presented in Section
\ref{sec:baseline} and a purely exploratory algorithm based on uniform
sampling on $\Xset$. In addition, we ran SMS-EGO
\citep{ponweiser2008multiobjective}, a multi-objective Bayesian optimization
algorithm based on hypervolume improvement that has been shown to compare
favorably to other indicator-based algorithms
\citep{wagner2010expected,binois2019gpareto}, in particular on the DTLZ2
problem.

\begin{table}[ht]\label{tab:configs}
\caption{Configurations and wall clock times (for a single iteration on
hartman) of the SUR algorithm. CKS relies on $10^4$ auxiliary points.}
\centering
\begin{tabular}{ccccc}
Configuration & $N$ & $M$ & clock time (KS) & clock time (CKS) \\
Coarse & 250 & 25 & 33s & 114s\\
Fine & 1000 & 100 & 528s & 3420s\\
\end{tabular}

\end{table}

All BO strategies (baselines and SUR) start with $n_0 = 2 \times d$
observations generated from an optimized LHD, followed by respectively $60$
(for DTLZ2) and $78$ (for hartman) infill points. Each strategy is run 10
times.

\begin{figure}[htbp]
\centering
\includegraphics[trim=0mm 0mm 0mm 0mm, width=0.16\textwidth, clip]{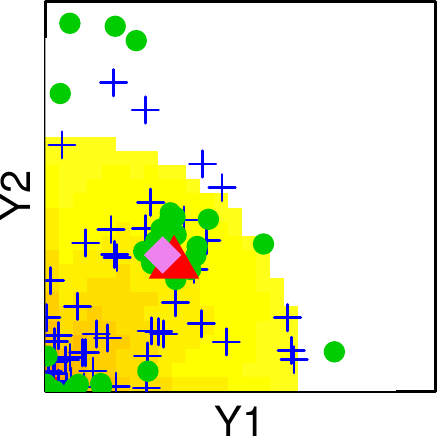}
\includegraphics[trim=0mm 0mm 0mm 0mm, width=0.16\textwidth, clip]{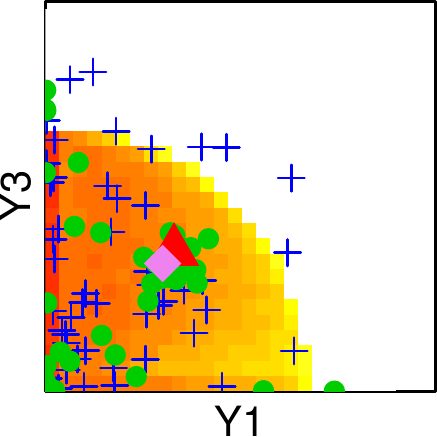}
\includegraphics[trim=0mm 0mm 0mm 0mm, width=0.16\textwidth, clip]{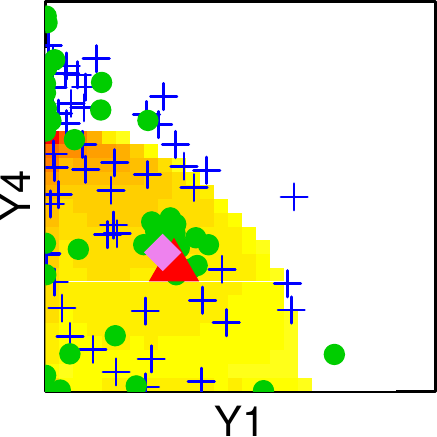}
\includegraphics[trim=0mm 0mm 0mm 0mm, width=0.16\textwidth, clip]{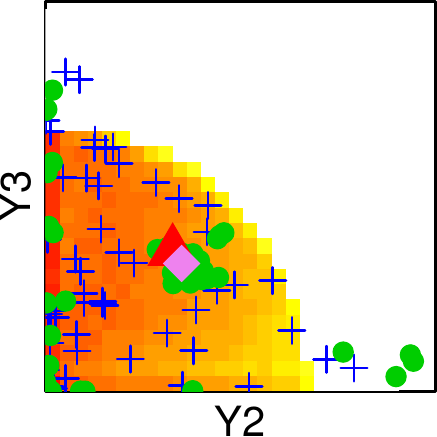}
\includegraphics[trim=0mm 0mm 0mm 0mm, width=0.16\textwidth, clip]{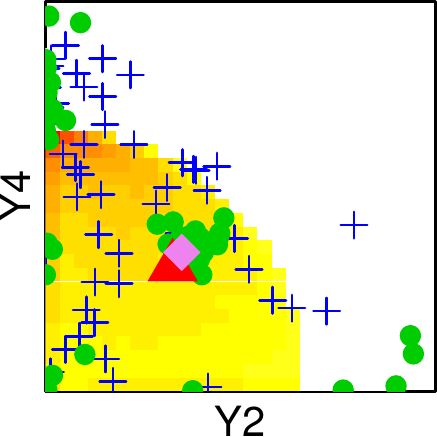}
\includegraphics[trim=0mm 0mm 0mm 0mm, width=0.16\textwidth, clip]{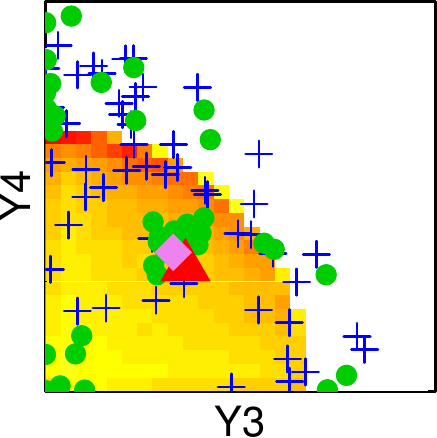}\\
\includegraphics[trim=0mm 0mm 0mm 0mm, width=0.16\textwidth, clip]{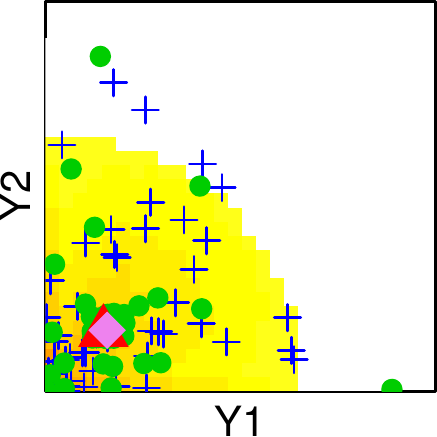}
\includegraphics[trim=0mm 0mm 0mm 0mm, width=0.16\textwidth, clip]{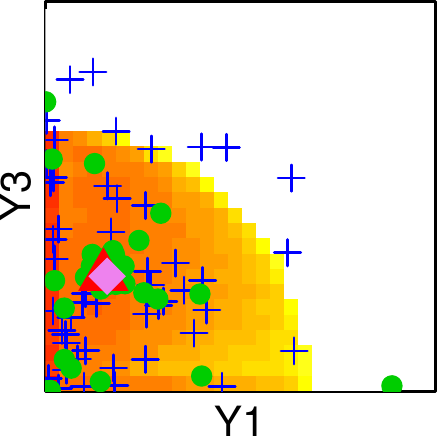}
\includegraphics[trim=0mm 0mm 0mm 0mm, width=0.16\textwidth, clip]{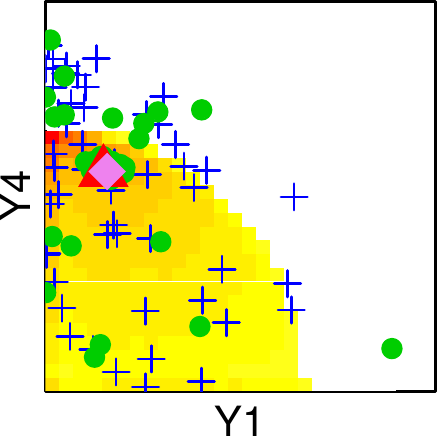}
\includegraphics[trim=0mm 0mm 0mm 0mm, width=0.16\textwidth, clip]{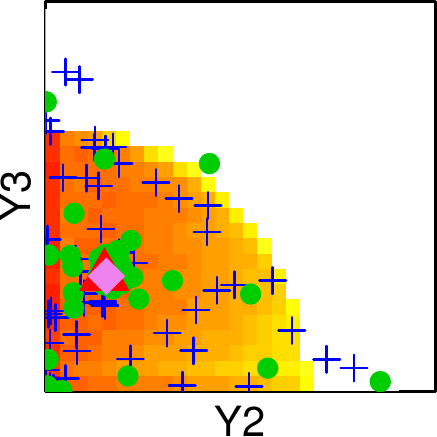}
\includegraphics[trim=0mm 0mm 0mm 0mm, width=0.16\textwidth, clip]{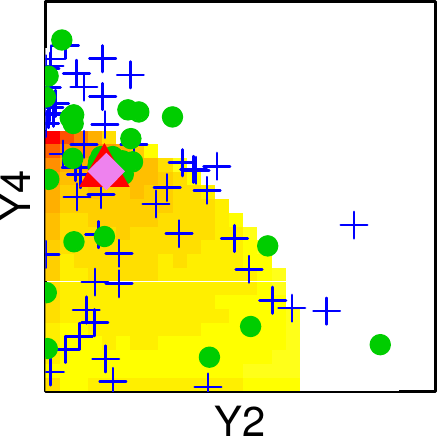}
\includegraphics[trim=0mm 0mm 0mm 0mm, width=0.16\textwidth, clip]{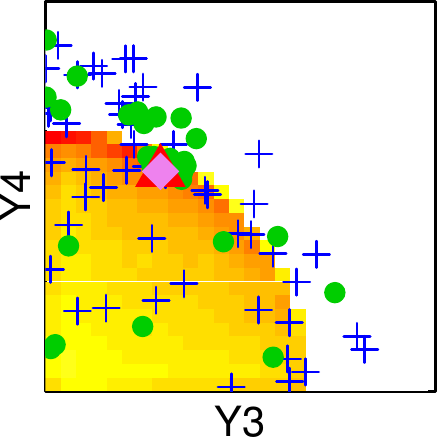}
 \caption{Results of one run for KS (top) and CKS (bottom), represented by
 using marginal 2D projections. The blue crosses are the initial design
 points, the green points the added designs, the red triangles the target true
 equilibria and pink diamonds the predicted ones. The heatmap represents the
 density of Pareto front points on the projected spaces.}
\label{fig:resDTLZ2}
\end{figure}

Results for one run of SUR are given in Figure \ref{fig:resDTLZ2} in the form
of projections on the marginal 2D spaces. Notice the central location of the
KS point on this problem (Figure \ref{fig:resDTLZ2}, first row), while the CKS
point leans toward areas of larger densities (for instance, second line, third
and fifth plots, CKS is close to the upper left corner). For KS, new points
are added close to the reference solution, but some are also more exploratory,
in particular near the individual minima to reduce uncertainty on the
$(\dis,\uto)$ line. For CKS, the behavior is more local, with points added
mostly around the reference solution.

\begin{figure}[htbp]
\centering
\includegraphics[trim=0mm 0mm 32mm 0mm, height=70mm, clip]{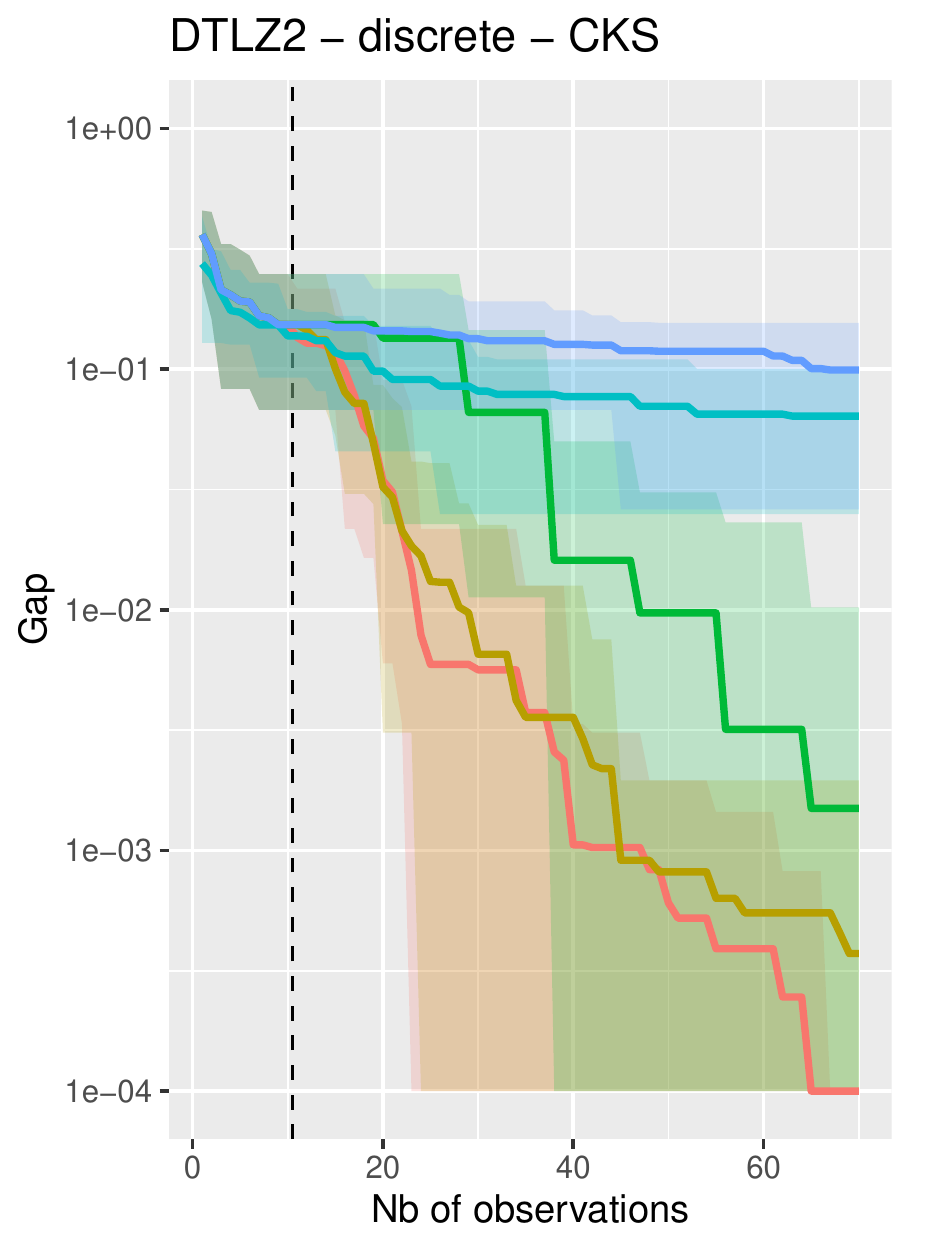}%
\includegraphics[trim=0mm 0mm 32mm 0mm, height=70mm, clip]{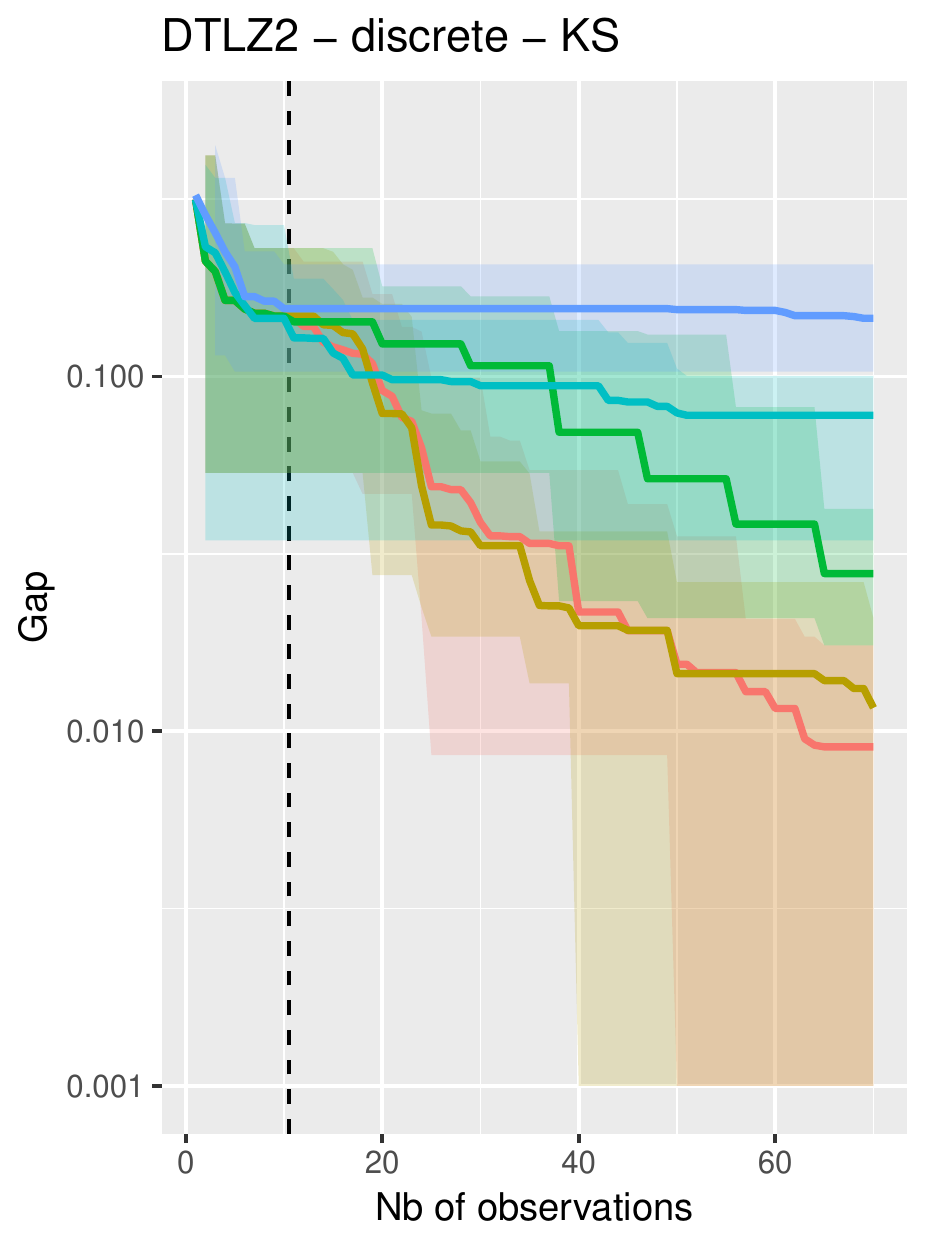}%
\includegraphics[trim=15mm 0mm 32mm 0mm, height=70mm, clip]{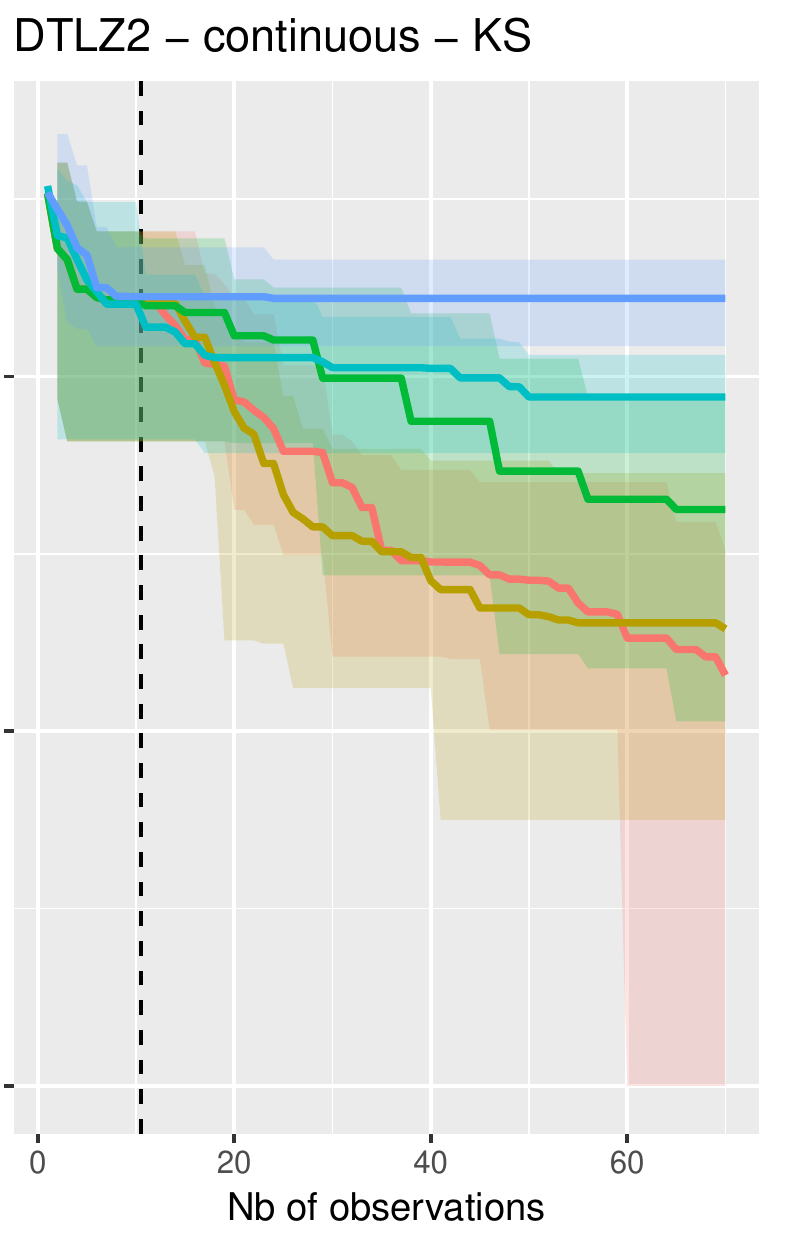}\\
\includegraphics[trim=0mm 0mm 32mm 0mm, height=70mm,  clip]{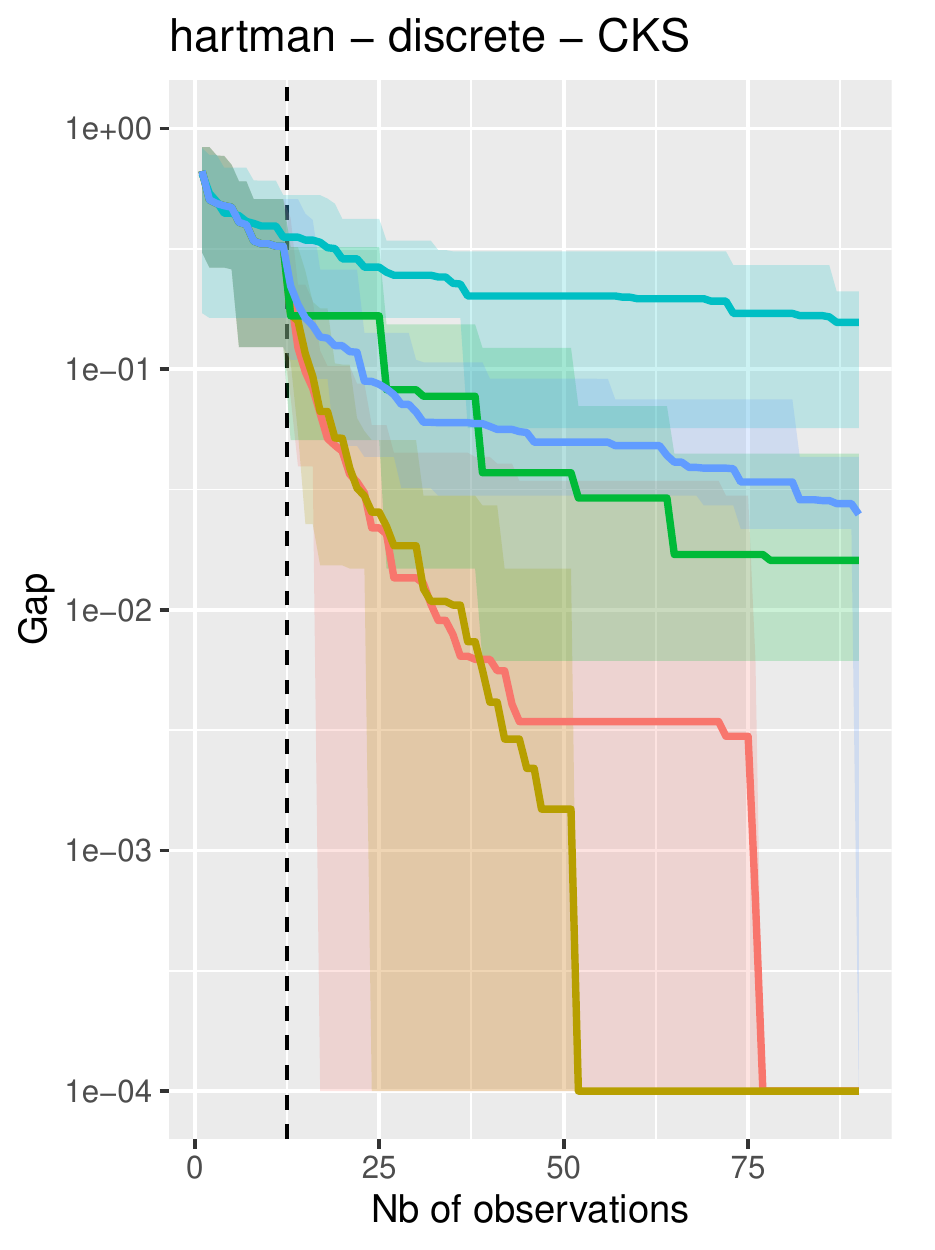}
\includegraphics[trim=0mm 0mm 32mm 0mm,height=70mm,  clip]{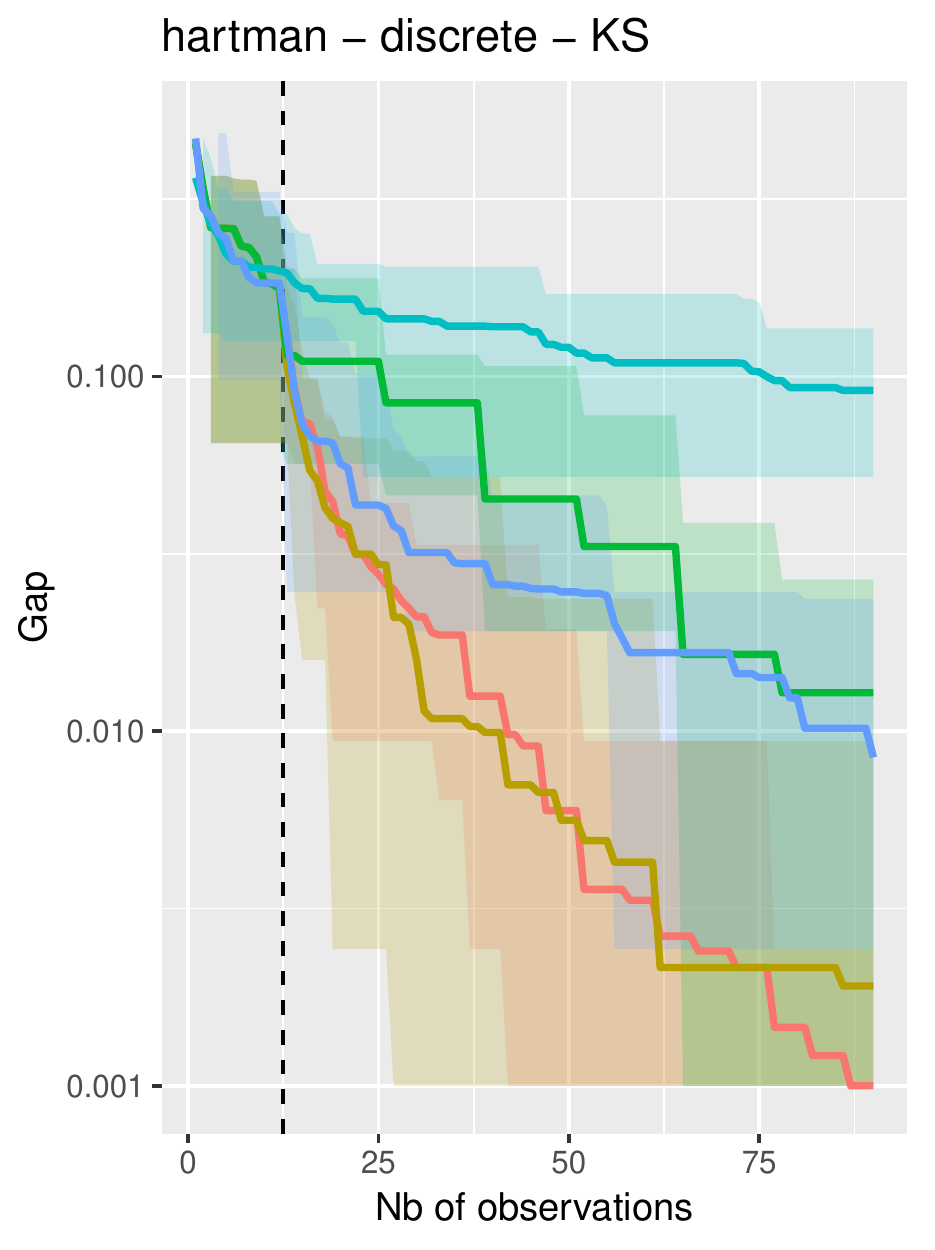}
\includegraphics[trim=15mm 0mm 32mm 0mm, height=70mm,  clip]{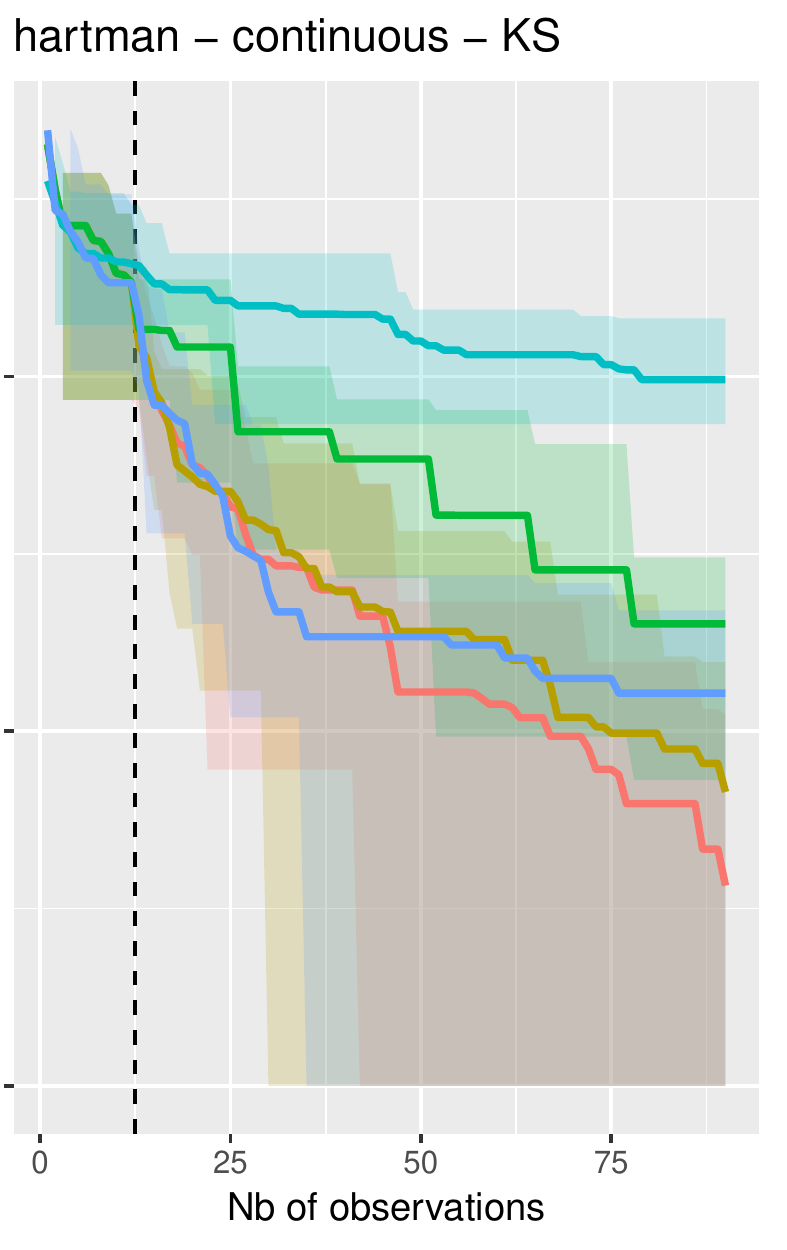}\\
\includegraphics[trim=0mm 0mm 0mm 0mm, height=5mm,  clip]{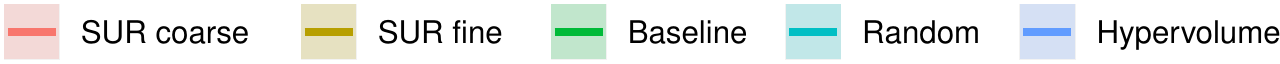}
 \caption{Optimality gap for KS (left) and CKS (right) on DTLZ2 (top) and
 hartman (bottom), for both discrete and continuous setups. Note that for
 readability, the Gap values are thresholded to $10^{-3}$ for KS and to
 $10^{-4}$ for CKS.}
\label{fig:results_toy}
\end{figure}

Convergence results are provided in Figure \ref{fig:results_toy} in terms of
optimality gap with respect to the minimum benefit ratio, that is:
\begin{equation*}
    G_{KS}(\s) = \min_{1 \leq i \leq \nobj} r^{(i)}(\s^*) - \min_{1 \leq i \leq \nobj} r^{(i)}(\s),
    \quad
    G_{CKS}(\s) = \min_{1 \leq i \leq \nobj} r_C^{(i)}(\s^\dagger) - \min_{1 \leq i \leq \nobj} r_C^{(i)}(\s),
\end{equation*}
with $\s$ the evaluated solution, $\s^*$ and $\s^{\dagger}$ the actual KS and
CKS solutions, respectively, $r$ the ratio based on the actual nadir and
utopia and $r_C$ (as in Eq.\ \ref{eq:defCKS}) computed using the actual
objective values at the auxiliary set.

On DTLZ2, we see clearly that the hypervolume alternative never seeks central
compromise solutions, and performs even worse than random search (this is
despite the fact that the algorithm performs very well in terms of the
hypervolume metric, see the empirical results in \citet{binois2019gpareto}).
The hypervolume approach is more competitive on hartman, which is due to the
fact that the Pareto set actually corresponds to a relatively small part of
the input space.

For both KS and CKS and both problems in the discrete case, the optimality gap
is approximately log-linear on average, with a couple of runs that find the
exact solution (gap $\leq 10^{-4}$) within a very restricted budget. On all
cases, SUR substantially outperform the baselines.

The difference between the coarse and fine SUR approaches is quite small on
most cases, the only significant difference being on hartman, CKS, where the
more expensive approach allows convergence to the exact solution much more
quickly than the coarse one. On average, we can conclude that the numerical
schemes to compute the SUR criterion provide a sufficient accuracy, and our
``coarse'' configuration results in reasonable running times (Figure
\ref{tab:configs}) compared to the cost of an expensive experiment such as
those described in the next subsections.

The effect of working on discretized versions of the problem is barely visible
on DTLZ2 results, but clearly affects the rate of convergence on the hartman
problem (although both SUR versions still converge to the actual solution).

Finally, an additional baseline is to reduce the number of objectives in order
to use efficiently the standard multi-objective algorithms. We report in
Appendix E additional experiments on the hartman case, that show that despite
some correlation between objectives, dropping all but two objectives has on
average a dramatic effect on the minimal ratios.

\subsection{Training of a convolutional neural network}

A growing need for BO methods emerges from machine learning applications, to
replace manual tuning of hyperparameters of increasingly complex methods. One
such example is with hyperparameters controlling the structure of a neural
network. Since such methods are integrated in products, accuracy is not the
only concern; and additional objectives have to be taken into account, such as
prediction times.

\subsubsection{Problem description}

\begin{figure}[htbp]
\centering
 \includegraphics[trim=0mm 0mm 0mm 0mm, width=\textwidth, clip]{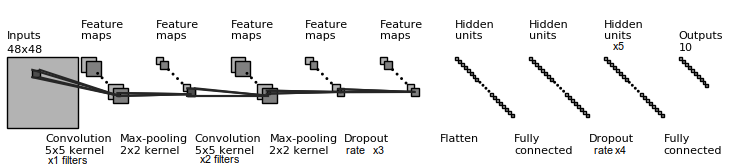}
 \caption{Architecture of the CNN on the MNIST data.}
\label{fig:CNN}
\end{figure}

We consider here the training of a convolutional neural network (CNN) on the
classical MNIST data \citep{lecun1998gradient}, with $60,000$ handwritten
digits for training and an extra $10,000$ for testing. We use the
\texttt{keras} package \citep{Allaire2018} to interface with the high-level
neural networks API Keras \citep{chollet2015keras} to create and train a CNN.
We follow a common structure for such a task, represented in Figure
\ref{fig:CNN}, with a first 2D convolutional layer, a first max pooling layer,
then a second 2D convolutional layer and a second max pooling layer. Max
pooling consists in keeping only the max over a window, introducing a small
amount of translational invariance. Dropout, that is, randomly cutting off
some neurons to increase robustness, is then applied before flattening to a
dense layer, followed by another dropout before the final dense layer. Because
of the dropout phases, the performance is random. This is handled by repeating
five times each experiment, which takes up to 30 minutes on a desktop with a
3.2 Ghz quad-core processor and 4 Go of RAM.

The six hyperparameters to tune, detailed in Table \ref{tab:mnist_in} in
Appendix C, along with their range of variation, are the number of filters and
dropout rates of each layer, plus the number of units of the last hidden layer
and number of epochs.

The training of the CNN is performed based on the categorical cross-entropy. A
validation data set is extracted from the training data to monitor
overfitting, representing 20\% of the initial size of the training data. The
progress can be monitored by taking into account the accuracy (i.e.,
proportion of properly classified data) and cross-entropy on the training
data, on the validation data or on the test data. Of these six corresponding
objectives, training and validation accuracy are extremely correlated with the
cross-entropy but they can be kept with our proposed methodology. Having both
validation and testing metrics is also useful, the validation ones are driving
the training, while the testing ones are less likely biased with data never
used in the training. As additional objectives, we considered the training
time of the CNN, as well as the prediction time on the testing data. This
latter is relevant, for instance, when using a pretrained CNN for a given
task. The eight objectives are summarized in Table \ref{tab:mnist_out} in
Appendix C.

We take a total budget of 100 evaluations, split in half between initial LHD
and sequential optimization using KS or CKS. GPs are trained by using Mat\'ern
5/2 kernels with an estimated linear trend. In this case, $N=500$ integration
points were selected out of $N_\text{large} = 10^6$ possible candidates,
renewed every iteration. The resulting time of each iteration is under a
minute for KS and less than 10 minutes for CKS.

\subsubsection{Results}

\begin{figure}[htbp]
\includegraphics[trim=22mm 22mm 14mm 20mm, width=0.85\textwidth, clip]{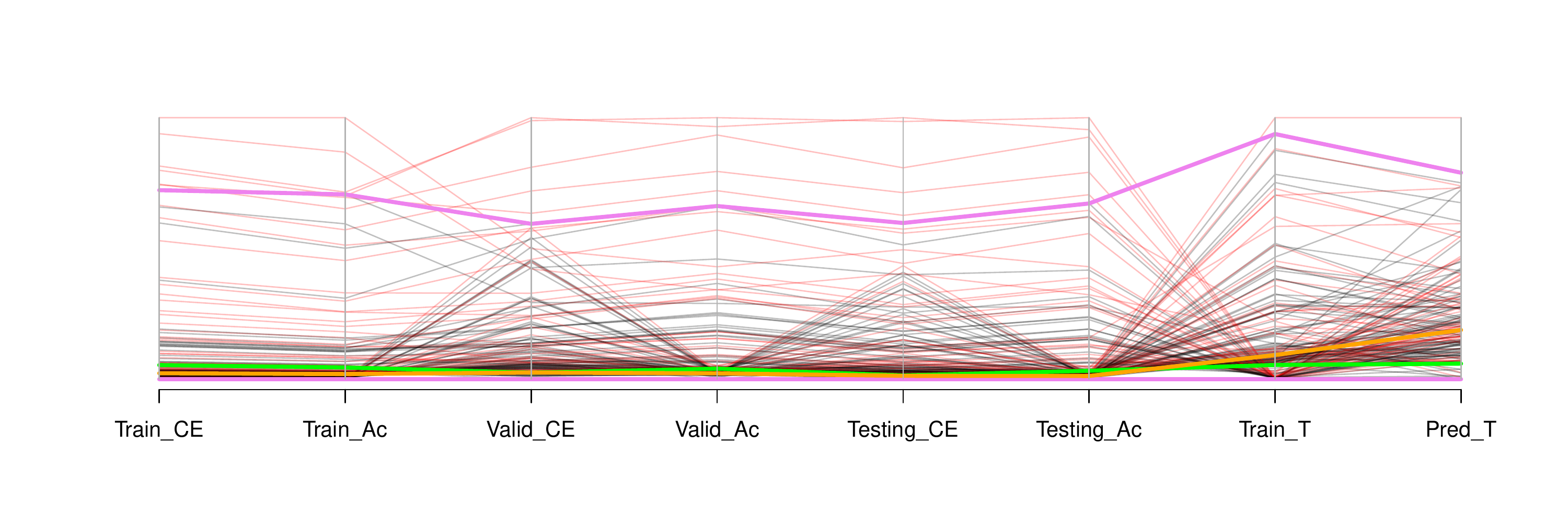}%
\includegraphics[trim=0mm 0mm 0mm 0mm, width=0.15\textwidth, clip]{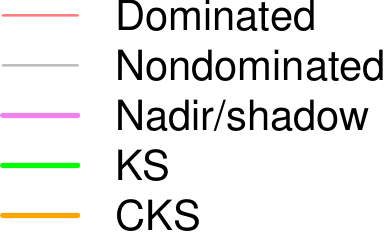}\\
\includegraphics[trim=22mm 10mm 14mm 20mm, width=0.85\textwidth, clip]{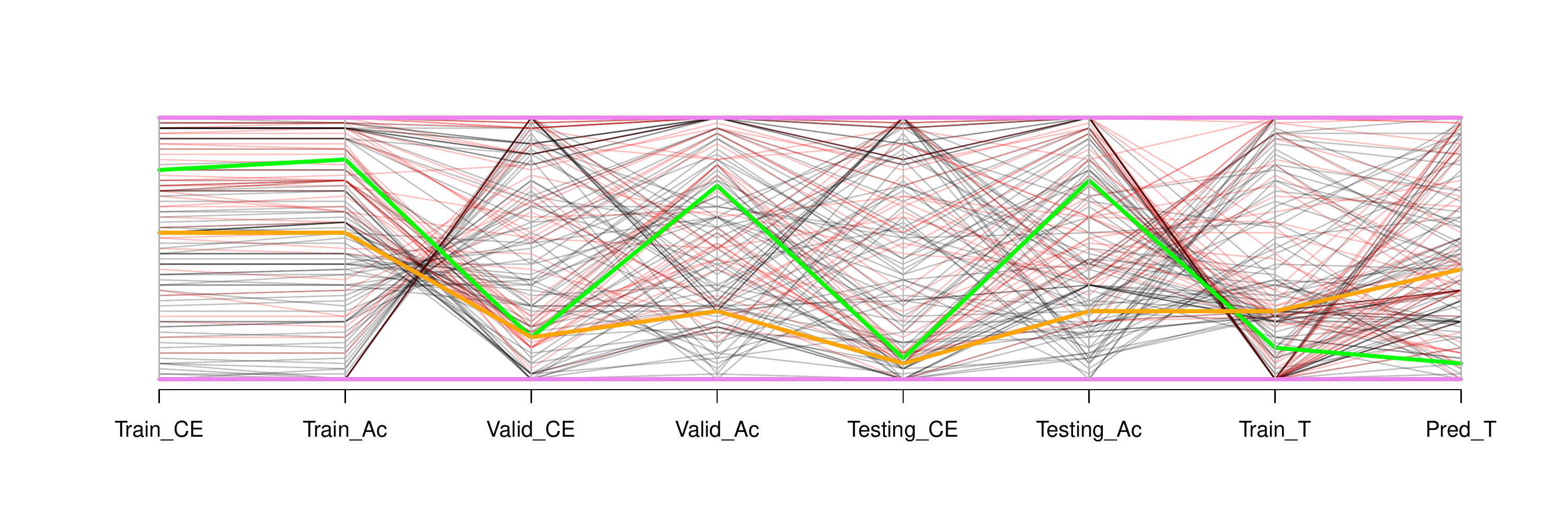}
 \caption{Parallel coordinates of the 150 observations, nadir and utopia
 points in the original (top) and copula (bottom) objective spaces (each line
 corresponds to a single observation, where objectives are rescaled on the
 graph). The colors highlight the dominated solutions and retained (KS and
 CKS) solutions.}
 \label{fig:parcoordscnn}
 \end{figure}

\begin{figure}[htbp]
\centering
\includegraphics[trim=0mm 4mm 5mm 13mm, width=\textwidth, clip]{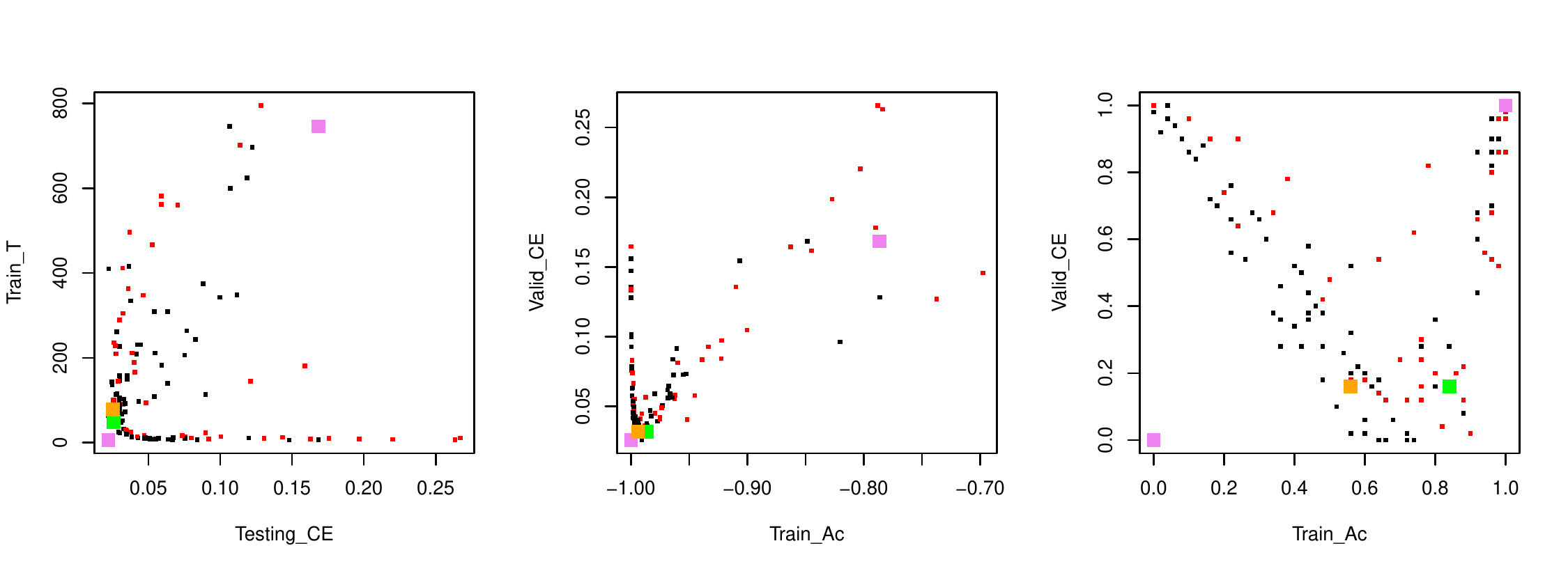}
 \caption{Projections on 2D spaces of the 150 observations, either in the
 original space (left, middle) or the copula space (right). The colors match
 the ones of Figure \ref{fig:parcoordscnn}.}
\label{fig:subpairs_cnn}
\end{figure}

Figure \ref{fig:parcoordscnn} represents the performance of networks obtained
during optimization, in parallel coordinates \citep[see ][]{li2017read} in
both original scale (top) and in the copula space (bottom). As minimization is
considered, minus the accuracies are reported. Selected 2D projections are
provided in Figure \ref{fig:subpairs_cnn}, while the whole set of projections
is in Appendix F. In the original objective coordinates, most observations
have low values close to objective minima, with performance objectives that
seem quite correlated. Nevertheless, 94 out of 150 observations are
nondominated, and when looking at ranks in the copula domain, the Pareto front
appears more complex, with conflicts between seemingly correlated objectives
as shown in Figure \ref{fig:subpairs_cnn}. These could correspond to
underfitting and overfitting. The nadir point is remarkably high in the
objective space, due to the leading trade-off between training time and
accuracy: bigger networks give better results but take more time. The KS and
CKS solutions are close to each other in the original objective space, but
much less so in the copula space, accounting for the concentration of low
values on the objectives. The KS solution favors quicker training and
prediction times than CKS, with seemingly marginally worse accuracy and
errors. But the properly scaled CKS solution recognizes that faster
architectures are actually quite extreme when looking at the ranks, thus
giving more weights to a small decrease in accuracy and errors, which is
significant in terms of ranks. As a result, the representation in parallel
coordinates is relatively flat for CKS (a flat line would indicate a perfectly
balanced solution) and quite irregular for KS, resulting in a much worse
minimal ratio.

\begin{figure}[htbp]
\centering
\includegraphics[trim=0mm 5mm 0mm 15mm, width=.49\textwidth, clip]{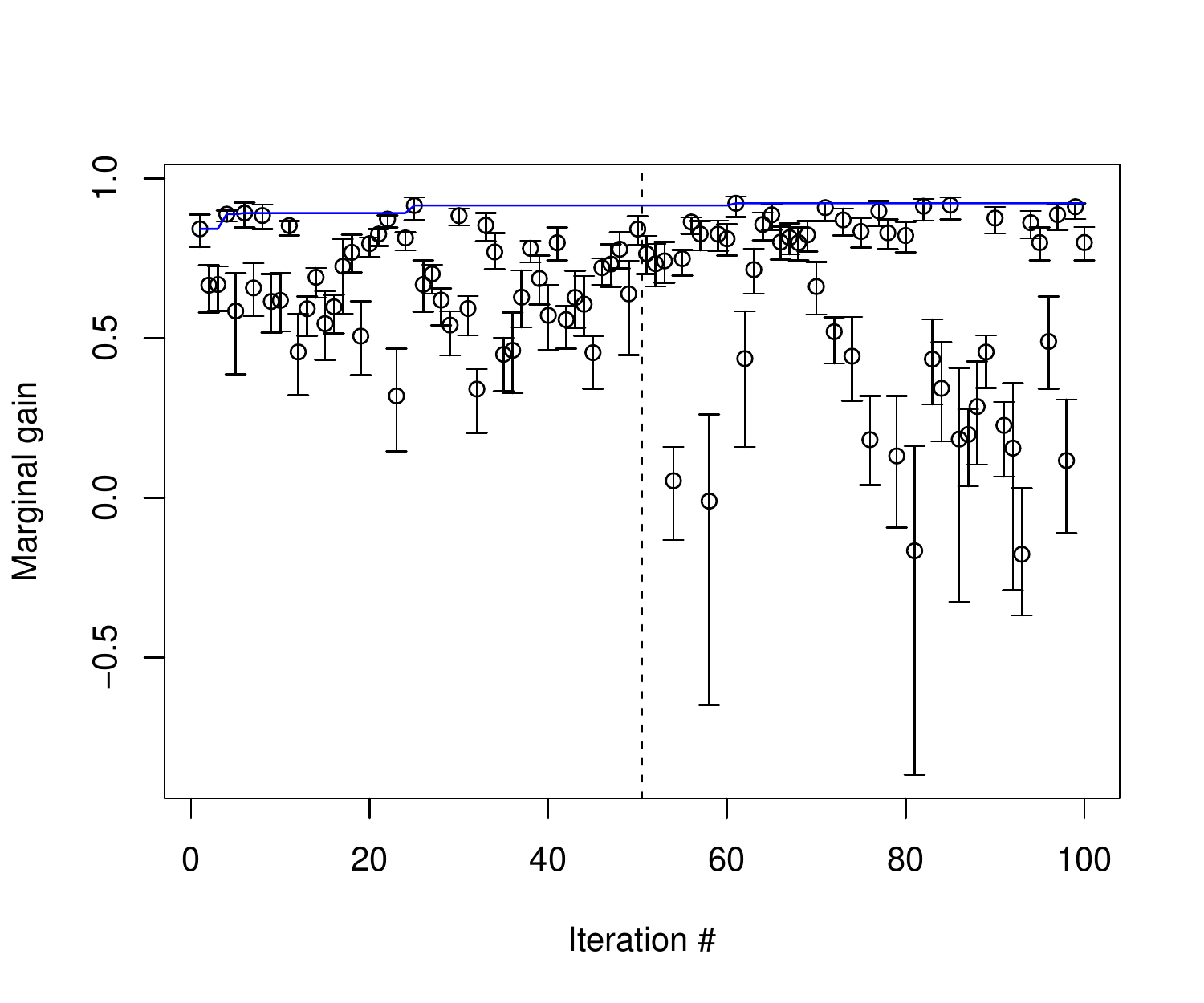}
\includegraphics[trim=0mm 5mm 0mm 15mm, width=.49\textwidth, clip]{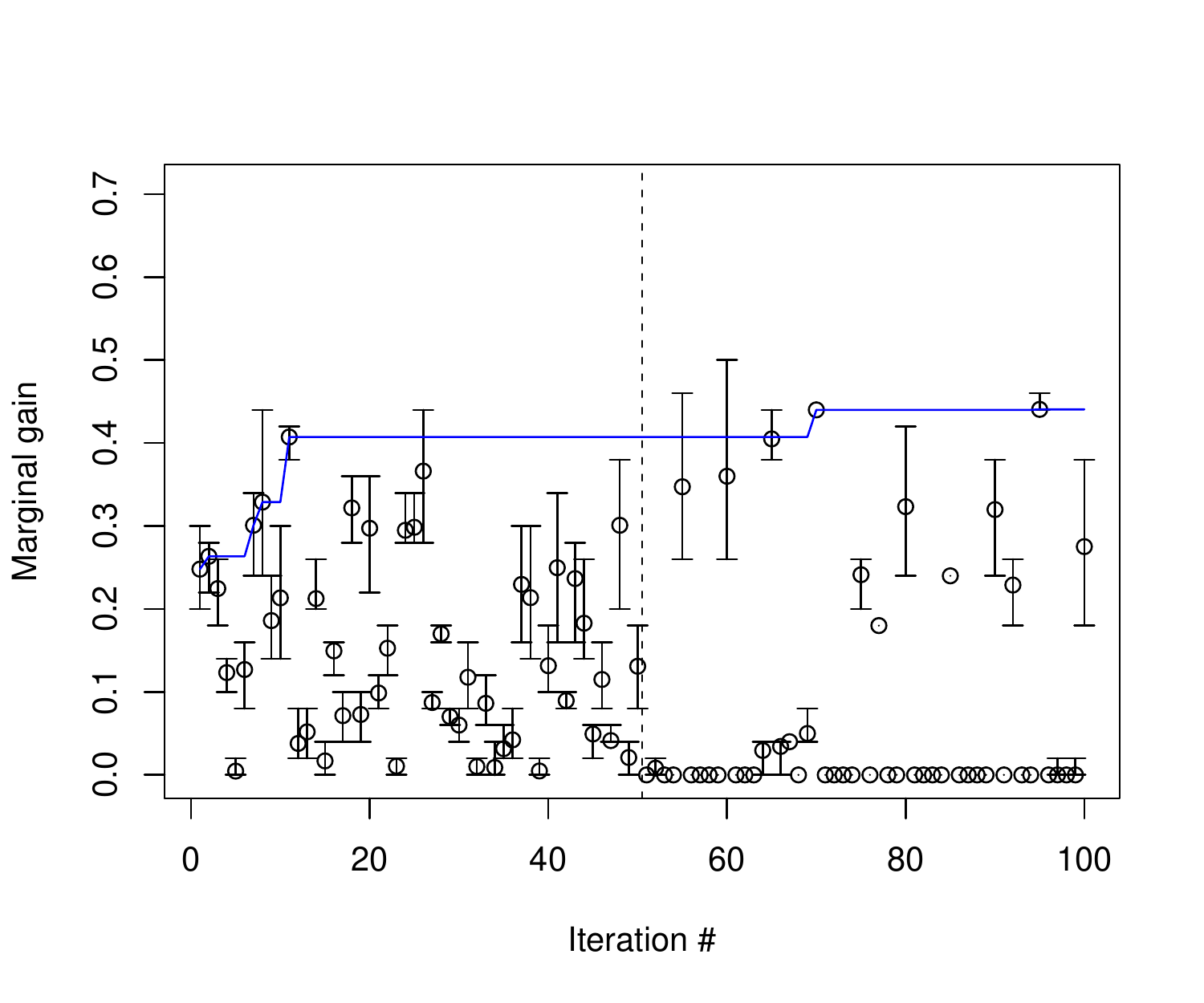}
 \caption{Evolution of the minimum ratio of the observations along the SUR
 algorithm for KS (left) and CKS (right). The vertical line shows where the
 SUR algorithm starts. Each graph corresponds to a specific ratio. The error
 bars account for uncertainty in the nadir and shadow, and the copula
 transformation, respectively.}
\label{fig:centrality_cnn}
\end{figure}

More insight is provided by Figure \ref{fig:centrality_cnn}, showing the
progress over the iterations of the minimum ratio (as in Equation
\ref{eq:defKS} or \ref{eq:defCKS}, respectively) of the observations. Due to
the computational cost, no ground truth is available, and in particular the
true nadir and shadow points are unknown, which renders uncertain the ratio
values. To estimate them as accurately and with as little bias as possible, we
collected the 150 observations (coming from the three optimization runs), fit
a GP model to all this set and use it to estimate nadir and shadow points, as
in Section \ref{sec:computingsur}. As a result, different ratios can be
computed using the different nadir-shadow pairs, which are represented using
error bars of Figure \ref{fig:centrality_gama}. A similar approach is followed
for the copula transformation to obtain uncertainty for CKS.

The seemingly small progress on KS (around 1\%) is related to the position of
the nadir, very high on each objective. In fact, the numerous evaluations with
almost zero (or negative) marginal gain correspond to the estimation of this
nadir, which is crucial for KS estimation. For CKS, the progress is more
visible, with a clear split between some exploitation steps (observations with
high minimal ratio) and exploration steps (with very small minimal ratio). It
is worth noticing that the degree of correlation between objectives did not
seem to be problematic. In fact, we illustrate more thoroughly the benefits of
keeping more objectives when they are all relevant in Appendix E.

Unless one is satisfied with the relative scaling of objectives of different
nature, the CKS solution has the advantage of being more robust and, in
addition, easier to estimate with a fixed nadir. The alternative to remain in
the original space is to define preferences via the disagreement point, as is
illustrated on the following test case.

For now, the six variables are taken as continuous. But these results could be
extended by increasing the number of variables, including categorical
variables affecting more profoundly the structure of the network, for instance
relying on the work of \cite{Roustant2018}.

\subsection{Calibration of an agent-based behavioral model}
Model calibration (sometimes referred to as inverse problem) consists of
adjusting input parameters so that the model outputs match real-life data. In
this experiment, we consider the calibration of the li-BIM model
\citep{libim}, implemented under the GAMA platform
\citep{taillandier2018building}, which exhibits several challenging features:
stochasticity, high numerical cost (approximately 30 minutes per run on a
desktop computer with a 3.60 GHz eight-core processor and 32 Go RAM), and a
large number of outputs.

\subsubsection{Problem description}
The Li-BIM model simulates the behavior of occupants in a building. It is
structured around the numerical modeling of the building and an evolved
occupational cognitive model developed with a belief-desire-intention (BDI)
architecture \citep{bourgais2017enhancing}. It simulates several quantities
that strongly depend on the occupants, such as thermal conditions, air
quality, lighting, etc. In the configuration considered here, three occupants
are simulated over a period of one year.

In order to reproduce realistic conditions, 13 parameters can be tuned,
related to either the occupant behavior or the building characteristics (see
Table \ref{tab:libim1} in Appendix C for details). Nine outputs ($G_1 \ldots
G_{9}$) should match some target values ($T_1 \ldots T_{9}$), chosen based on
records or surveys (see Table \ref{tab:libim2} in Appendix C). Since the model
is stochastic, we consider as objectives the squared expected relative
differences between the outputs and targets:
\begin{equation}
 y_i = \log \left( \left[ \esp \left( \frac{G_i - T_i}{T_i}\right) \right]^2 + \delta \right)
\end{equation}
The logarithm transformation is useful here to bring more contrast for values
close to zero, attenuated by a small $\delta$ (in our experiments,
$\delta=0.01$). In practice, we use estimates based on eight repeated runs.
Note that such an objective focuses on the average behavior without
considering the variability. As an alternative criterion, one may invert the
square and expectation.

To solve this 13-variable 9-objective problem, we proceed as follow. An
initial 100-point optimized LHD is generated, which is used to fit GPs
(constant trend, Mat\'ern 5/2 anisotropic covariance). From this initial
design, both KS and CKS SUR strategies are conveyed independently with 100
additional points for each.

In addition, a third solution is sought by using KS with a partly prespecified
disagreement point $\mathbf{d}$ to account for preferences (see Eq.\
\ref{eq:constraints}), so that the average error on five of the outputs does
not exceed a certain percentage (either 50\% or 30\%), the other outputs being
unconstrained. To do so, we use $$\tilde d_i = \min (N_i, c_i), \quad 1 \leq i
\leq 9,$$ $N_i = \max_{\x \in \mathbb{X}^*}\ y^{(i)}(\x)$ being the nadir i-th
coordinate, and $$c = \log ([0.5, 0.5, +\infty, +\infty, 0.3, 0.5, +\infty,
0.5, +\infty]^2).$$ We refer to this strategy as KSpref.

Importantly, the GPs are used to fit the expected values of the outputs of the
model ($\esp G_i$) instead of the objectives ($y_i$). This greatly improves
the prediction quality of the GPs (as $G_i$ is smoother than $(G_i-T_i)^2$),
while allowing us to convey our strategy almost without modification: on
Sections \ref{sec:sur} and \ref{sec:comput}, the drawings $\Yr$ and $\Fr$ are
obtained by first generating drawings $\Gr$ of $G$, then transforming them
($\log \Gr^2$).

We used $N^* = 1,000$ points for $\Xset^*$, taken from $\Xset_\text{large}$, a
$2 \times 10^5$ space-filling design, both renewed at each iteration.

\subsubsection{Results}
\begin{figure}[htbp]
    \includegraphics[trim=22mm 22mm 14mm 20mm, width=0.85\textwidth, clip]{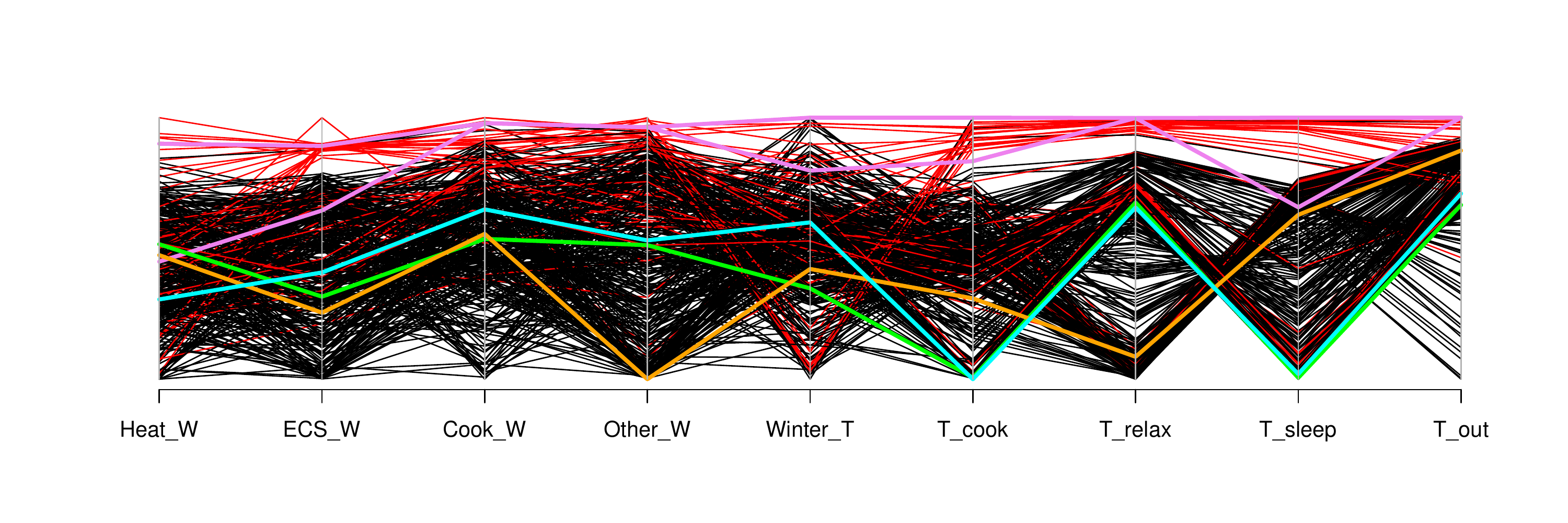}%
    \includegraphics[trim=0mm 0mm 0mm 0mm, width=0.15\textwidth, clip]{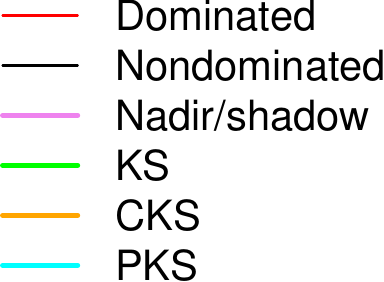}\\
    \includegraphics[trim=22mm 10mm 14mm 20mm, width=0.85\textwidth, clip]{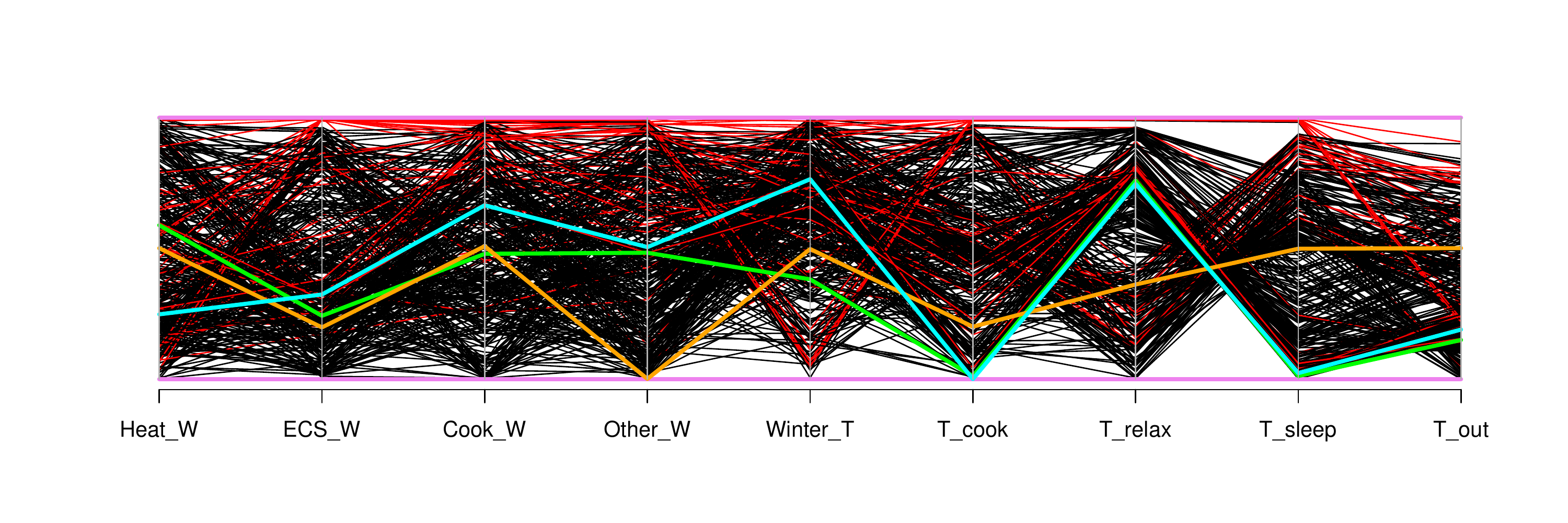}
     \caption{Parallel coordinates of the 400 observations, nadir and utopia
     points in the original (top) and copula (bottom) objective spaces (each
     line corresponds to a single observation, where objectives are rescaled
     on the graph). The colors highlight the dominated solutions and retained
     (KS, CKS and PKS) solutions.}
    \label{fig:parcoordsgama}
    \end{figure}

\begin{figure}[htbp]
\centering
\includegraphics[trim=0mm 3mm 5mm 13mm, width=\textwidth, clip]{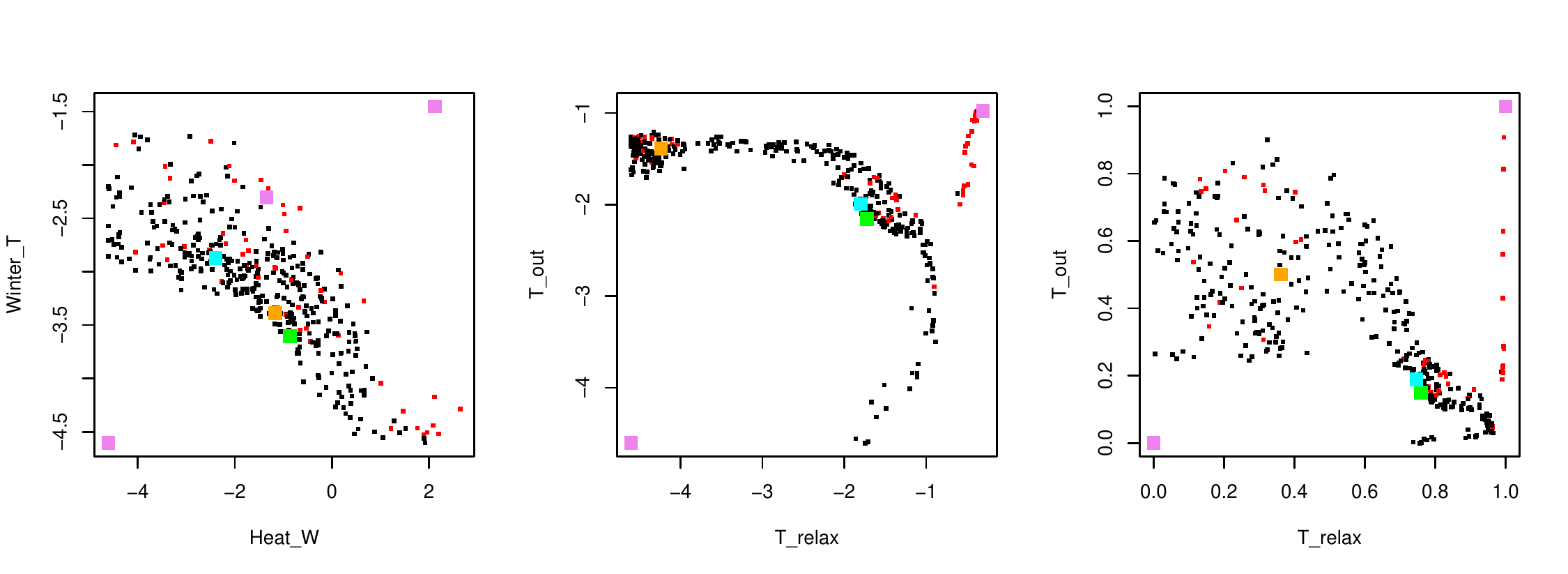}
\caption{Projections on 2D spaces of the 400 observations, either in the
original space (left, middle) or the copula space (right). The colors match
the ones of Figure \ref{fig:parcoordsgama}.}
\label{fig:subpairs_gama}
\end{figure}

The resulting designs of experiments and solutions are reported in graphical
form in Figure \ref{fig:parcoordsgama}, that shows parallel coordinates of the
observations  and Figure \ref{fig:subpairs_gama}, that show a couple of 2D
projections of the observations (the full set of 2D projections is available
in Appendix G).

As a preliminary observation, of the 400 points computed during this
experiment, only 55 were dominated. This result illustrates the exponential
growth of Pareto sets with the number of objectives. As a consequence, in
Figure \ref{fig:parcoordsgama} the non-dominated solutions (in black) cover
most of the objective ranges, without apparent structure, and the nadir (in
purple) is almost equal to the individual maxima.

Figure \ref{fig:subpairs_gama} shows that two pairs of objectives are strongly
conflicting. Indeed, the error on Heat\_W (the energy spent on heating) is in
competition with the one on Winter\_T (the room temperature in winter), which
is intuitive, and T\_relax and T\_out, that are both time spent on different
activities. The fact that no solution reconciles both pairs of objectives
shows an intrinsic discrepancy between the data and the model, and the
calibration task amounts to finding the best balance between the errors. Other
pairwise antagonisms are visible in Appendix, in particular between the last 4
objectives (all corresponding to times spent on activities), or between
T\_cook and W\_cook (time and energy corresponding to cooking).

On Figure \ref{fig:parcoordsgama}, we see first that the three solutions are
similar for some objectives (Heat\_W, ECS\_W, Cook\_W, Winter\_T) and differ
substantially for others (in particular Other\_W, Time\_relax, Time\_sleep,
Time\_out). The largest difference between the KS and KS with predefined
disagreement point (PKS) appears on objectives Heat\_W and Winter\_W. This is
explained by Figure \ref{fig:subpairs_gama} (left), where we can see that
imposing a constraint on Heat\_W shifts the trade-off between those
objectives. The large difference between CKS and the two other solutions is
explained by Figure \ref{fig:subpairs_gama} (middle and right). While KS is
very central in the original space, the large concentration of solutions for
large Time\_out values induces a strong difference in the copula space, where
CKS is central and the other solutions are not.

\begin{figure}[htbp]
\centering
\includegraphics[trim=0mm 5mm 0mm 15mm, width=.32\textwidth, clip]{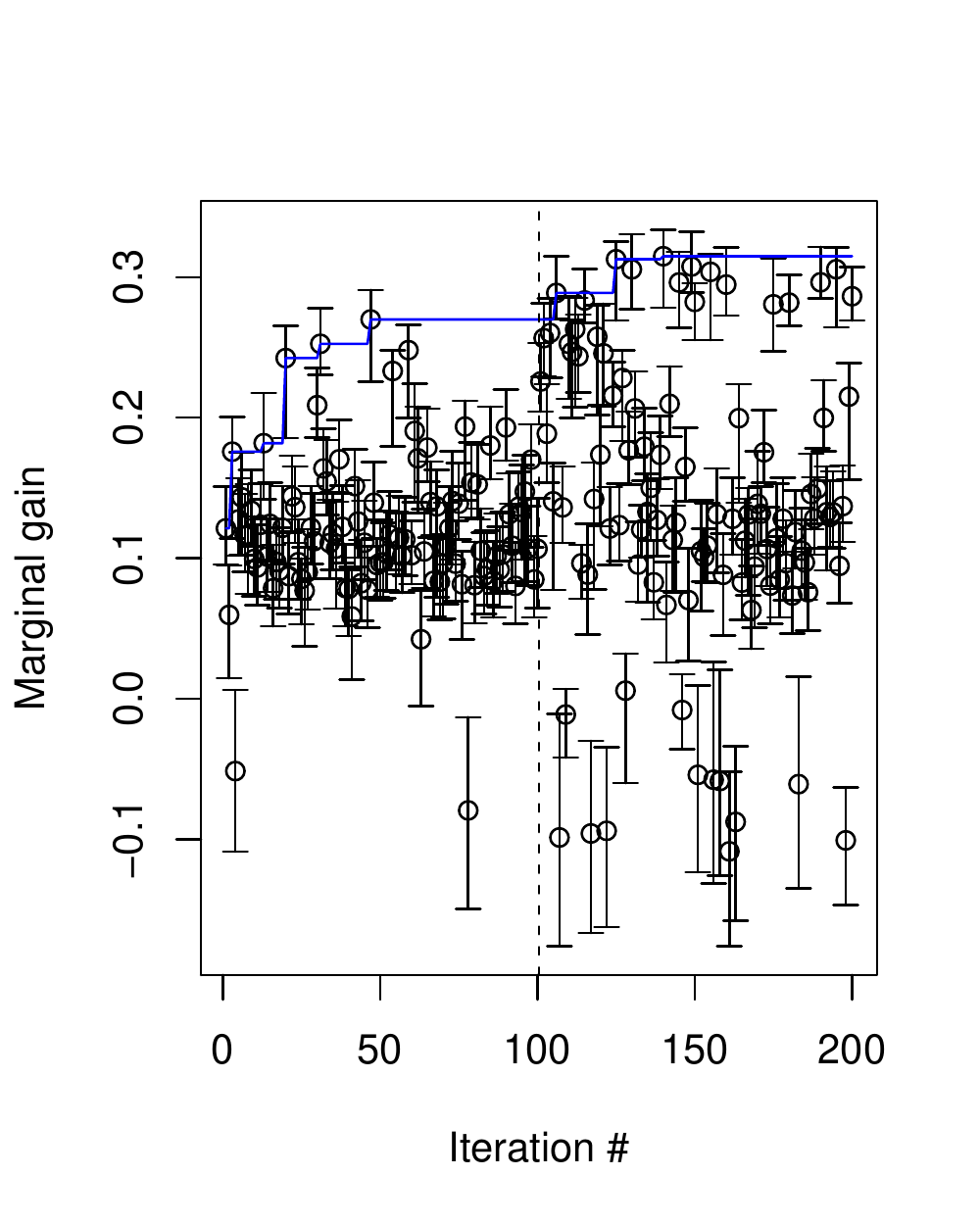}
\includegraphics[trim=0mm 5mm 0mm 15mm, width=.32\textwidth, clip]{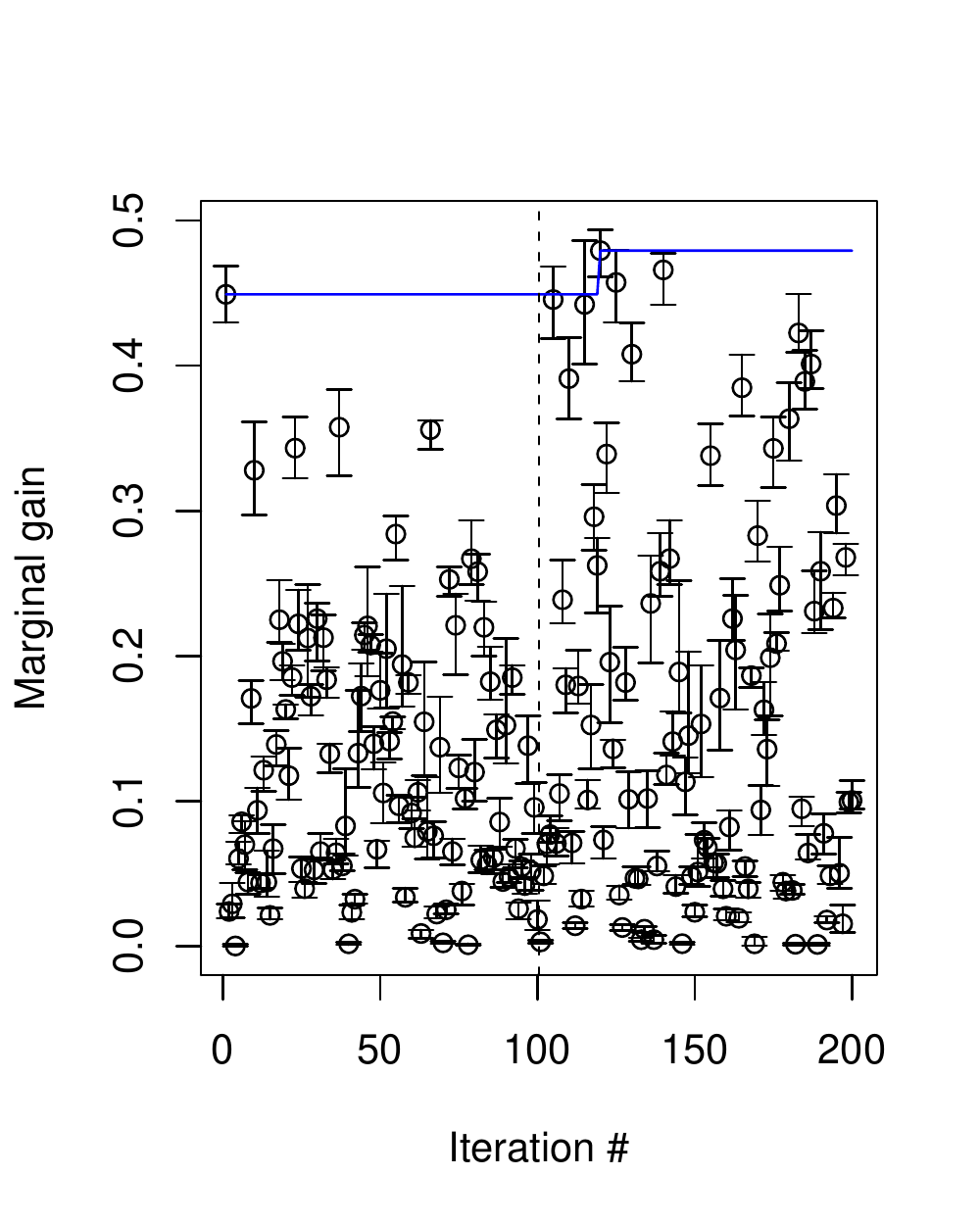}
\includegraphics[trim=0mm 5mm 0mm 15mm, width=.32\textwidth, clip]{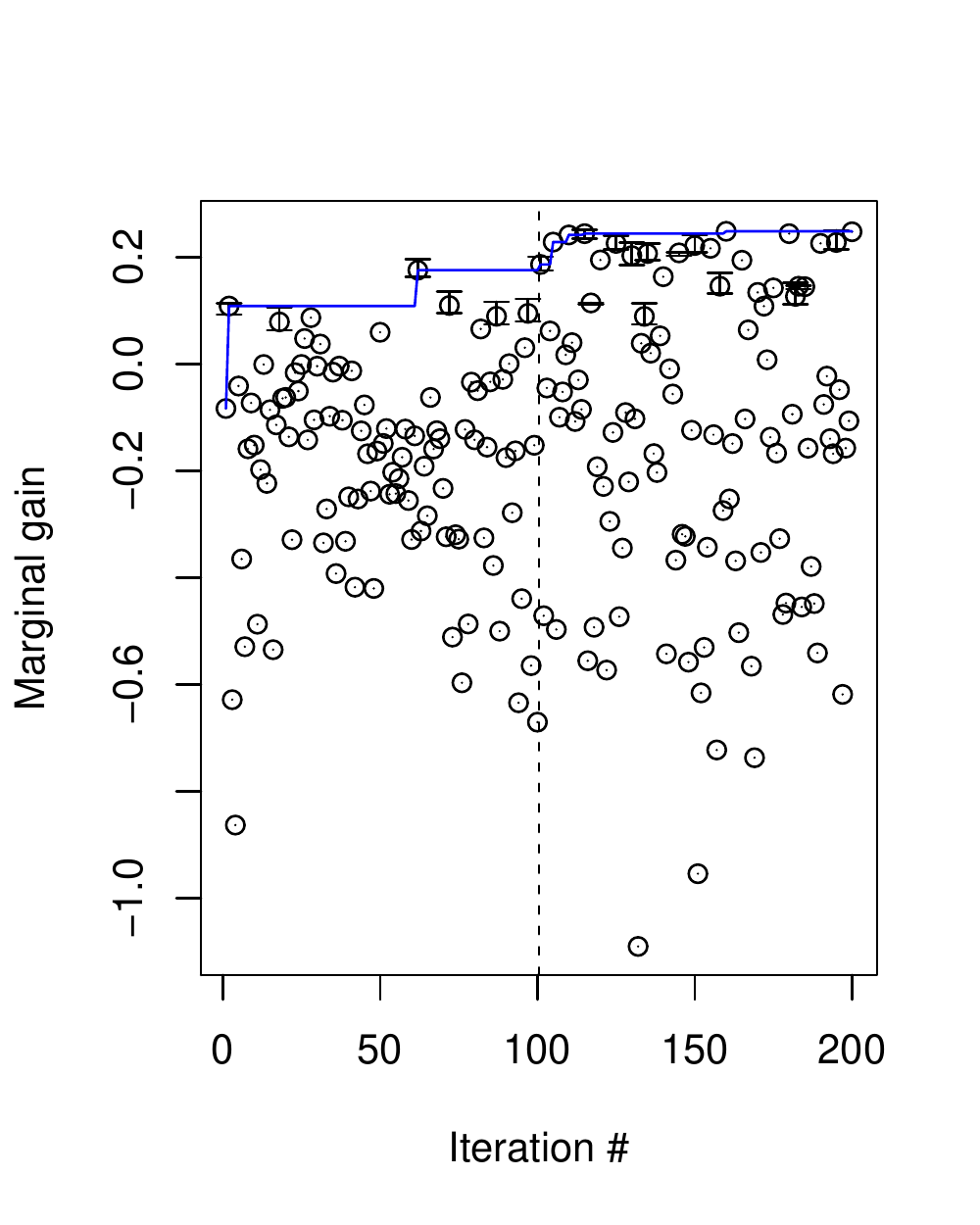}
 \caption{Evolution of the minimum ratio of the observations along the SUR
 algorithm for KS (left), CKS (middle) and KS with preferences (right). The
 vertical line shows where the SUR algorithm starts. Each graph corresponds to
 a specific ratio. The error bars account for uncertainty in the nadir and
 shadow and the copula transformation, respectively.}
\label{fig:centrality_gama}
\end{figure}

Finally, Figure \ref{fig:centrality_gama} shows the performance of our
algorithm, by representing the minimum ratio (as in Equation \ref{eq:defKS} or
\ref{eq:defCKS}, respectively) of the observations along the SUR iterations.
As in the CNN case, no ground truth is available, and we use a similar
approach to account for the uncertainty in the nadir-shadow line and the
copula estimates.

We can see that on all cases, the SUR algorithm improves substantially over
the space-filling observations. Contrary to the toy problems, here SUR worked
better on KS (with and without preferences) than on CKS (which might be simply
explained by the presence of a very good initial solution for CKS). One may
finally notice that the uncertainty is highest for KS and smallest for KS with
preferences. This indicates that finding the nadir may be the most difficult
of the learning tasks for this problem.

\section{Conclusion and future work}\label{sec:conclusion}

In this article, we tackled many-objective problems by looking for a single
well-balanced solution, originating from game theory. Two alternatives to this
solution have been proposed, either by imposing preferences on some objective
values (by specifying a disagreement point)  or by working in the space of
copulas, the latter solution being insensitive to monotonic transformations of
the objectives. Looking for these solutions is in general a complex learning
task. We proposed a tailored algorithm based on the stepwise uncertainty
reduction paradigm, that automatically performs a trade-off between the
different learning tasks (estimating the ideal point, the nadir point, and the
marginals or exploring locally the space next to the estimated solution). We
tested our algorithm on three different problems with growing complexity and
found that well-balanced solutions could be obtained despite severely
restricted budgets.

Choosing between the alternatives seems highly problem-dependent. On DTLZ2, KS
appears as a more ``central'' (hence desirable) solution, while on the CNN
tuning problem CKS is clearly a better choice, and on the calibration problem
both solutions seem equally good. This difference may be imputed to scaling of
the objectives prior to optimization. On DTLZ2, all objectives behave
similarly, while on CNN some are strongly heavy-tailed. On the other hand, CKS
is dependent on the problem formulation (for instance by considering as an
input $x \in [0, 1]$ or $\log(x) \in [-\infty, 0]$ could lead to very
different results). Hence, a reasonable choice would be KS when the objectives
formulation makes them directly comparable, and CKS when the input formulation
is pertinent (in the sense that all parts of the input space have a roughly
comparable effect on the outputs). As CKS is in general less likely to favor a
particular objective over others, it may be used for more exploratory studies,
and KS for a finer design on a pre-explored problem. Incorporating user
preferences, as in the calibration problem, was proved easy and efficient and
may direct the choice toward KS in such a case.

Many potential lines of future works remain. First, we may consider
batch-sequential strategies instead of one observation at a time. This
approach was not necessary in our test problems, since parallel computation
was used to repeat simulations and average out noise. The SUR strategy
naturally adapts to this case \citep{chevalier2014fast}. In practice, one may
also have to deal with asynchronous returns of batches, and possibly the use
of several black-boxes with different costs, as in
\citet{hernandez2016general}. An unused degree of freedom here is the number
of replications to handle the noise, which might improve substantially the
practical efficiency \citep{jalali2017comparison, binois2019replication}. To
do so, one may combine the presented approach with an efficient scheme for
designing replications and estimating noise, as done by
\citet{binois2016practical}.

Another direction is to consider alternative equilibria. A promising idea is
to use a set of disagreement points instead of a single one: this could result
in more robust solutions and potentially a small set of Pareto-optimal
solutions, which might be preferred by decision-makers. Combining the KS
solution with Nash games has been suggested in the game theory literature.
\citet{conley1991bargaining} proposed using Nash equilibria as disagreement
points,  provided there exists a natural, or some relevant, splitting of the
decision variable among the players.  In the multicriteria Nash game framework
\citep{ghose1989solution},  noncooperative players have to handle individual
vector payoffs, so that to each request by other players, they have to respond
with some payoff, rationally selected from their own Pareto front; a good
candidate would be the KS solution. Both alternatives were poorly investigated
in practice, mainly because their potential of application  seems to be
hindered by the lack of efficient tools that allow for an acceptable
implementation.  We think that our algorithmic framework could apply to both
cases and provide a first solution to such problems.

Copula spaces have been used here mostly for rescaling. Taking advantage of
the predictive capacity of copulas  to accelerate the estimation of the Pareto
front might accelerate substantially the search for the CKS solution.
Furthermore, combining efficiently the two types of metamodels (GPs and
copulas), as done by \cite{wilson2010copula}, may lead to new theoretical
advances and algorithms. A complementary line of work would be to study the
influence of the input distribution on CKS in order to define (or infer)
distributions that lead to solutions that are independent of (or at least
robust with respect to) the input problem formulation.

We leave to future work theoretical considerations on the convergence of the
approach, following for instance the recent works on information-directed
sampling by \citet{russo2014learning} or on SUR by
\citet{bect2019supermartingale}.

\acks{
The work of MB is supported by the U.S.\ Department of Energy, Office of
Science, Office of  Advanced Scientific Computing Research under Contract No.\
DE-AC02-06CH11357. We thank Gail Pieper for her useful language editing. }

\appendix
\section*{Appendix A: Quantities related to Gaussian processes}
\paragraph{GP moments}
We provide below the equations of the moments of a GP conditioned on $n$
(noisy) observations $\f = (f_1, \ldots, f_n)$. Assuming a kernel function
$\sigma$ and a mean function $m(\x)$, we have
\begin{eqnarray*}
\mu_n(\x) & =& m(\x) + \lambda(\x) \left(\f - m(\x) \right),  \\
\sigma_n^2(\x, \x') & =& \sigma(\x, \x') - \lambda(\x) \sigma(\x', \X_n),
\end{eqnarray*}
where 
\begin{itemize}
 \item $\lambda(\x) := \sigma(\x, \X_n)^\top \sigma(\X_n, \X_n)^{-1}$,
 \item  $\sigma(\x, \X_n) := (\sigma(\x, \x_1), \dots, \sigma(\x, \x_n))^\top$ and 
 \item $\sigma(\X_n, \X_n) := (\sigma(\x_i, \x_j) + \tau_i^2 \delta_{i=j})_{1 \leq i,j \leq n}$,
\end{itemize}
$\delta$ standing for the Kronecker function.

Commonly, $\sigma$ belongs to a parametric family of covariance functions such
as the Gaussian and Mat\'ern kernels, based on hypotheses about the smoothness
of $y$. Corresponding hyperparameters are often obtained as maximum likelihood
estimates; see e.g., \citet{Rasmussen2006} for the corresponding details.

\paragraph{Expected improvement}
Denote $f_{\min} = \min_{1 \leq i \leq n}(f_i)$ the minimum of the observed
values. The expected improvement is the expected positive difference between
$f_{\min}$ and the new potential observation $Y_n(\x)$:
\begin{eqnarray*}
 EI(\x) &=& \esp \left( \max\left((0, f_{\min} - Y_n(\x) \right) \right) \\
        &=& (f_{\min} - \mu_n(\x))\Phi\left( \frac{f_{\min} - \mu_n(\x)}{\sigma_n(\x,\x) } \right) + \sigma_n^2(\x) \phi\left( \frac{f_{\min} - \mu_n(\x)}{\sigma_n(\x,\x) } \right),
\end{eqnarray*}
where $\phi$ and $\Phi$ are respectively the PDF and CDF of the standard Gaussian variable.

For maximization, we use, with $f_{\max} = \max_{1 \leq i \leq n}(f_i)$:
\begin{eqnarray*}
 \EIinv(\x) &=& \esp \left( \max\left((0, Y_n(\x) - f_{\max} \right) \right).
\end{eqnarray*}

\paragraph{Probability to belong to a box ($p_{\text{box}}$)}

Let $LB \in \Rset^\nobj$ and $UB \in \Rset^\nobj$ such that $\forall 1 \leq i
\leq \nobj, LB^{i} < UB^{i}$ define a box in the objective space. Defining
$\boldsymbol{\Psi} = \left[\Psi(\Yr_1), \ldots, \Psi(\Yr_M) \right]$ the
$\nobj \times M$ matrix of simulated KS solutions, we use
\begin{equation*}
 \forall 1 \leq i \leq \nobj \qquad LB^{i} = \min \boldsymbol{\Psi}_{i, 1 \ldots M} \quad \text{ and } \quad UB^{i} = \max \boldsymbol{\Psi}_{i, 1 \ldots M}.
\end{equation*}

Then, the probability to belong to the box is
\begin{equation*}
 p_{\text{box}}(\x) = \prod_{i=1}^\nobj \left[ \Phi\left(\frac{UB^{i} - \mu_n^i(\x)}{\sigma_n^i(\x,\x)} \right) - \Phi\left(\frac{\mu_n^i(\x) - LB^{i}}{\sigma_n^i(\x, \x)} \right) \right].
\end{equation*}

\paragraph{Probability of nondomination}
Let $\Xset_n^*$ be the subset of nondominated observations. The probability of
non-domination is
\begin{equation*}
 \PND(\x) = \prob\left( \forall \x^* \in \Xset_n^*, \exists k \in \{1, \ldots, \nobj\} \text{ such that } \nonumber Y_n^{(k)}(\x) \leq Y_n^{(k)}(\x^*) \right).
\end{equation*}
Using the GP equations for $Y_n$, one can compute $p_{ND}(\x)$ in closed form,
for $p \leq 3$. We refer to the work of \citet{couckuyt2014fast} for the
formulas expressed in an efficient form. For a larger number of objectives,
this probability must be computed by Monte-Carlo.

\section*{Appendix B: CNN lists of inputs and outputs}

\begin{table}[hp]
\centering
\caption{List of inputs for the CNN training problem.}
\begin{tabular}{clcc}
      & Description                                     & Min & Max \\ \hline
$x_1$ & number of filters of first convolutional layer  & 2             & 100           \\
$x_2$ & number of filters of second convolutional layer & 2             & 100           \\
$x_3$ & first dropout rate                              & 0             & 0.5           \\
$x_4$ & second dropout rate                             & 0             & 0.5           \\
$x_5$ & number of units of dense hidden layer           & 10            & 1000          \\
$x_6$ & number of epochs                                & 2             & 50
\end{tabular}

\label{tab:mnist_in}
\end{table}

\begin{table}[hp]
\centering
\caption{List of outputs for the CNN training problem}
\begin{tabular}{lc}
Name                     & Unit \\ \hline
Training cross-entropy   &  -   \\
Training accuracy        &  -   \\
Validation cross-entropy &  -   \\
Validation accuracy      &  -   \\
Testing cross-entropy    &  -   \\
Testing accuracy         &  -   \\
Training time            & s    \\
Prediction time          & s   
\end{tabular}

\label{tab:mnist_out}
\end{table}

\section*{Appendix C: li-BIM list of inputs and outputs}
\begin{table}[ht]
\centering
\caption{List of inputs of the li-BIM model.}
 \begin{tabular}{cclcc}
  Name & Unit & Meaning & Min & Max\\
    \hline
  Cth\_h  & J/K & Thermal capacity of the dwelling & $10^6$ & $10^7$\\
  RD & W & Average power of the relaxing devices & 150 & 1000 \\
  HW & W &   Power of the boiler to produce hot water & 500 & 3000\\
  CD & W &   Average power of the cooking devices & 100 & 300\\
  PmaxHeat & kW &   Maximum power of the boiler to heat dwelling & 500 & 3000\\
  SensitiveCold & [0-1] &   Sensitivity of the occupant to cold temperature & 0.2 & 1.0\\
  SensitiveWarm & [0-1] &   Sensitivity of the occupant to warm temperature & 0.5 & 1\\
  NbHfreshair & hours & Average number of hours between two outings & 16 & 36 \\
  NbHtire & hours &   Average number of hours between two sleeps & 8 & 60\\
  NbHhungry & hours &   Average number of hours between two meals & 6 & 16 \\
  NbHdirty & hours &   Average number of hours between two showers/baths & 20 & 48 \\
  Deltamodif & hours & Time before new action if previous action insufficient & 5 & 30 \\
  Thermal\_effort & Celsius & Max difference to ideal temperature before acting & 2 & 15 \\
 \end{tabular}

\label{tab:libim1}
\end{table}

\begin{table}[ht]
\centering
\caption{List of outputs of the li-BIM model.}
 \begin{tabular}{cclc}
  Name & Unit & Meaning & Target \\
  \hline
  Heat\_W & kWh & Total energetic consumption of heating devices & 1384\\
  ECS\_W & kWh & Total energetic consumption of hot water devices & 1198\\
  Cook\_W & kWh & Total energetic consumption of cooking devices & 306\\
  Other\_W & kWh & Total energetic consumption of other devices & 1751\\
  Winter\_T & \degree C & Average temperature during winter & 21.8\\
  Time\_cook & hours & Average time spent cooking & 0.9\\
  Time\_relax & hours & Average time spent relaxing & 3.7\\
  Time\_sleep & hours & Average time spent sleeping & 8.35\\
  Time\_out & hours & Average time spent outside the building & 0.58
  \end{tabular}
 
\label{tab:libim2}
\end{table} 
      
\FloatBarrier
\section*{Appendix D: hartman toy problem}
The original hartman function, defined on $[0,1]^6$, is:
\begin{equation*}
    h(\mathbf{x}) = \frac{-1}{1.94} \left[2.58 + \sum_{i=1}^4{ C_i \exp \left( -\sum_{j=1}^6{  a_{ji} \left( x_j - p_{ji} \right)^2  } \right) } \right]
\end{equation*}
with:
$\mathbf{C} = [1.0, 1.2, 3.0, 3.2]$,

$\mathbf{a} = \left[ \begin{array}{cccc}
	         10.00 &  0.05&  3.00& 17.00 \\
           3.00& 10.00&  3.50&  8.00\\
           17.00& 17.00&  1.70&  0.05\\
           3.50&  0.10& 10.00& 10.00\\
           1.70&  8.00& 17.00&  0.10\\
           8.00& 14.00&  8.00& 14.00
\end{array} \right]$, 
$\mathbf{p} = \left[ \begin{array}{cccc}
           0.1312& 0.2329& 0.2348& 0.4047\\
           0.1696& 0.4135& 0.1451& 0.8828\\
           0.5569& 0.8307& 0.3522& 0.8732\\
           0.0124& 0.3736& 0.2883& 0.5743\\
           0.8283& 0.1004& 0.3047& 0.1091\\
           0.5886& 0.9991& 0.6650& 0.0381
\end{array} \right].$

To obtain six objectives, we rescale and translate the inputs, so that:
\begin{eqnarray*}
    f_1(\x) &=& -\log(- h([\tilde x_1, \tilde x_2, \tilde x_3, \tilde x_4, \tilde x_5, \tilde x_6])\\
    f_2(\x) &=& -\log(- h(\nicefrac{1}{2} + [\tilde x_6, \tilde x_5, \tilde x_4, \tilde x_3, \tilde x_2, \tilde x_1])\\
    f_3(\x) &=& -\log(- h([\tilde x_2, \nicefrac{1}{2} + \tilde x_4, \tilde x_6, \tilde x_1, \nicefrac{1}{2} + \tilde x_3, \nicefrac{1}{2} + \tilde x_5])\\
    f_4(\x) &=& -\log(- h([\nicefrac{1}{2} +\tilde x_5, \tilde x_3, \nicefrac{1}{2} +\tilde x_1, \nicefrac{1}{2} +\tilde x_2, \tilde x_4, \tilde x_6])\\
    f_5(\x) &=& -\log(- h([\tilde x_3, \tilde x_6, \tilde x_1, \nicefrac{1}{2} +\tilde x_4, \nicefrac{1}{2} +\tilde x_2, \nicefrac{1}{2} +\tilde x_5])\\
    f_6(\x) &=& -\log(- h([\nicefrac{1}{2} +\tilde x_4, \nicefrac{1}{2} +\tilde x_2, \nicefrac{1}{2} +\tilde x_6, \tilde x_5, \tilde x_3, \tilde x_1]),
\end{eqnarray*}
with $\tilde \x = \nicefrac{\x}{2}$.

\section*{Appendix E: Effect of taking more objectives}

To avoid many-objective challenges, a natural way is to select fewer
objectives that are taken as more important, or looking at correlations. In
this section we illustrate that even if discarding objectives gives some
improvements in the solutions over the selected ones, overall more may be lost
on the discarded objectives. To this end, we consider the 6 objectives hartman
problem.

\begin{figure}[htbp]
\centering
\includegraphics[width=\textwidth]{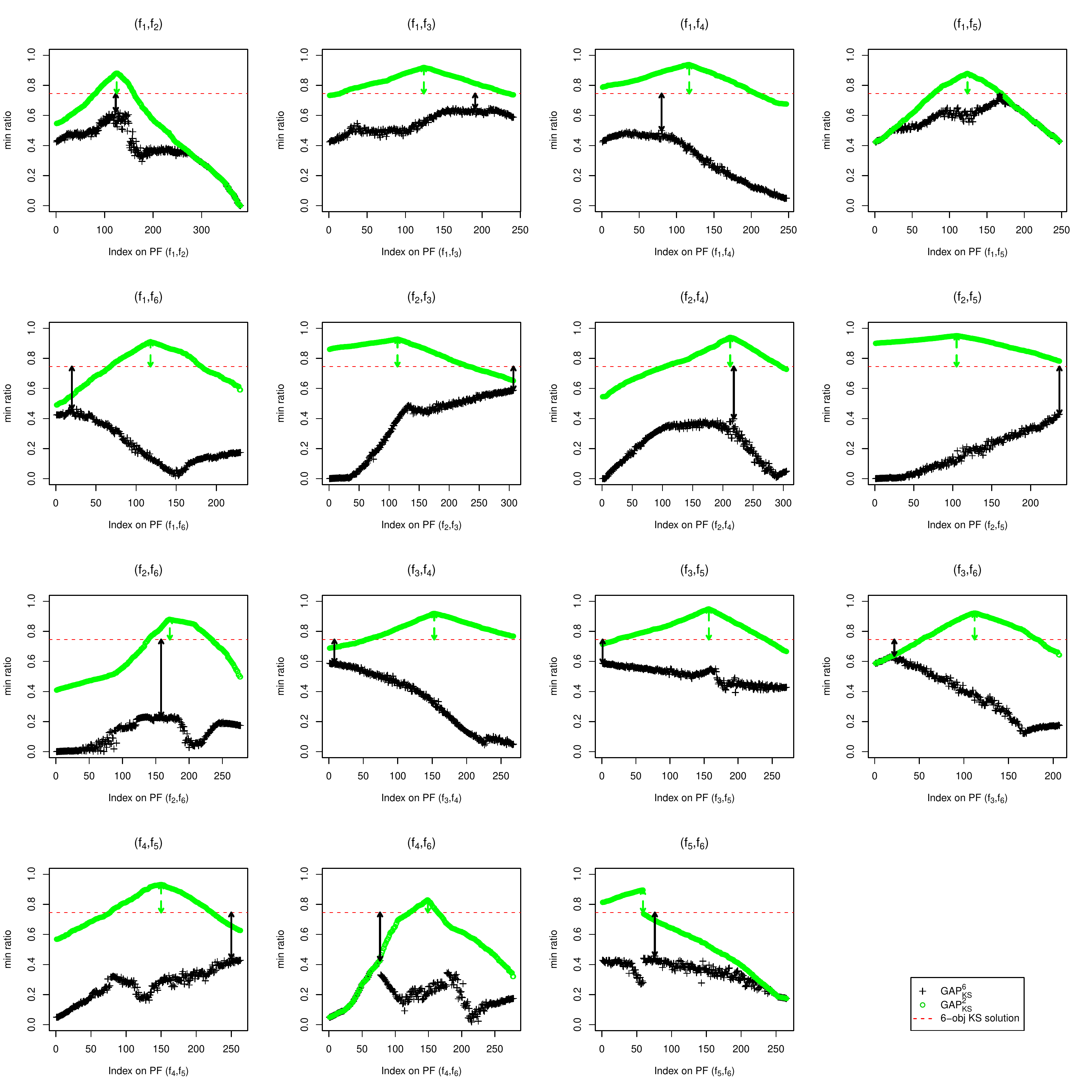}
\caption{Benefit ratios of Pareto optimal solutions for combinations of two
objectives, over the same two objectives (green $\circ$) and over all
objectives (black $+$). The corresponding performance of the 6 objectives KS
solution is depicted by the red dashed line. The black arrows (resp. dashed
green) marks the difference between the best max min ratio over 6 (resp. 2)
objectives bi-objective Pareto fronts solution and the KS.}
\label{fig:interestMany}
\end{figure}

\begin{figure}[htbp]
\centering
\includegraphics[width=\textwidth]{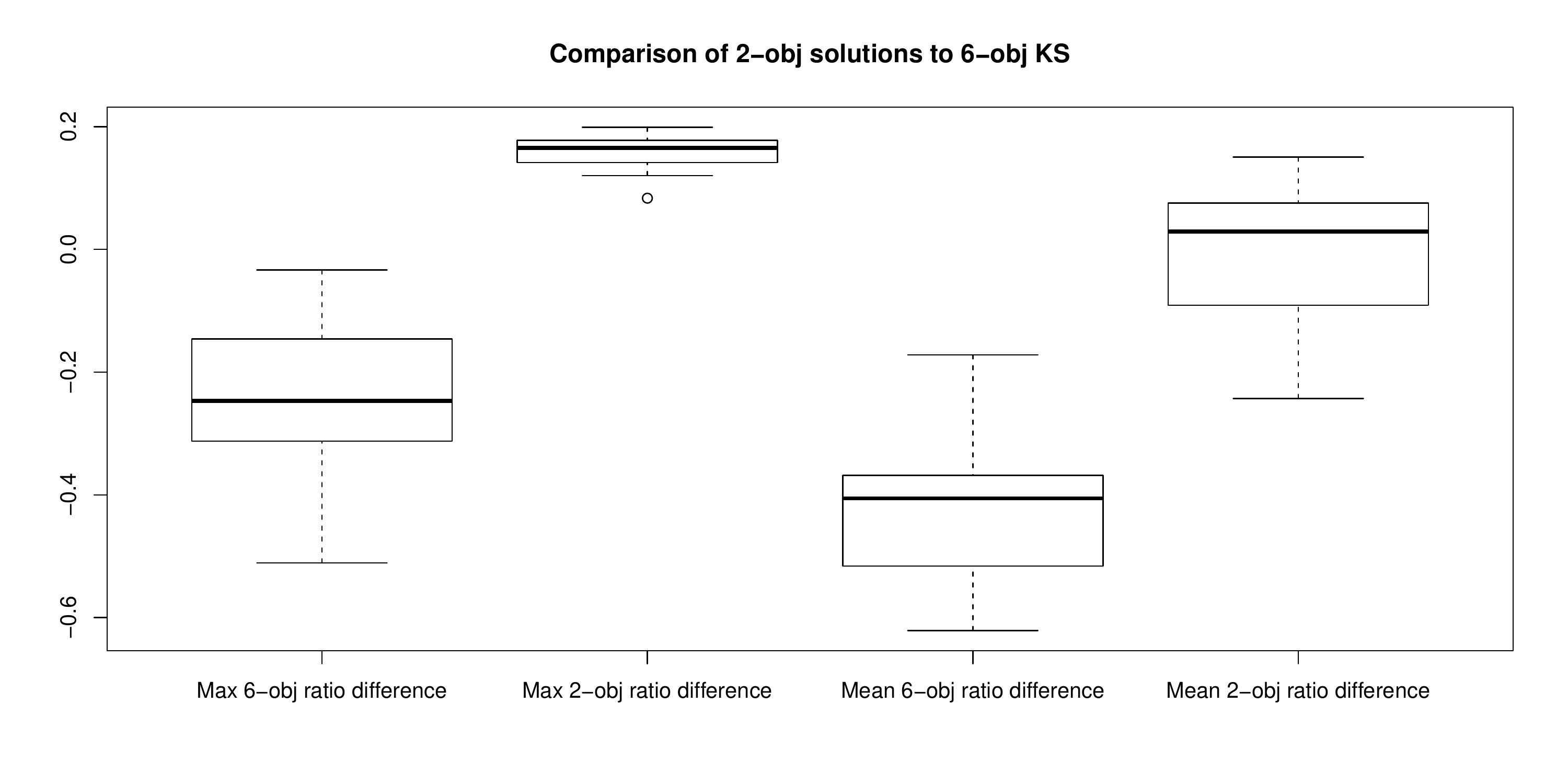}
\caption{Loss of discarding 4 objectives on their benefit ratio (black minus
KS) versus gain on the remaining two (green minus KS), either for the max
(i.e., KS on 2-objective Pareto fronts, depicted by arrows in Figure
\ref{fig:interestMany}) or the mean .}
\label{fig:boxMany}
\end{figure}

The setup is as follows: a reference KS solution is computed on this toy
example, then we estimate the Pareto front for all combinations of two
objectives. The results are given in Figure \ref{fig:interestMany} and
summarized in Figure \ref{fig:boxMany}. The best ratios over two objectives
are better than the KS ones, but on average more is lost on the remaining
objectives, advocating to keep all objectives and perhaps use preferences
rather that discarding objectives. Looking uniformly at points on those
2-objective Pareto fronts, the gain is less important while the loss on the
remaining objective is even larger.

\FloatBarrier

\section*{Appendix F: additional results on the convolutional neural network problem}
\begin{figure}[htbp]
\centering
\includegraphics[trim=5mm 5mm 5mm 5mm, width=\textwidth, clip]{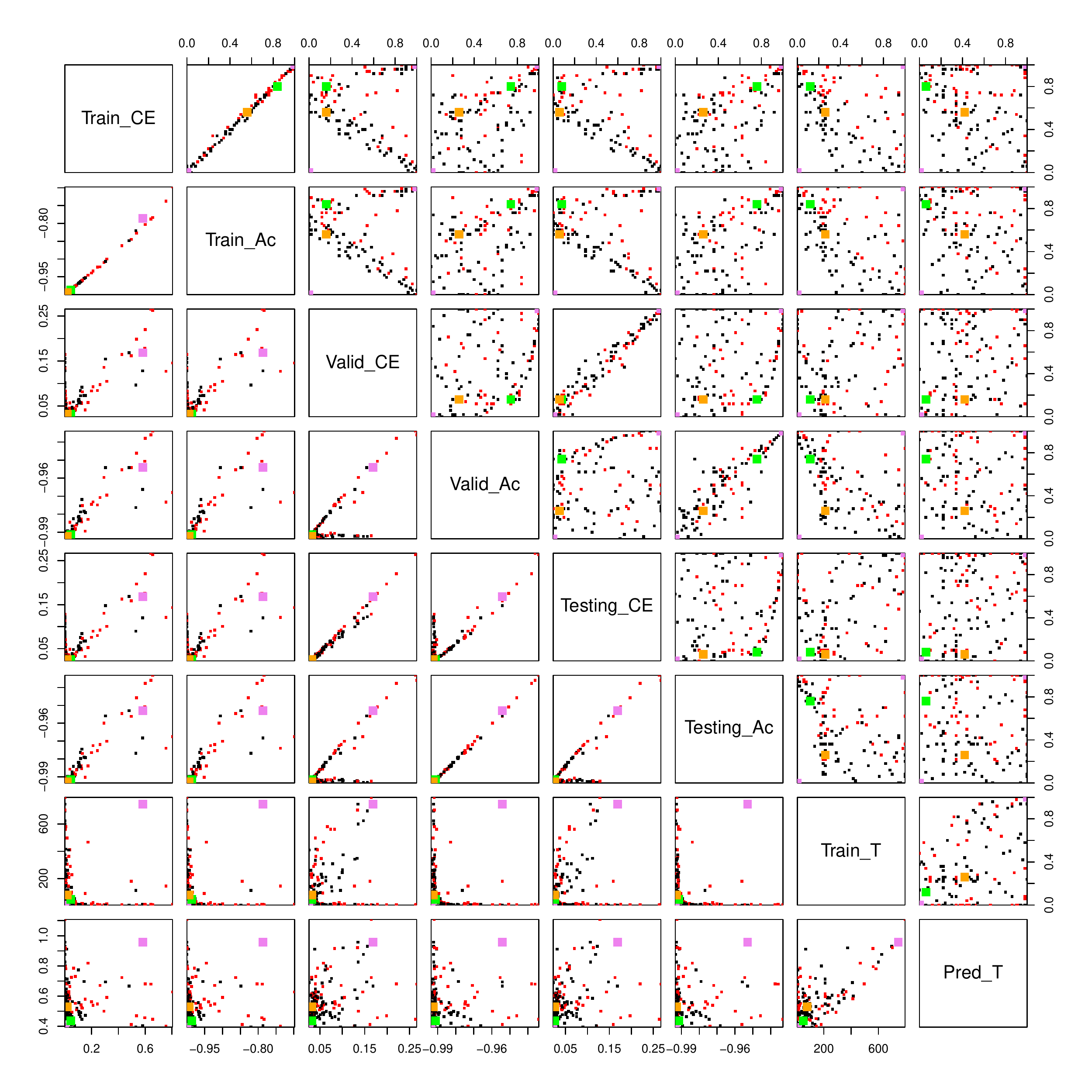}
 \caption{Projections on 2D spaces of the 150 observations, either in the
 original space (lower left triangle) or the copula space (upper right
 triangle). The colors match the ones of Figure \ref{fig:parcoordscnn}. One
 may observe in general the central position of KS in the original space and
 the central position of CKS in the copula space.}
\label{fig:allpairs_cnn}
\end{figure}

\FloatBarrier

\section*{Appendix G: additional results on the calibration problem}
\begin{figure}[htbp]
\centering
\includegraphics[trim=5mm 5mm 5mm 5mm, width=\textwidth, clip]{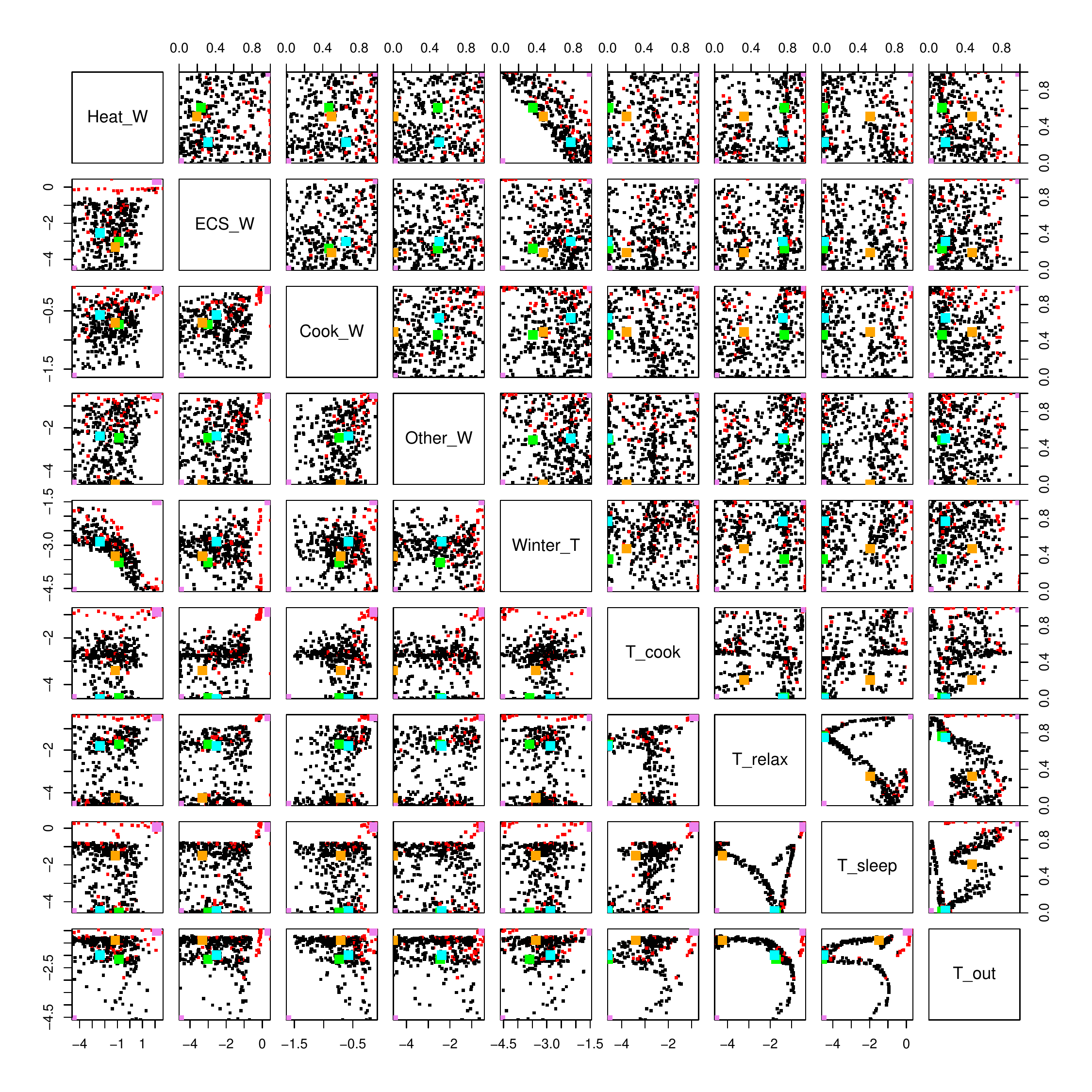}
 \caption{Projections on 2D spaces of the 400 observations, either in the
 original space (lower left triangle) or the copula space (upper right
 triangle). The colors match the ones of Figure \ref{fig:parcoordsgama}. One
 may observe in general the central position of KS in the original space and
 the central position of CKS in the copula space.}
\label{fig:allpairs_gama}
\end{figure}

\FloatBarrier

\bibliography{UQ_KS}

\begin{thebibliography}{81}
\providecommand{\natexlab}[1]{#1}
\providecommand{\url}[1]{\texttt{#1}}
\expandafter\ifx\csname urlstyle\endcsname\relax
  \providecommand{\doi}[1]{doi: #1}\else
  \providecommand{\doi}{doi: \begingroup \urlstyle{rm}\Url}\fi

\bibitem[Abou El~Majd et~al.(2010)Abou El~Majd, Desideri, and
  Habbal]{el2010optimisation}
Badr Abou El~Majd, Jean-Antoine Desideri, and Abderrahmane Habbal.
\newblock Optimisation de forme fluide-structure par un jeu de {N}ash.
\newblock \emph{Revue Africaine de la Recherche en Informatique et
  Math{\'e}matiques Appliqu{\'e}es}, 13:\penalty0 3--15, 2010.

\bibitem[Allaire and Chollet(2018)]{Allaire2018}
Joseph~J. Allaire and François Chollet.
\newblock \emph{keras: R Interface to 'Keras'}, 2018.
\newblock URL \url{https://keras.rstudio.com}.
\newblock R package version 2.1.4.

\bibitem[Asafuddoula et~al.(2015)Asafuddoula, Ray, and
  Sarker]{asafuddoula2015decomposition}
Md~Asafuddoula, Tapabrata Ray, and Ruhul Sarker.
\newblock A decomposition-based evolutionary algorithm for many objective
  optimization.
\newblock \emph{IEEE Transactions on Evolutionary Computation}, 19\penalty0
  (3):\penalty0 445--460, 2015.

\bibitem[Bader and Zitzler(2011)]{bader2011hype}
Johannes Bader and Eckart Zitzler.
\newblock Hype: An algorithm for fast hypervolume-based many-objective
  optimization.
\newblock \emph{Evolutionary {C}omputation}, 19\penalty0 (1):\penalty0 45--76,
  2011.

\bibitem[Bechikh et~al.(2010)Bechikh, Said, and Ghedira]{bechikh2010estimating}
Slim Bechikh, Lamjed~Ben Said, and Khaled Ghedira.
\newblock Estimating nadir point in multi-objective optimization using mobile
  reference points.
\newblock In \emph{IEEE Congress on Evolutionary Computation}, pages 1--9.
  IEEE, 2010.

\bibitem[Bect et~al.(2012)Bect, Ginsbourger, Li, Picheny, and
  Vazquez]{bect2012sequential}
Julien Bect, David Ginsbourger, Ling Li, Victor Picheny, and Emmanuel Vazquez.
\newblock Sequential design of computer experiments for the estimation of a
  probability of failure.
\newblock \emph{Statistics and Computing}, 22\penalty0 (3):\penalty0 773--793,
  2012.

\bibitem[Bect et~al.(2019)Bect, Bachoc, Ginsbourger,
  et~al.]{bect2019supermartingale}
Julien Bect, Fran{\c{c}}ois Bachoc, David Ginsbourger, et~al.
\newblock A supermartingale approach to gaussian process based sequential
  design of experiments.
\newblock \emph{Bernoulli}, 25\penalty0 (4A):\penalty0 2883--2919, 2019.

\bibitem[Bergstra and Bengio(2012)]{Bergstra2012}
James Bergstra and Yoshua Bengio.
\newblock Random search for hyper-parameter optimization.
\newblock \emph{The Journal of Machine Learning Research}, 13:\penalty0
  281--305, 2012.

\bibitem[Bergstra et~al.(2011)Bergstra, Bardenet, Bengio, and
  K{\'e}gl]{bergstra2011algorithms}
James~S Bergstra, R{\'e}mi Bardenet, Yoshua Bengio, and Bal{\'a}zs K{\'e}gl.
\newblock Algorithms for hyper-parameter optimization.
\newblock In \emph{Advances in {N}eural {I}nformation {P}rocessing {S}ystems},
  pages 2546--2554, 2011.

\bibitem[Binois and Picheny(2019)]{binois2019gpareto}
Micka{\"e}l Binois and Victor Picheny.
\newblock {GPareto}: An {R} package for {G}aussian-process-based
  multi-objective optimization and analysis.
\newblock \emph{Journal of Statistical Software}, 89\penalty0 (1):\penalty0
  1--30, 2019.

\bibitem[Binois et~al.(2015)Binois, Rulli\`ere, and Roustant]{Binois2015b}
Micka\"el Binois, Didier Rulli\`ere, and Olivier Roustant.
\newblock On the estimation of {P}areto fronts from the point of view of copula
  theory.
\newblock \emph{Information Sciences}, 324:\penalty0 270 -- 285, 2015.

\bibitem[Binois et~al.(2018)Binois, Gramacy, and
  Ludkovski]{binois2016practical}
Micka{\"e}l Binois, Robert~B Gramacy, and Mike Ludkovski.
\newblock Practical heteroscedastic {G}aussian process modeling for large
  simulation experiments.
\newblock \emph{Journal of Computational and Graphical Statistics}, 27\penalty0
  (4):\penalty0 808--821, 2018.

\bibitem[Binois et~al.(2019)Binois, Huang, Gramacy, and
  Ludkovski]{binois2019replication}
Micka{\"e}l Binois, Jiangeng Huang, Robert~B Gramacy, and Mike Ludkovski.
\newblock Replication or exploration? {S}equential design for stochastic
  simulation experiments.
\newblock \emph{Technometrics}, 61\penalty0 (1):\penalty0 7--23, 2019.

\bibitem[Bourgais et~al.(2017)Bourgais, Taillandier, and
  Vercouter]{bourgais2017enhancing}
Mathieu Bourgais, Patrick Taillandier, and Laurent Vercouter.
\newblock Enhancing the behavior of agents in social simulations with emotions
  and social relations.
\newblock In \emph{18th workshop on {M}ulti-{A}gent-{B}ased Simulation-MABS
  2017}, 2017.

\bibitem[Bozbay et~al.(2012)Bozbay, Dietrich, and Peters]{bozbay2012bargaining}
Irem Bozbay, Franz Dietrich, and Hans Peters.
\newblock Bargaining with endogenous disagreement: The extended
  {K}alai--{S}morodinsky solution.
\newblock \emph{Games and Economic Behavior}, 74\penalty0 (1):\penalty0
  407--417, 2012.

\bibitem[Charnes and Cooper(1977)]{charnes1977goal}
Abraham Charnes and William~Wager Cooper.
\newblock Goal programming and multiple objective optimizations: Part 1.
\newblock \emph{European journal of operational research}, 1\penalty0
  (1):\penalty0 39--54, 1977.

\bibitem[Cheng et~al.(2017)Cheng, Li, Tian, Zhang, Yang, Jin, and
  Yao]{cheng2017benchmark}
Ran Cheng, Miqing Li, Ye~Tian, Xingyi Zhang, Shengxiang Yang, Yaochu Jin, and
  Xin Yao.
\newblock A benchmark test suite for evolutionary many-objective optimization.
\newblock \emph{Complex \& Intelligent Systems}, 3\penalty0 (1):\penalty0
  67--81, 2017.

\bibitem[Chevalier et~al.(2014)Chevalier, Bect, Ginsbourger, Vazquez, Picheny,
  and Richet]{chevalier2014fast}
Cl{\'e}ment Chevalier, Julien Bect, David Ginsbourger, Emmanuel Vazquez, Victor
  Picheny, and Yann Richet.
\newblock Fast parallel kriging-based stepwise uncertainty reduction with
  application to the identification of an excursion set.
\newblock \emph{Technometrics}, 56\penalty0 (4):\penalty0 455--465, 2014.

\bibitem[Chevalier et~al.(2015)Chevalier, Emery, and
  Ginsbourger]{chevalier2015fast}
Cl{\'e}ment Chevalier, Xavier Emery, and David Ginsbourger.
\newblock Fast update of conditional simulation ensembles.
\newblock \emph{Mathematical Geosciences}, 47\penalty0 (7):\penalty0 771--789,
  2015.

\bibitem[Chollet et~al.(2015)]{chollet2015keras}
Fran\c{c}ois Chollet et~al.
\newblock Keras.
\newblock \url{https://github.com/fchollet/keras}, 2015.

\bibitem[Chugh et~al.(2017)Chugh, Sindhya, Hakanen, and Miettinen]{Chugh2017}
Tinkle Chugh, Karthik Sindhya, Jussi Hakanen, and Kaisa Miettinen.
\newblock A survey on handling computationally expensive multiobjective
  optimization problems with evolutionary algorithms.
\newblock \emph{Soft Computing}, pages 1--30, 2017.

\bibitem[Chugh et~al.(2018)Chugh, Jin, Miettinen, Hakanen, and
  Sindhya]{Chugh2016}
Tinkle Chugh, Yaochu Jin, Kaisa Miettinen, Jussi Hakanen, and Karthik Sindhya.
\newblock A surrogate-assisted reference vector guided evolutionary algorithm
  for computationally expensive many-objective optimization.
\newblock \emph{IEEE Transactions on Evolutionary Computation}, 22\penalty0
  (1):\penalty0 129--143, 2018.

\bibitem[Conley and Wilkie(1991)]{conley1991bargaining}
John~P Conley and Simon Wilkie.
\newblock The bargaining problem without convexity: {E}xtending the egalitarian
  and {K}alai-{S}morodinsky solutions.
\newblock \emph{Economics Letters}, 36\penalty0 (4):\penalty0 365--369, 1991.

\bibitem[Couckuyt et~al.(2014)Couckuyt, Deschrijver, and
  Dhaene]{couckuyt2014fast}
Ivo Couckuyt, Dirk Deschrijver, and Tom Dhaene.
\newblock Fast calculation of multiobjective probability of improvement and
  expected improvement criteria for {P}areto optimization.
\newblock \emph{Journal of Global Optimization}, 60\penalty0 (3):\penalty0
  575--594, 2014.

\bibitem[Das and Dennis(1998)]{das1998normal}
Indraneel Das and John~E Dennis.
\newblock Normal-boundary intersection: A new method for generating the
  {P}areto surface in nonlinear multicriteria optimization problems.
\newblock \emph{SIAM Journal on Optimization}, 8\penalty0 (3):\penalty0
  631--657, 1998.

\bibitem[Deb et~al.(2002)Deb, Pratap, Agarwal, and Meyarivan]{deb2002fast}
Kalyanmoy Deb, Amrit Pratap, Sameer Agarwal, and TAMT Meyarivan.
\newblock A fast and elitist multiobjective genetic algorithm: {NSGA}-{II}.
\newblock \emph{IEEE transactions on Evolutionary Computation}, 6\penalty0
  (2):\penalty0 182--197, 2002.

\bibitem[Deb et~al.(2006)Deb, Chaudhuri, and Miettinen]{deb2006towards}
Kalyanmoy Deb, Shamik Chaudhuri, and Kaisa Miettinen.
\newblock Towards estimating nadir objective vector using evolutionary
  approaches.
\newblock In \emph{Proceedings of the 8th annual conference on Genetic and
  evolutionary computation}, pages 643--650. ACM, 2006.

\bibitem[Desautels et~al.(2014)Desautels, Krause, and
  Burdick]{desautels2014parallelizing}
Thomas Desautels, Andreas Krause, and Joel~W Burdick.
\newblock Parallelizing exploration-exploitation tradeoffs in gaussian process
  bandit optimization.
\newblock \emph{The Journal of Machine Learning Research}, 15\penalty0
  (1):\penalty0 3873--3923, 2014.

\bibitem[D\'esid\'eri et~al.(2014)D\'esid\'eri, Duvigneau, and
  Habbal]{JAD13-AIAA-Special}
Jean-Antoine D\'esid\'eri, R\'egis Duvigneau, and Abderrahmane Habbal.
\newblock \emph{Computational Intelligence in Aerospace Sciences, V. M. Becerra
  and M. Vassile Eds.}, volume 244 of \emph{Progress in Astronautics and
  Aeronautics}, chapter Multi-Objective Design Optimization Using Nash Games.
\newblock AIAA, 2014.

\bibitem[Diggle and Ribeiro(2007)]{Diggle2007}
Peter Diggle and Paulo~Justiniano Ribeiro.
\newblock \emph{Model-Based Geostatistics}.
\newblock Springer, 2007.

\bibitem[Dixon and Szeg{\"o}(1978)]{dixon1978towards}
Laurence Charles~Ward Dixon and Giorgio~P Szeg{\"o}.
\newblock \emph{Towards global optimisation}, volume~2.
\newblock North-Holland Amsterdam, 1978.

\bibitem[Fedorov(1972)]{fedorov1972theory}
Valerii~Vadimovich Fedorov.
\newblock \emph{Theory of Optimal Experiments}.
\newblock Elsevier, 1972.

\bibitem[Ghose and Prasad(1989)]{ghose1989solution}
Debasish Ghose and UR~Prasad.
\newblock Solution concepts in two-person multicriteria games.
\newblock \emph{Journal of Optimization Theory and Applications}, 63\penalty0
  (2):\penalty0 167--189, 1989.

\bibitem[Hakanen and Knowles(2017)]{hakanen2017using}
Jussi Hakanen and Joshua~D Knowles.
\newblock On using decision maker preferences with {ParEGO}.
\newblock In \emph{International Conference on Evolutionary Multi-Criterion
  Optimization}, pages 282--297. Springer, 2017.

\bibitem[Hennig and Schuler(2012)]{Hennig2012}
Philipp Hennig and Christian~J Schuler.
\newblock Entropy search for information-efficient global optimization.
\newblock \emph{The Journal of Machine Learning Research}, 13:\penalty0
  1809--1837, 2012.

\bibitem[Hern{\'a}ndez-Lobato et~al.(2016{\natexlab{a}})Hern{\'a}ndez-Lobato,
  Hernandez-Lobato, Shah, and Adams]{Hernandez-Lobato2015}
Daniel Hern{\'a}ndez-Lobato, Jose Hernandez-Lobato, Amar Shah, and Ryan Adams.
\newblock Predictive entropy search for multi-objective {B}ayesian
  optimization.
\newblock In \emph{International Conference on Machine Learning}, pages
  1492--1501, 2016{\natexlab{a}}.

\bibitem[Hern{\'a}ndez-Lobato et~al.(2016{\natexlab{b}})Hern{\'a}ndez-Lobato,
  Gelbart, Adams, Hoffman, and Ghahramani]{hernandez2016general}
Jos{\'e}~Miguel Hern{\'a}ndez-Lobato, Michael~A Gelbart, Ryan~P Adams,
  Matthew~W Hoffman, and Zoubin Ghahramani.
\newblock A general framework for constrained {B}ayesian optimization using
  information-based search.
\newblock \emph{Journal of Machine Learning Research}, 17\penalty0
  (160):\penalty0 1--53, 2016{\natexlab{b}}.

\bibitem[Hougaard and Tvede(2003)]{hougaard2003nonconvex}
Jens~Leth Hougaard and Mich Tvede.
\newblock Nonconvex n-person bargaining: efficient maxmin solutions.
\newblock \emph{Economic Theory}, 21\penalty0 (1):\penalty0 81--95, 2003.

\bibitem[Ishibuchi et~al.(2008)Ishibuchi, Tsukamoto, and
  Nojima]{ishibuchi2008evolutionary}
Hisao Ishibuchi, Noritaka Tsukamoto, and Yusuke Nojima.
\newblock Evolutionary many-objective optimization: A short review.
\newblock In \emph{IEEE Congress on Evolutionary Computation, 2008. CEC 2008.},
  pages 2419--2426. IEEE, 2008.

\bibitem[Ishibuchi et~al.(2016)Ishibuchi, Setoguchi, Masuda, and
  Nojima]{ishibuchi2016performance}
Hisao Ishibuchi, Yu~Setoguchi, Hiroyuki Masuda, and Yusuke Nojima.
\newblock Performance of decomposition-based many-objective algorithms strongly
  depends on pareto front shapes.
\newblock \emph{IEEE Transactions on Evolutionary Computation}, 21\penalty0
  (2):\penalty0 169--190, 2016.

\bibitem[Jalali et~al.(2017)Jalali, Van~Nieuwenhuyse, and
  Picheny]{jalali2017comparison}
Hamed Jalali, Inneke Van~Nieuwenhuyse, and Victor Picheny.
\newblock Comparison of kriging-based algorithms for simulation optimization
  with heterogeneous noise.
\newblock \emph{European Journal of Operational Research}, 261\penalty0
  (1):\penalty0 279--301, 2017.

\bibitem[Jones et~al.(1998)Jones, Schonlau, and Welch]{jones1998efficient}
Donald~R Jones, Matthias Schonlau, and William~J Welch.
\newblock Efficient global optimization of expensive black-box functions.
\newblock \emph{Journal of Global Optimization}, 13\penalty0 (4):\penalty0
  455--492, 1998.

\bibitem[Junker(2004)]{junker2004preference}
Ulrich Junker.
\newblock Preference-based search and multi-criteria optimization.
\newblock \emph{Annals of Operations Research}, 130\penalty0 (1-4):\penalty0
  75--115, 2004.

\bibitem[Kalai and Smorodinsky(1975)]{KSE1975}
Ehud Kalai and Meir Smorodinsky.
\newblock Other solutions to {N}ash's bargaining problem.
\newblock \emph{Econometrica}, 43:\penalty0 513--518, 1975.

\bibitem[Knowles(2006)]{Knowles2006}
Joshua Knowles.
\newblock {ParEGO}: a hybrid algorithm with on-line landscape approximation for
  expensive multiobjective optimization problems.
\newblock \emph{IEEE Transactions on Evolutionary Computation}, 10\penalty0
  (1):\penalty0 50--66, February 2006.

\bibitem[Kukkonen and Lampinen(2007)]{Kukkonen2007}
Saku Kukkonen and Jouni Lampinen.
\newblock Ranking-dominance and many-objective optimization.
\newblock In \emph{IEEE Congress on Evolutionary Computation, 2007. CEC 2007.},
  pages 3983--3990. IEEE, 2007.

\bibitem[LeCun et~al.(1998)LeCun, Bottou, Bengio, and
  Haffner]{lecun1998gradient}
Yann LeCun, L{\'e}on Bottou, Yoshua Bengio, and Patrick Haffner.
\newblock Gradient-based learning applied to document recognition.
\newblock \emph{Proceedings of the IEEE}, 86\penalty0 (11):\penalty0
  2278--2324, 1998.

\bibitem[Li et~al.(2017)Li, Zhen, and Yao]{li2017read}
Miqing Li, Liangli Zhen, and Xin Yao.
\newblock How to read many-objective solution sets in parallel coordinates
  [educational forum].
\newblock \emph{IEEE Computational Intelligence Magazine}, 12\penalty0
  (4):\penalty0 88--100, 2017.

\bibitem[McKay et~al.(1979)McKay, Beckman, and Conover]{mckay1979comparison}
Michael~D McKay, Richard~J Beckman, and William~J Conover.
\newblock Comparison of three methods for selecting values of input variables
  in the analysis of output from a computer code.
\newblock \emph{Technometrics}, 21\penalty0 (2):\penalty0 239--245, 1979.

\bibitem[Miettinen(2012)]{miettinen2012nonlinear}
Kaisa Miettinen.
\newblock \emph{Nonlinear multiobjective optimization}, volume~12.
\newblock Springer Science \& Business Media, 2012.

\bibitem[Nelsen(2006)]{Nelsen2006}
Roger~B Nelsen.
\newblock \emph{An Introduction to Copulas}.
\newblock Springer, 2006.

\bibitem[Niederreiter(1988)]{niederreiter1988low}
Harald Niederreiter.
\newblock Low-discrepancy and low-dispersion sequences.
\newblock \emph{Journal of Number Theory}, 30\penalty0 (1):\penalty0 51--70,
  1988.

\bibitem[Oakley(2004)]{oakley2004estimating}
Jeremy Oakley.
\newblock Estimating percentiles of uncertain computer code outputs.
\newblock \emph{Journal of the Royal Statistical Society: Series C (Applied
  Statistics)}, 53\penalty0 (1):\penalty0 83--93, 2004.

\bibitem[Omelka et~al.(2009)Omelka, Gijbels, and Veraverbeke]{Omelka2009}
Marek Omelka, Ir{\`e}ne Gijbels, and No{\"e}l Veraverbeke.
\newblock Improved kernel estimation of copulas: weak convergence and
  goodness-of-fit testing.
\newblock \emph{The Annals of Statistics}, 37\penalty0 (5B):\penalty0
  3023--3058, 2009.

\bibitem[Parr(2013)]{parr2013improvement}
James Parr.
\newblock \emph{Improvement criteria for constraint handling and multiobjective
  optimization}.
\newblock PhD thesis, University of Southampton, 2013.

\bibitem[Picheny(2013)]{picheny2013multi}
Victor Picheny.
\newblock Multiobjective optimization using {G}aussian process emulators via
  stepwise uncertainty reduction.
\newblock \emph{Statistics and Computing}, pages 1--16, 2013.

\bibitem[Picheny and Binois(2018)]{picheny2017}
Victor Picheny and Micka{\"e}l Binois.
\newblock \emph{{GPGame}: Solving Complex Game Problems using {G}aussian
  Processes}, 2018.
\newblock URL \url{http://CRAN.R-project.org/package=GPGame}.
\newblock R package version 1.1.0.

\bibitem[Picheny et~al.(2013)Picheny, Wagner, and
  Ginsbourger]{picheny2013benchmark}
Victor Picheny, Tobias Wagner, and David Ginsbourger.
\newblock A benchmark of kriging-based infill criteria for noisy optimization.
\newblock \emph{Structural and Multidisciplinary Optimization}, 48\penalty0
  (3):\penalty0 607--626, 2013.

\bibitem[Picheny et~al.(2019)Picheny, Binois, and Habbal]{picheny2016bayesian}
Victor Picheny, Micka{\"e}l Binois, and Abderrahmane Habbal.
\newblock A {B}ayesian optimization approach to find {N}ash equilibria.
\newblock \emph{Journal of Global Optimization}, 73\penalty0 (1):\penalty0
  171--192, 2019.

\bibitem[Ponweiser et~al.(2008)Ponweiser, Wagner, Biermann, and
  Vincze]{ponweiser2008multiobjective}
Wolfgang Ponweiser, Tobias Wagner, Dirk Biermann, and Markus Vincze.
\newblock Multiobjective optimization on a limited budget of evaluations using
  model-assisted s-metric selection.
\newblock In \emph{International Conference on Parallel Problem Solving from
  Nature}, pages 784--794. Springer, 2008.

\bibitem[{R Core Team}(2018)]{R2016}
{R Core Team}.
\newblock \emph{R: A Language and Environment for Statistical Computing}.
\newblock R Foundation for Statistical Computing, Vienna, Austria, 2018.
\newblock URL \url{https://www.R-project.org/}.

\bibitem[Rasmussen and Williams(2006)]{Rasmussen2006}
Carl~E. Rasmussen and Christopher Williams.
\newblock \emph{{Gaussian Processes for Machine Learning}}.
\newblock MIT Press, 2006.
\newblock URL \url{http://www.gaussianprocess.org/gpml/}.

\bibitem[Roustant et~al.(2018)Roustant, Padonou, Deville, Cl{\'e}ment, Perrin,
  Giorla, and Wynn]{Roustant2018}
Olivier Roustant, Esperan Padonou, Yves Deville, Alo{\"\i}s Cl{\'e}ment,
  Guillaume Perrin, Jean Giorla, and Henry Wynn.
\newblock Group kernels for {G}aussian process metamodels with categorical
  inputs.
\newblock \emph{arXiv preprint arXiv:1802.02368}, 2018.

\bibitem[Russo and Van~Roy(2014)]{russo2014learning}
Daniel Russo and Benjamin Van~Roy.
\newblock Learning to optimize via information-directed sampling.
\newblock In \emph{Advances in Neural Information Processing Systems}, pages
  1583--1591, 2014.

\bibitem[Schonlau et~al.(1998)Schonlau, Welch, and Jones]{schonlau1998global}
Matthias Schonlau, William~J Welch, and Donald~R Jones.
\newblock Global versus local search in constrained optimization of computer
  models.
\newblock \emph{Lecture Notes-Monograph Series}, pages 11--25, 1998.

\bibitem[Shahriari et~al.(2016)Shahriari, Swersky, Wang, Adams, and
  de~Freitas]{Shahriari2016}
Bobak Shahriari, Kevin Swersky, Ziyu Wang, Ryan~P Adams, and Nando de~Freitas.
\newblock Taking the human out of the loop: A review of {B}ayesian
  optimization.
\newblock \emph{Proceedings of the IEEE}, 104\penalty0 (1):\penalty0 148--175,
  2016.

\bibitem[Singh et~al.(2011)Singh, Isaacs, and Ray]{singh2011pareto}
Hemant~Kumar Singh, Amitay Isaacs, and Tapabrata Ray.
\newblock A {P}areto corner search evolutionary algorithm and dimensionality
  reduction in many-objective optimization problems.
\newblock \emph{IEEE Transactions on Evolutionary Computation}, 15\penalty0
  (4):\penalty0 539--556, 2011.

\bibitem[Smithson et~al.(2016)Smithson, Yang, Gross, and
  Meyer]{smithson2016neural}
Sean~C Smithson, Guang Yang, Warren~J Gross, and Brett~H Meyer.
\newblock Neural networks designing neural networks: Multi-objective
  hyper-parameter optimization.
\newblock In \emph{2016 IEEE/ACM International Conference on Computer-Aided
  Design (ICCAD)}, pages 1--8. IEEE, 2016.

\bibitem[Srinivas et~al.(2012)Srinivas, Krause, Kakade, and
  Seeger]{srinivas2012information}
Niranjan Srinivas, Andreas Krause, Sham~M Kakade, and Matthias Seeger.
\newblock Information-theoretic regret bounds for {G}aussian process
  optimization in the bandit setting.
\newblock \emph{Information Theory, IEEE Transactions on}, 58\penalty0
  (5):\penalty0 3250--3265, 2012.

\bibitem[Svenson(2011)]{Svenson2011}
Joshua~D. Svenson.
\newblock \emph{Computer Experiments: Multiobjective Optimization and
  Sensitivity Analysis}.
\newblock PhD thesis, The Ohio State University, 2011.

\bibitem[Tabatabaei et~al.(2019)Tabatabaei, Hartikainen, Sindhya, Hakanen, and
  Miettinen]{tabatabaei2019interactive}
Mohammad Tabatabaei, Markus Hartikainen, Karthik Sindhya, Jussi Hakanen, and
  Kaisa Miettinen.
\newblock An interactive surrogate-based method for computationally expensive
  multiobjective optimisation.
\newblock \emph{Journal of the Operational Research Society}, 70\penalty0
  (6):\penalty0 898--914, 2019.

\bibitem[Taillandier et~al.(2017)Taillandier, Micolier, and Taillandier]{libim}
Franck Taillandier, Alice Micolier, and Patrick Taillandier.
\newblock Li-bim (version 1.0.0), 2017.

\bibitem[Taillandier et~al.(2018)Taillandier, Gaudou, Grignard, Huynh,
  Marilleau, Caillou, Philippon, and Drogoul]{taillandier2018building}
Patrick Taillandier, Benoit Gaudou, Arnaud Grignard, Quang-Nghi Huynh, Nicolas
  Marilleau, Philippe Caillou, Damien Philippon, and Alexis Drogoul.
\newblock Building, composing and experimenting complex spatial models with the
  gama platform.
\newblock \emph{GeoInformatica}, pages 1--24, 2018.

\bibitem[Thiele et~al.(2009)Thiele, Miettinen, Korhonen, and
  Molina]{thiele2009preference}
Lothar Thiele, Kaisa Miettinen, Pekka~J Korhonen, and Julian Molina.
\newblock A preference-based evolutionary algorithm for multi-objective
  optimization.
\newblock \emph{Evolutionary computation}, 17\penalty0 (3):\penalty0 411--436,
  2009.

\bibitem[Villemonteix et~al.(2009)Villemonteix, Vazquez, and
  Walter]{villemonteix2009informational}
Julien Villemonteix, Emmanuel Vazquez, and Eric Walter.
\newblock An informational approach to the global optimization of
  expensive-to-evaluate functions.
\newblock \emph{Journal of Global Optimization}, 44\penalty0 (4):\penalty0
  509--534, 2009.

\bibitem[Wagner et~al.(2010)Wagner, Emmerich, Deutz, and
  Ponweiser]{wagner2010expected}
Tobias Wagner, Michael Emmerich, Andr{\'e} Deutz, and Wolfgang Ponweiser.
\newblock On expected-improvement criteria for model-based multi-objective
  optimization.
\newblock In \emph{International Conference on Parallel Problem Solving from
  Nature}, pages 718--727. Springer, 2010.

\bibitem[Walter and Pronzato(1997)]{walter1997identification}
Eric Walter and Luc Pronzato.
\newblock \emph{Identification of Parametric Models from Experimental Data}.
\newblock Springer Verlag, 1997.

\bibitem[Wierzbicki(1979)]{wierzbicki1979use}
Andrzej~P Wierzbicki.
\newblock The use of reference objectives in multiobjective
  optimization-theoretical implications and practical experience.
\newblock 1979.

\bibitem[Wierzbicki(1980)]{wierzbicki1980use}
Andrzej~P Wierzbicki.
\newblock The use of reference objectives in multiobjective optimization.
\newblock In \emph{Multiple criteria decision making theory and application},
  pages 468--486. Springer, 1980.

\bibitem[Wilson and Ghahramani(2010)]{wilson2010copula}
Andrew~G Wilson and Zoubin Ghahramani.
\newblock Copula processes.
\newblock In \emph{Advances in Neural Information Processing Systems}, pages
  2460--2468, 2010.

\bibitem[Zhang et~al.(2010)Zhang, Liu, Tsang, and Virginas]{Zhang2010}
Qingfu Zhang, Wudong Liu, E.~Tsang, and B.~Virginas.
\newblock Expensive multiobjective optimization by {MOEA/D} with {G}aussian
  process model.
\newblock \emph{IEEE Transactions on Evolutionary Computation}, 14\penalty0
  (3):\penalty0 456--474, 2010.

\end{thebibliography}
\bibliographystyle{unsrt}
\end{document}